\documentclass[11pt]{article}
\usepackage[a4paper,margin=1.12in]{geometry}
\usepackage{amsmath,amssymb,amsthm,mathtools}
\usepackage{enumitem,microtype,xcolor,needspace}
\usepackage[colorlinks=true,linkcolor=blue!45!black,citecolor=blue!45!black,urlcolor=blue!55!black,hypertexnames=false]{hyperref}
\hypersetup{pdftitle={Critical one-arm probabilities for metric-graph loop soups in dimensions d >= 3},pdfauthor={Yifan Gao}}
\numberwithin{equation}{section}
\newtheorem{theorem}{Theorem}[section]
\newtheorem{proposition}[theorem]{Proposition}
\newtheorem{lemma}[theorem]{Lemma}

\theoremstyle{definition}
\newtheorem{remark}[theorem]{Remark}

\newtheorem{exploration}[theorem]{Algorithm}
\title{Critical one-arm probabilities for metric-graph\\
loop soups in dimensions \(d\ge 3\)}
\author{Yifan Gao\thanks{Institute for Theoretical Sciences, Westlake University,
\href{mailto:gaoyifan75@westlake.edu.cn}{\textcolor{black}{gaoyifan75@westlake.edu.cn}}}}
\date{}

\begin{document}
\maketitle

\begin{abstract}
We study the critical Brownian loop soup on the metric graph
\(\widetilde{\mathbb Z}^d\).
Using a new unified approach, we determine the one-arm probability up to
multiplicative constants in all dimensions \(d\ge3\).
In particular, at the upper critical dimension \(d_c=6\),
we show that the one-arm probability is comparable to
\(r^{-2}(\log r)^{\frac12}\),
and moreover, with high probability a box of radius \(r\)
contains on the order of \(\log r\) macroscopic clusters, each with
volume of order \(r^4(\log r)^{-\frac12}\).

For \(d\ne6\), we recover the one-arm estimates of Cai and Ding,
\(r^{-(d-2)/2}\) for \(d=3,4,5\) and \(r^{-2}\) for
\(d>6\)~\cite{CD24,CDhigh25},
and the three-dimensional estimate of
Drewitz, Pr\'evost and Rodriguez~\cite{DPR24}.
Our proof adapts the entropy method of Dewan and Muirhead~\cite{DM23}
to the loop-soup setting. This yields a renormalization inequality for
the one-arm probability, from which we derive the upper bounds.

Independently, Rodriguez, Shi, and Zhang~\cite{RSZ26} have recently used a closely related entropy
approach to recover the one-arm estimates for \(d>6\) and improve the upper
bound for \(d=6\). They also announce the same sharp
one-arm estimate for \(d=6\) in forthcoming work.
\end{abstract}

\tableofcontents

\section{Introduction}

Understanding the critical behavior of percolation models is a central
problem in statistical physics.
At criticality, many observables are expected to obey power laws,
and determining their critical exponents remains a challenge.
There is typically an upper critical dimension \(d_c\) associated with the model.
It is predicted that
below \(d_c\), the corresponding critical exponents depend on the dimension
and satisfy the hyperscaling
relations, while above \(d_c\), they are expected to
take their mean-field values.
Moreover, at the upper critical dimension \(d=d_c\), logarithmic
corrections to the mean-field power laws are expected to
appear~\cite{EGG78}, but rigorous results are available for
only a few models; see, e.g.,~\cite{BBS14phi,BBS15wsaw,Hut25H,Hut25III}.

In this paper, we study the percolation model formed by clusters of a Poisson
collection of Brownian loops on \(\widetilde{\mathbb Z}^d\), the metric
graph of \(\mathbb Z^d\), for \(d\ge3\).
The intensity \(\frac12\) is critical for this loop
percolation model~\cite{Lup16,CDL24}.
Moreover, at the critical intensity \(\frac12\), the loop soup can be coupled with the Gaussian
free field (GFF) so that the occupation field is half of the squared field
and the loop clusters coincide with its sign clusters~\cite{LJ11,Lup16}.
The isomorphism expresses the critical two-point function explicitly
in terms of the Green function~\cite{Lup16}.
Thus loop-soup percolation, or equivalently GFF level-set
percolation, provides a tractable model for studying critical behavior
in the presence of long-range dependence.

We focus on the critical one-arm probability \(\theta_d(r)\), that is, the
probability that the origin is connected to the boundary of a box of
radius \(r\) in the critical loop soup, or equivalently, in a GFF sign
cluster; see~\eqref{eq:one-arm} below.
The upper critical dimension of this model is \(d_c=6\), which is also
the conjectured upper critical dimension of Bernoulli
percolation~\cite{Tou74,Hut25survey}.
Werner formulated various conjectures on the
scaling limits and one-arm exponents of critical loop-soup clusters in different dimensions in~\cite{Wer21}.
These conjectures have guided subsequent research, which we briefly review below.
Away from dimension six, the predicted one-arm exponents were confirmed in~\cite{CD24,CDhigh25} by the sharp estimate
\begin{equation}\label{eq:known-one-arm-estimates}
 \theta_d(r)\asymp
 \begin{cases}
  r^{-(d-2)/2}, & d=3,4,5,\\
  r^{-2}, & d>6.
 \end{cases}
\end{equation}
The three-dimensional estimate was also obtained independently
in~\cite{DPR24}, as part of a result for a broader class of transient graphs.
Earlier works established non-optimal one-arm bounds and developed
exploration and capacity methods that underlie these sharp
estimates~\cite{LW18,DW20,DPR22,DPR23,DPR26}.
Building on these developments, we give a new unified proof of
\eqref{eq:known-one-arm-estimates} in Theorem~\ref{thm:main}, and the same approach also determines
the logarithmic correction in dimension six up to multiplicative
constants.

For \(d=6\), the polynomial exponent \(2\) was already identified by
the lower bound of~\cite[Theorem~5]{DW20} and the upper bound
of~\cite[Theorem~1.1]{CD24},
\[
 r^{-2}\lesssim \theta_6(r)\lesssim r^{-2}
       \exp\!\big((\log r)^{\frac12}\log\log r\big).
\]
It was further conjectured that
\(\theta_6(r)\asymp r^{-2}(\log r)^\delta\) for some \(\delta>0\),
without specifying the logarithmic exponent~\cite[Remark~1.5]{CD24}\footnote{Cai and Ding note that this conjecture was inspired by a discussion with Tom Hutchcroft.}.
In Theorem~\ref{thm:main}, we prove that
\begin{equation}\label{eq:six-one-arm-asymptotic}
 \theta_6(r)\asymp r^{-2}(\log r)^{\frac12}.
\end{equation}
This establishes the conjectured estimate with
\(\delta=\frac12\).
To our knowledge, the value \(\frac12\) had not previously been predicted
for this model. 
This also seems to be the first rigorous determination
of a nonzero logarithmic correction exponent for a percolation model
with upper critical dimension six.
(Update: Independently, Rodriguez, Shi, and Zhang~\cite{RSZ26} use a closely related entropy
approach to recover the one-arm estimates for \(d>6\), improve the upper bound
for \(d=6\) to order \(r^{-2}\log r\), and establish analogous bounds on a class of
more general graphs. They also announce that the sharp estimate
\eqref{eq:six-one-arm-asymptotic} will be established in forthcoming work.)
Note that the exponent \(\frac12\) is smaller than \(\frac{11}{21}\), the
corresponding exponent conjectured for critical Bernoulli percolation
in dimension six~\cite{FACS04}.
A related result for the critical infinite-volume FK--Ising model on
\(\mathbb Z^6\) bounds the one-arm probability between constant
multiples of \(r^{-2}\) and \(r^{-2}\log r\), leaving the logarithmic
exponent undetermined~\cite[Theorem~1.8, (1.19)]{EGPS25}.

Our approach is to adapt the entropy method~\cite{DM23,DM23PB}
to the loop-soup setting and derive a renormalization inequality for the one-arm probabilities,
from which we obtain the sharp upper bounds in~\eqref{eq:known-one-arm-estimates}
for $d\neq6$.
Related strategies have been used, for example, for critical Bernoulli
percolation~\cite{Hut22,EGPS25perc} and the Ising/FK--Ising
model~\cite{EGPS25}.
As expected, the case of $d=6$ is much more delicate,
since the direct bound from the triangle diagram is only of order
\(r^{-2}\log r\).
To obtain the sharp logarithmic factor, we combine the entropy
method with a refined capacity estimate for the cluster of the origin.
We discuss the complete proof strategy in detail in Section~\ref{sec:strategy}.

Using the one-arm estimate \eqref{eq:six-one-arm-asymptotic}, we also obtain a quantitative description of
macroscopic clusters in dimension six.
Theorem~\ref{thm:six-large-clusters} shows that the loop soup in a box of
radius \(r\) contains on the order of \(\log r\) distinct clusters
with diameter at least \(r/8\) and vertex volume of order
\(r^4(\log r)^{-\frac12}\), with probability at least
\(1-O((\log r)^{-1})\). This complements the conjectural picture of macroscopic clusters proposed in~\cite{Wer21}.
For comparison, we review related results in other dimensions after
Theorem~\ref{thm:six-large-clusters}.

The literature on loop soups is extensive, and we can only mention a few related works here.
Recent developments on switching identities and the intensity-doubling
conjecture provide new descriptions of loop-soup 
connectivity~\cite{Wer25switch,LW25double}.
Quantitative aspects of cluster geometry have also been
investigated in~\cite{Wer21,GN24,CDQM25,CDIIC25,DPRvol24,GJ24,CDcut25,CDgap25}.
For related developments on loop soups and GFF
percolation in the plane, see~\cite{LW04,SW12,DL18,DW20,DWW22,JLQ25cross,BGNQ25,RZ25anomalous,BGL26}.
Other related models include discrete loop-soup percolation on
\(\mathbb Z^d\)~\cite{LJL13,CS16,Cha17,BDKM26}, GFF level-set
percolation~\cite{BLM87,RS13,DPR18,DGRS23,GRS22}, and the continuum
Brownian loop soup in \(\mathbb R^3\)~\cite{JL26,Gao26}.
The critical behavior of these related models nevertheless remains less
well understood in dimensions \(d\ge3\).

\paragraph{Notation.} The metric graph
\(\widetilde{\mathbb Z}^d\) is obtained by replacing each nearest-neighbor edge \(\{x,y\}\) of
\(\mathbb Z^d\) with an interval of length \(d\), whose endpoints are identified with \(x\) and \(y\).
Each edge has conductance \(c_{xy}=1/(2d)\).
We extend the lattice coordinates linearly along each edge and use
a tilde to distinguish metric-graph objects from their discrete
counterparts. For \(z\in\widetilde{\mathbb Z}^d\), write \(|z|=\max_{1\le i\le d}|z_i|\) for the
\(\ell^\infty\)-norm of \(z\).
For a nonempty set \(A\subset\widetilde{\mathbb Z}^d\), let
\(\operatorname{diam}(A)=\sup_{x,y\in A}|x-y|\) denote its diameter.
For \(A\subset\mathbb Z^d\), write \(|A|\) for its cardinality.
For \(r>0\), define the closed metric-graph box and its boundary by
\(\widetilde B_r=\{z\in\widetilde{\mathbb Z}^d:|z|\le r\}\) and
\(\partial\widetilde B_r=\{z\in\widetilde{\mathbb Z}^d:|z|=r\}\).
The vertex box is \(B_r=\widetilde B_r\cap\mathbb Z^d\), and
\(\partial B_r\) denotes its inner boundary. We write
\(B(x,r)=x+B_r\) and \(\widetilde B(x,r)=x+\widetilde B_r\).
For nonnegative quantities, \(f\lesssim g\) and \(f=O(g)\) mean
\(f\le Cg\) for some constant \(C>0\). We write \(f\gtrsim g\)
for \(g\lesssim f\), and \(f\asymp g\) when both inequalities hold.

\subsection{Main results}\label{sec:main-results}

Let \(\mathcal L\) be the Brownian loop soup on
\(\widetilde{\mathbb Z}^d\) at intensity \(\frac12\), and let \(\phi\)
be the metric-graph Gaussian free field coupled with \(\mathcal L\)
as in Lemma~\ref{thm:lupu-isomorphism}.
For \(x\in\mathbb Z^d\), let \(\widetilde{\mathcal C}(x)\) be the cluster of \(\mathcal L\) that contains \(x\).
Write \(\mathcal C(x)=\widetilde{\mathcal C}(x)\cap\mathbb Z^d\) for its discrete set.
The one-arm probability is defined by
\begin{equation}\label{eq:one-arm}
 \theta_d(r):= \mathbb P(0\xleftrightarrow{\mathcal L} \partial\widetilde B_r)
 = \mathbb P(\widetilde{\mathcal C}(0)\cap \partial\widetilde B_r\neq \emptyset)
 =2\mathbb P\!\left(
 0\xleftrightarrow{\{\phi\ge0\}} \partial\widetilde B_r\right)
\end{equation}
for all $r\ge1$.
We write \(A\xleftrightarrow{C}B\) when a path in \(C\) connects
\(A\) to \(B\). When \(C\) is a loop collection, paths are taken in the trace of
the union of finite chains of loops.
The last equality in~\eqref{eq:one-arm} follows from the coupling
in Lemma~\ref{thm:lupu-isomorphism} and the symmetry of the GFF.

\Needspace{18\baselineskip}
Our first main result determines the one-arm probability up to multiplicative constants in every dimension \(d\ge3\) using a 
new unified proof strategy.

\begin{theorem}\label{thm:main}
For every \(d\ge3\) with \(d\ne6\),
\eqref{eq:known-one-arm-estimates} holds for all \(r\ge1\).
For \(d=6\), there is \(r_0<\infty\) such that
\begin{equation}\label{eq:main-six}
 \theta_6(r)\asymp r^{-2}(\log r)^{\frac12},\qquad r\ge r_0.
\end{equation}
\end{theorem}

As mentioned, the lower bounds
\(\theta_d(r)\gtrsim r^{-(d-2)/2}\) for
\(3\le d\le6\) have been established previously in~\cite[Theorem~5]{DW20},
and the complete up-to-constants estimates for any $d\neq 6$ were later obtained in~\cite{CD24,CDhigh25}.
The $d=3$ case was obtained independently in~\cite{DPR24}, which also considered general transient graphs. 
The estimate \eqref{eq:main-six} is new, addressing the remaining critical-dimensional case $d=6$.

Theorem~\ref{thm:main} also yields information about the geometry of
macroscopic clusters in dimension six.
We fix some notation before stating the result.
For an integer \(r\ge1\), consider the critical loop soup
\(\mathcal L_{\widetilde B_r}\) in \(\widetilde B_r\), for \(d=6\).
For \(x\in B_r\), let \(\widetilde{\mathcal C}^r(x)\) be the
cluster containing $x$ in \(\mathcal L_{\widetilde B_r}\), and write
\(\mathcal C^r(x)=\widetilde{\mathcal C}^r(x)\cap\mathbb Z^6\).
Let \(\mathfrak C^r=\{\mathcal C^r(x):x\in B_r\}\) be the
family of distinct vertex clusters, which partitions \(B_r\).

\begin{theorem}
\label{thm:six-large-clusters}
There are constants \(a,c,C>0\), \(A\ge1\), and
\(r_0<\infty\) such that, for every integer \(r\ge r_0\),
if we set the thresholds \(v_r^-=ar^4(\log r)^{-\frac12}\) and
\(v_r^+=Ar^4(\log r)^{-\frac12}\), then the cluster count
\(N_r=\#\{\mathcal C\in\mathfrak C^r:
  v_r^-\le |\mathcal C|\le v_r^+,\,
  \operatorname{diam}(\mathcal C)\ge r/8\}\) satisfies
\begin{equation}\label{eq:six-large-cluster-count}
 \mathbb P_{\widetilde B_r}\!\left(c\log r\le N_r\le C\log r\right)
 \ge1-\frac{C}{\log r}.
\end{equation}
\end{theorem}

We recall some related results in dimensions \(d\ne6\).
For \(d>6\),
there are at least \(cr^{d-6}/(\log r)^2\) macroscopic clusters,
and every cluster has volume at most \(Cr^4\log r\), with high probability~\cite[Section~5.2]{Wer21}.
For \(d=3,4,5\), the number of macroscopic crossing clusters is
bounded in probability~\cite[proof of Lemma~7.1]{CDQM25}.
Conditional on a one-arm event to radius \(Cr\) for some large constant \(C\),
the volume \(|\mathcal C(0)\cap B_r|\) has order \(r^{(d+2)/2}\) for
\(d=3,4,5\) and order \(r^4\) for
\(d>6\)~\cite[Theorem~1.2]{CDIIC25}.
In concurrent work~\cite[Theorem~1.1,
Remark~5.2(5)]{DPRvol24}, for \(d=3,4,5\), the largest cluster has volume
of order \(r^{(d+2)/2}\), and the
expected number of macroscopic clusters with volumes in a fixed window
of this order is bounded.
For \(d>6\), the same work gives lower and upper bounds for the
largest cluster volume on the scales \(r^4\) and \(r^4\log r\),
respectively, and at least \(cr^{d-6}/\log r\) clusters of volume
at least \(cr^4\) with probability tending to
one~\cite[(1.6), Appendix~A]{DPRvol24}.
We also refer to~\cite{Aiz97} for corresponding results in Bernoulli percolation.

Theorem~\ref{thm:six-large-clusters} completes this picture at the upper
critical dimension \(d_c=6\) and identifies the predicted logarithmic corrections precisely, that is,
with high probability, there are on the
order of \(\log r\) macroscopic clusters in \(\widetilde B_r\), each with volume of order
\(r^4(\log r)^{-\frac12}\).

Moreover, the up-to-constants estimates for the (heterochromatic) two-arm probability have been proved for
\(d\ge3\) with \(d\ne6\)~\cite{CDtwo25},
and the scaling limits of critical clusters have recently been established
for \(d=3,4,5\)~\cite{CDscale26}.
The results of this paper can also be used to deduce corresponding
results for \(d=6\), and we plan to pursue them in forthcoming works. 
More precisely, the two-arm probability for \(d=6\) should be of order \(r^{-4}(\log r)^{\frac12}\),
and the scaling limit (under the normalization provided by Theorem~\ref{thm:six-large-clusters}) 
should be described by the integrated super-Brownian excursion, as predicted 
by Werner for $d>6$ in~\cite{Wer21,Wer25switch}.
We refer to~\cite{CCHS26,BH26,CCFK26} for recent progress on multipoint
connectivity and scaling limits in high-dimensional critical Bernoulli percolation.

\subsection{Proof strategy}\label{sec:strategy}

For completeness, we also give an independent proof of the one-arm lower
bounds in Lemma~\ref{lem:one-arm-lower}, which are sharp for \(d\ne6\).
We now focus on the more involved one-arm upper bounds. 
The proof centers on the renormalization
inequality~\eqref{eq:strategy-common-recursion}.
We first explain its derivation, then describe the geometric estimates
used to close the argument in each dimension.
Following the idea of~\cite{CD24}, we work with the weighted one-arm probability
\begin{equation}\label{eq:wr}
 w_d(r):=\mathbb E[\sqrt{2\widehat{\mathcal L}_0}\mathbf{1}_{\smash[t]{\{0\xleftrightarrow{\mathcal L}\partial\widetilde B_r\}}}]
 =\mathbb E[|\phi_0|\mathbf{1}_{\smash[t]{\{0\xleftrightarrow{\mathcal L}\partial\widetilde B_r\}}}],
\end{equation}
where \(\widehat{\mathcal L}_0\) is the total occupation time of the
loop soup at the origin, and the second equality follows from the isomorphism of
Lemma~\ref{thm:lupu-isomorphism}.
For every \(d\ge3\), both \(|\phi_0|\) and the one-arm event
\(\{0\xleftrightarrow{\mathcal L} \partial\widetilde B_r\}\) are increasing in
the loop soup. By the FKG inequality~\cite[Theorem~20.4]{LP18}
and~\eqref{eq:phi0}, we have
\begin{equation}\label{eq:weighted-to-one-arm}
 w_d(r)\ge\mathbb E|\phi_0|\,
       \mathbb P(0\xleftrightarrow{\mathcal L}\partial\widetilde B_r)
 =\sqrt{2G(0,0)/\pi}\,\theta_d(r).
\end{equation}
Thus it suffices to prove the desired upper bounds for \(w_d(r)\).
As shown in~\cite{CD24}, the additional weight \(|\phi_0|\) in \eqref{eq:wr} gives
\(w_d(r)\) a harmonic representation,
which is useful for us to apply the induction; see Lemma~\ref{lem:separator}.

\Needspace{4\baselineskip}
\medskip\noindent\textbf{The renormalization inequality.}
Set \(\beta=(d-2)/2\) and fix \(0<\delta<\frac18\). Let
\(\mathcal L^{\le\delta r}\) be \(\mathcal L\) restricted to loops of diameter
at most \(\delta r\), and let \(\widehat{\mathcal L}^{\le\delta r}_x\)
be its occupation field at $x$. Set
\begin{equation}\label{eq:wdr}
 w_d^{\le\delta r}(r)=\mathbb E[\sqrt{2\widehat{\mathcal L}^{\le\delta r}_0}\mathbf{1}_{\smash[t]{\{0\xleftrightarrow{\mathcal L^{\le\delta r}}\partial\widetilde B_r\}}}].
\end{equation}
Following the idea of~\cite{Gao26}, we first remove large loops by showing that
\begin{equation}\label{eq:large-loop-decomposition}
 w_d(r)\le w_d^{\le\delta r}(r)+C_\delta r^{-\beta}.
\end{equation}
The second term in~\eqref{eq:large-loop-decomposition} is the error
from removing loops of diameter greater than \(\delta r\), 
which can be controlled by \cite[Proposition~5]{LW18}. 
Thus the correct one-arm exponent \(\beta=(d-2)/2\) for \(d<6\)
already appears in the above decomposition~\eqref{eq:large-loop-decomposition}. 

To bound \(w_d^{\le\delta r}(r)\), we compare it with a loop soup
subject to additional killing outside a macroscopic neighborhood of the origin.
We use the killing rate
\begin{equation}\label{eq:protected-killing-rate}
 \boldsymbol{k}(x)=\kappa^2r^{-2}
 \mathbf{1}_{\smash[t]{\{|x|\ge r/8\}}},
 \qquad r\ge1,\quad\kappa\ge1.
\end{equation}
Write \(w_d^{\boldsymbol{k}}(r)\) for the corresponding weighted
one-arm probability in the loop soup with this killing (see \eqref{eq:wdk}).
This killing suppresses connections across the outer annulus and
gives the contraction estimate below,
\begin{equation}\label{eq:contraction1}
       w_d^{\boldsymbol{k}}(r)\lesssim e^{-c\kappa}w_d(r/4).
\end{equation}
At the same time, the box \(B_{r/8}\) is left free of killing.
Since \(\delta<1/8\), every loop of diameter at most \(\delta r\)
that visits the origin stays in this protected region.
Its occupation at the origin is therefore unchanged by the killing,
which is essential for the following entropy estimate. 

We compare \(w_d^{\le\delta r}(r)\) with
\(w_d^{\boldsymbol{k}}(r)\) by adapting the relative entropy method
of~\cite{DM23,DM23PB} to the loop-soup exploration.
By Proposition~\ref{lem:query},
the resulting entropy cost \(\Delta(r)\) satisfies
\begin{equation}\label{eq:exploration-entropy-1}
 \begin{aligned}
 \Delta(r)
 \lesssim \kappa^4r^{-4}
 \int_{\substack{\operatorname{diam}\gamma\le\delta r\\\boldsymbol{J}(\gamma)>0}}
 \mathcal V(\gamma)^2R(\gamma)\,
                         \mu_{\widetilde B_{2r}}^{\mathrm{loop}}(d\gamma).
 \end{aligned}
\end{equation}
Here \(\ell_x(\gamma)\) is the local time of \(\gamma\) at
\(x\in\mathbb Z^d\),
\(\mathcal V(\gamma):=\sum_{x\in\mathbb Z^d}\ell_x(\gamma)\)
is its total vertex local time, and
\(\boldsymbol{J}(\gamma):=\sum_{x\in\mathbb Z^d}\boldsymbol{k}(x)\ell_x(\gamma)\)
is its accumulated killing cost.
The factor \(R(\gamma)\) is the weighted probability
that the exploration reveals \(\gamma\), which is controlled by the weighted probability that the origin is
connected to \(\gamma\); see \eqref{eq:Rgamma} below.
As in the OSSS inequality~\cite{OSSS05}, revealment enters because
of the randomized
exploration algorithm.
By a variant of Pinsker's inequality~\cite[Lemma~2.3]{DM23}, we obtain
\begin{equation}\label{eq:weighted-comparison}
 w_d^{\le\delta r}(r)\le2w_d^{\boldsymbol{k}}(r)+2\Delta(r).
\end{equation}
Combining \eqref{eq:large-loop-decomposition}, \eqref{eq:contraction1}, and \eqref{eq:weighted-comparison}, we get
the desired renormalization inequality
\begin{equation}\label{eq:strategy-common-recursion}
 w_d(r)\le C_\delta r^{-\beta}
       +Ce^{-c\kappa}w_d(r/4)+C\Delta(r).
\end{equation}
The three terms account, respectively, for large-loop removal,
the contraction produced by killing, and the entropy cost of the
comparison.
The remaining task is to bound \(\Delta(r)\) sharply enough to
iterate~\eqref{eq:strategy-common-recursion}.

\medskip\noindent\textbf{Dimensions three to five.}
Consider a loop \(\gamma\) contributing to
\eqref{eq:exploration-entropy-1}, with
\(s<\operatorname{diam}\gamma\le2s\) and \(s\le\delta r\).
Since \(\boldsymbol J(\gamma)>0\), the loop visits a vertex outside
\(B_{r/8}\). Its weighted revealment probability satisfies
\begin{equation}\label{eq:Rgamma}
       R(\gamma)\lesssim \mathbb E[|\phi_0|\mathbf{1}_{\smash[t]{\{0\xleftrightarrow{\mathcal L}\gamma\}}}]
       \lesssim \Big(\frac{s}{r}\Big)^{d-2}w_d(2s).
\end{equation}
The integral of \(\mathcal V(\gamma)^2\) over loops in this diameter
range is bounded by a constant multiple of \(r^ds^{4-d}\).
Combining the two estimates, multiplying by \(\kappa^4r^{-4}\), and
summing over dyadic scales, we obtain
\begin{equation}\label{eq:deltar}
        \Delta(r)\lesssim \kappa^4r^{-2}
       \sum_{j=0}^{\lfloor\log_2(\delta r)\rfloor}
              2^{2j}w_d(2^{j+1}).
\end{equation}
Substituting this bound into \eqref{eq:strategy-common-recursion}, we obtain the following renormalization inequality
\begin{equation}
 w_d(r)\le C_\delta r^{-\beta}
       +Ce^{-c\kappa}w_d(r/4)
       +C\kappa^4r^{-2}
       \sum_{j=0}^{\lfloor\log_2(\delta r)\rfloor}
              2^{2j}w_d(2^{j+1}).
\end{equation}
Since \(2>\beta\), we choose \(\kappa\) large and then
\(\delta\) small and prove \(w_d(r)\lesssim r^{-\beta}\) by induction
over scales.
In fact, there is a shorter way to get the proof for \(d=3\); see Section~\ref{subsec:three-dimensional}.

\medskip\noindent\textbf{The triangle diagram above dimension six.}
Above dimension six, we control the entropy cost using the
finite-volume triangle diagram ($G$ is the Green function)
\[
 \nabla_{B_{2r}}
 :=\sup_{x\in B_{2r}}\sum_{y,z\in B_{2r}}
       G_{B_{2r}}(x,y)G_{B_{2r}}(y,z)G_{B_{2r}}(z,x).
\]
By the union bound and the weighted two-point function, 
we bound the weighted revealment probability by 
\[
R(\gamma)\lesssim \sum_{x\in B_{2r}} \mathbf{1}_{\smash[t]{\{x\in\gamma\}}} \mathbb E[|\phi_0|
\mathbf{1}_{\smash[t]{\{0\xleftrightarrow{\mathcal L}x\}}}] \lesssim 
\sum_{x\in B_{2r}} \mathbf{1}_{\smash[t]{\{x\in\gamma\}}}   G(0,x).
\]
Moreover, uniformly in \(x\in B_{2r}\), we have 
\[
 \int\mathcal V(\gamma)^2\mathbf{1}_{\smash[t]{\{x\in\gamma\}}}
           \,\mu_{\widetilde B_{2r}}^{\mathrm{loop}}(d\gamma)
 \le 2\nabla_{B_{2r}}.
\]
Combining these bounds with
\(\sum_{x\in B_{2r}}G(0,x)\lesssim r^2\), we obtain
\begin{equation}\label{eq:strategy-triangle-entropy}
 \Delta(r)\lesssim \kappa^4r^{-4} \nabla_{B_{2r}} \sum_{x\in B_{2r}}  G(0,x) 
 \lesssim \kappa^4r^{-2}\nabla_{B_{2r}}.
\end{equation}
For \(d>6\), the triangle diagram is bounded uniformly in \(r\).
Since \(\beta>2\), both error terms in
\eqref{eq:strategy-common-recursion} are \(O(r^{-2})\) for fixed
parameters. Choosing \(\kappa\) sufficiently large and iterating,
we obtain \(w_d(r)\lesssim r^{-2}\); see
Section~\ref{sec:above-six}.
However, in dimension six, by the estimate
\(\nabla_{B_{2r}}\lesssim\log r\), we obtain only
\(w_6(r)\lesssim r^{-2}\log r\).
We next explain how to improve the logarithmic factor.

\medskip\noindent\textbf{One-arm estimate in dimension six.}
Dimension six is the borderline case for both preceding estimates.
To improve the triangle-diagram bound, we use finer information
about the geometry of \(\mathcal C(0)\) to estimate the entropy
cost directly.
Write \(\mathcal C_r(0):=\mathcal C(0)\cap B_r\).
The key input is the following bound on the weighted expected capacity
of \(\mathcal C_r(0)\),
\begin{equation}\label{eq:Mrr}
     M(r)=\mathbb E\!\left[
  |\phi_0|\operatorname{cap}(\mathcal C_r(0))\right] \lesssim r^2 (\log r)^{-1/2}.   
\end{equation}
By Lemma~\ref{lem:capacity-arm-transfer}, we have
\(\theta_6(r)\gtrsim M(r)^{-1}\). Together with~\eqref{eq:Mrr}, this proves the
lower bound in Theorem~\ref{thm:main}.
Before turning to the upper bound, we explain how to obtain this
capacity estimate.

Expanding $\operatorname{cap}(\mathcal C_r(0))$ as a sum of equilibrium masses reduces the
problem to estimating
\(\mathbb E[|\phi_0|e_{\mathcal C_r(0)}(x);
0\xleftrightarrow{\mathcal L}x]\) for \(x\in B_r\).
Using the BFS--Dynkin isomorphism, we represent this weighted
expectation by inserting an independent bridge between \(0\) and \(x\),
\begin{equation}\label{eq:BFS}
       \mathbb E\!\big[|\phi_0|e_{\mathcal C_r(0)}(x); 0\xleftrightarrow{\mathcal L}x\big] 
=\int\mu_{\widetilde{\mathbb Z}^6}^{x,0}(d\gamma)\,
       \mathbb E\!\left[
       e_{K^\gamma\cap B_r}(x)(2\Gamma^\gamma(x))^{-\frac12}\right].
\end{equation}
Here \(\mu_{\widetilde{\mathbb Z}^6}^{x,0}\) is the path
measure from \(x\) to \(0\), \(\Gamma^\gamma\) is the occupation
field after adding the independent bridge \(\gamma\) to the loop soup,
and \(K^\gamma\) is the discrete set of the resulting cluster
containing \(x\); see~\eqref{eq:six-local-bridge-insertion}.
This representation leads us to consider the weighted escape probability
\[
 b(R):=\mathbb E\!\left[
       e_{\mathcal K_R}(0)(2\Gamma_R(0))^{-\frac12}\right],
\]
where \(\Gamma_R\) and \(\mathcal K_R\) are defined similarly by adding
an independent metric-graph Brownian path, started at \(0\) and
stopped on reaching \(\partial\widetilde B_R\), to the loop soup
\(\mathcal L_{\widetilde B_R}\).
The Markov property and monotonicity bound the right-hand side
of~\eqref{eq:BFS} by \(G(0,x)b(|x|/8)\).
Thus~\eqref{eq:Mrr} follows by summing over \(x\) and applying
the following weighted escape estimate,
\begin{equation}\label{eq:br}
        b(R)\lesssim (\log R)^{-\frac12}.
\end{equation}
We prove~\eqref{eq:br} by a multiscale renormalization argument
in Section~\ref{sec:escape-proof}, and refer to that section for details.

We next explain how to use \eqref{eq:br} to derive the one-arm
upper bound on \(\theta_6(r)\).
To refine the entropy cost, this time we condition on the cluster in
the outer annulus first, denoted by
\(\widehat{\mathcal C}_r(0)=\mathcal C(0)\cap
(B_{2r}\setminus B_{r/16})\), and then estimate the total mass of
revealed loops that intersect $\widehat{\mathcal C}_r(0)$.
Summing according to the duration of loops, we get
\begin{equation}\label{eq:V1}
       \int \mathcal V(\gamma)^2
       \mathbf 1_{\{\gamma\cap \widehat{\mathcal C}_r(0)\ne\varnothing\}}
       \,\mu_{\widetilde B_{2r}}^{\mathrm{loop}}(d\gamma)
       \lesssim\sum_{j=0}^{\lceil\log_2r\rceil}
              \sum_{Q\in\mathcal Q_{2^j}}\operatorname{cap}(\widehat{\mathcal C}_r(0)\cap Q),
\end{equation}
where \(\mathcal Q_{2^j}\) is a partition of \(\mathbb Z^6\) into boxes
of side length \(2^j\).
Using the estimate~\eqref{eq:br}, we obtain
\begin{equation}\label{eq:V2}
       \mathbb E\!\Big[|\phi_0|
       \sum_{Q\in\mathcal Q_{2^j}}
       \operatorname{cap}(\widehat{\mathcal C}_r(0)\cap Q)\Big] \lesssim r^2 (j+1)^{-\frac12}.
\end{equation}
Combining \eqref{eq:exploration-entropy-1}, \eqref{eq:V1}, and \eqref{eq:V2}, we conclude that 
\[
 \Delta(r)
 \lesssim \kappa^4r^{-4} \sum_{j=0}^{\lceil\log_2r\rceil} r^2 (j+1)^{-\frac12}\lesssim \kappa^4r^{-2} (\log r)^{\frac12}.
\]
Substituting this improved entropy bound into
\eqref{eq:strategy-common-recursion} and using induction, we prove the sharp
one-arm upper bound.

\medskip\noindent\textbf{Organization.}
Section~\ref{sec:inputs} collects preliminary estimates and proves the one-arm
lower bound for \(d\ne6\).
Section~\ref{sec:exploration} proves the one-arm upper bounds for \(d\ne6\),
starting with a direct three-dimensional proof and then giving a
unified proof by a refined entropy bound.
Section~\ref{sec:capacity} proves the sharp one-arm bounds in dimension six.
Section~\ref{sec:escape-proof} establishes the weighted escape estimate~\eqref{eq:br}.
Section~\ref{sec:large-clusters} proves the macroscopic-cluster count in
Theorem~\ref{thm:six-large-clusters}.

\medskip\noindent\textbf{Acknowledgments.}
We are grateful to Zhenhao Cai, Jian Ding, and Xinyi Li for inspiring discussions
and helpful comments.
AI tools assisted with the exploration and refinement of some technical arguments
in the proof of Proposition~\ref{s6:escape-bound}, suggested using the BFS--Dynkin
isomorphism \eqref{eq:bridge} in place of the switching identity~\cite{Wer25switch} for simplicity, and helped with manuscript
preparation and language editing.
The main ideas of the proof, including the ``small/large'' loop-decomposition strategy inspired by
the author's recent work~\cite{Gao26} and the adaptation of the entropy method
from~\cite{DM23}, were developed by the author.
The author takes full responsibility for the content of this paper.

\medskip

We briefly record the chronology of the present paper in relation to the
independent work of Rodriguez, Shi, and Zhang~\cite{RSZ26}.
The results of~\cite{Gao26}, which primarily concern subcritical one-arm
probabilities for loop soups in dimensions \(d\ge3\), were obtained in late April 2026.
As stated in the paragraph immediately following~\cite[Theorem~1.10]{Gao26},
``In this paper, we focus on the loop soup in a unified framework
without using the GFF, and leave the critical case to future work.''
The present paper continues that work and addresses the critical case.
Its proof ideas were developed in late July 2026, and the details, including the
critical case \(d=6\), were written down in August 2026.
A draft was first shared with a small group of colleagues on September 24, 2026.
The author was unaware of~\cite{RSZ26} until its first version appeared on arXiv
on November 3, 2026.
Finally, we believe that the methods developed here can also be applied to more general graphs.

\section{Preliminaries}\label{sec:inputs}
\subsection{Loop soup and Gaussian free field}\label{subsec:basic}

Let \(Z\) be the Brownian motion on \(\widetilde{\mathbb Z}^d\),
as constructed in~\cite[Section~2]{Lup16}.
Its local times \(L_t^z\) are occupation densities relative to
the length measure \(\mathfrak m\) on \(\widetilde{\mathbb Z}^d\), extended continuously to vertices.
The vertex trace \(X\) of \(Z\), parametrized by accumulated
vertex local time, is a continuous-time simple random walk with
exponential holding times of mean one. Its jump chain \(Y\) is a
discrete-time simple random walk. 
The Green functions of \(Z\), \(X\), and \(Y\) agree at vertices,
which is given by
\[
 G(x,y)=\int_0^\infty\mathbb P_x(X_t=y)\,dt
       =\sum_{n\ge0}\mathbb P_x(Y_n=y).
\]
By the Green asymptotic in~\cite[Theorem~4.3.1]{LL10}, we have
\begin{equation}\label{eq:green-estimates}
 G(x,y)\asymp(1+|x-y|)^{2-d},\qquad
 \sum_{y\in B(x,r)}G(x,y)\lesssim r^2.
\end{equation}
For a finite discrete set \(V\), let \(H_V\) and \(\widehat H_V\)
be the first entrance time and the first positive return time of
\(Y\), respectively. Define the equilibrium measure and capacity of \(V\) by
\begin{equation}\label{eq:capacity}
 e_V(x)=\mathbf{1}_{\smash[t]{\{x\in V\}}}\mathbb P_x(\widehat H_V=\infty),
 \qquad \operatorname{cap}(V)=\sum_xe_V(x).
\end{equation}
The following capacity estimate and the last-exit decomposition
in~\cite[Propositions~6.5.2 and~4.6.4]{LL10} will be used throughout the paper,
\begin{equation}\label{eq:capacity-hitting-estimates}
 \operatorname{cap}(B_r)\asymp r^{d-2},\qquad
 \mathbb P_x(H_V<\infty)=\sum_{y\in V}G(x,y)e_V(y).
\end{equation}
We use the same entrance notation for continuous-time processes when
specified, and write \(G_V\) for the Green function of \(X\) in discrete set \(V\).
The extension of capacity to metric-graph sets in~\cite[Section~2]{DPR24}
agrees with~\eqref{eq:capacity} on discrete sets. 
For \(V\subset\mathbb Z^d\), let \(V^+\) consist of \(V\)
and all its nearest neighbors in $\mathbb Z^d$.
By subadditivity and translation invariance of capacity, we have
\begin{equation}\label{eq:neighbor-capacity-comparison}
 \operatorname{cap}(V)
 \le\operatorname{cap}(V^+)\le(2d+1)\operatorname{cap}(V).
\end{equation}

Write \(\widetilde p(t;x,y)\) for the transition density of the Brownian motion \(Z\)
relative to \(\mathfrak m\), and \(\widetilde{\mathbb P}_t^{x,y}\) for its bridge law
from \(x\) to \(y\) in time \(t>0\).
Following~\cite[Section~2]{Lup16}, the Brownian loop measure
on \(\widetilde{\mathbb Z}^d\) is defined by
\[
 \mu^{\mathrm{loop}}:=\int_{\widetilde{\mathbb Z}^d}\mathfrak m(dz)\int_0^\infty
       \frac{\widetilde p(t;z,z)}{t}\,
       \widetilde{\mathbb P}_t^{z,z}(d\gamma)\,dt.
\]
The critical loop soup \(\mathcal L\) on \(\widetilde{\mathbb Z}^d\) is the Poisson point process with
intensity measure \(\frac12\mu^{\mathrm{loop}}\). 
For an open subset \(\widetilde D\) of \(\widetilde{\mathbb Z}^d\), let
\(\mu_{\widetilde D}^{\mathrm{loop}}\) be the restriction of
\(\mu^{\mathrm{loop}}\) to loops which lie in \(\widetilde D\),
and let \(\mathcal L_{\widetilde D}\) consist of these loops in
\(\mathcal L\). By the restriction property,
\(\mathcal L_{\widetilde D}\) is a Poisson point process with intensity measure
\(\frac12\mu_{\widetilde D}^{\mathrm{loop}}\).
Two loops belong to the same cluster if they can be joined by a finite
chain of overlapping loops. We identify each cluster with the union
of its loop traces.

For a loop \(\gamma\) and a point \(z\in\widetilde{\mathbb Z}^d\),
let \(\ell_z(\gamma)\) be its total local time
at \(z\). We define the occupation fields of \(\mathcal L\) and \(\mathcal L_{\widetilde D}\) by
\(\widehat{\mathcal L}_z=\sum_{\gamma\in\mathcal L}\ell_z(\gamma)\) and 
\(\widehat{\mathcal L}_{\widetilde D,z}=\sum_{\gamma\in\mathcal L_{\widetilde D}}\ell_z(\gamma)\), respectively.
Let \(G_{\widetilde D}\) be the Green function of \(Z\)
in \(\widetilde D\), which agrees with the previous
definition when restricted to discrete sets.
We also recall the following identity
from~\cite[Section~2.4, Proposition~5]{LJ11}
and~\cite[Proposition~2.1]{Lup16}. For every open domain
\(\widetilde D\subset\widetilde{\mathbb Z}^d\) and
\(x,y\in\widetilde D\),
\begin{equation}\label{eq:loop-second-occupation}
 \mu_{\widetilde D}^{\mathrm{loop}}[\ell_x\ell_y]
 =G_{\widetilde D}(x,y)^2.
\end{equation}

We will also need to use loop soups with killings later.
For \(k\ge0\), killing at rate \(k(x)\) per unit local time at
each vertex \(x\) produces a loop soup \(\mathcal L^k\) with intensity measure
\begin{equation}\label{eq:kill}
 \frac12\exp\!\left\{-\sum_{x\in\mathbb Z^d}
             k(x)\ell_x(\gamma)\right\}\mu^{\mathrm{loop}}(d\gamma).
\end{equation}
We use superscripts \(k\) for the corresponding laws, kernels,
and occupation fields.

\Needspace{8\baselineskip}
Let \(\phi_{\widetilde D}\) be the Gaussian
free field (GFF) on \(\widetilde D\) with Dirichlet boundary conditions, that is,
a centered Gaussian field on \(\widetilde D\) with covariance \(G_{\widetilde D}\), with 
zero boundary values. We suppress the subscript when 
\(\widetilde D=\widetilde{\mathbb Z}^d\). In particular,
\begin{equation}\label{eq:phi0}
 \mathbb E[|\phi_0|]=\sqrt{2G(0,0)/\pi},\qquad
 \mathbb E[\phi_0^2]=G(0,0).
\end{equation}
We use the following coupling of the loop soup and the GFF in~\cite{Lup16}.

\begin{lemma}[Isomorphism]\label{thm:lupu-isomorphism}
For any open subset \(\widetilde D\) of \(\widetilde{\mathbb Z}^d\),  
the loop soup \(\mathcal L_{\widetilde D}\) and the GFF \(\phi_{\widetilde D}\)
admit a coupling such that, almost surely,
\begin{equation}\label{eq:lupu-occupation-identity}
 \widehat{\mathcal L}_{\widetilde D,z}=\frac12\phi_{\widetilde D,z}^2
 \quad\text{for every }z\in \widetilde D.
\end{equation}
Moreover, the loop clusters coincide with the connected components of
\(\{z:\phi_{\widetilde D,z}\ne0\}\). Conditional on the loop soup, the signs on
these components are independent and uniform in \(\{-1,1\}\).
The same assertions hold after adding nonnegative killing.
\end{lemma}

We also record the following two-point function from~\cite[Proposition~5.2]{Lup16},
for every \(d\ge3\) and \(x,y\in \widetilde D\),
\begin{equation}\label{eq:two-point}
 \mathbb P_{\widetilde D}(x\xleftrightarrow{\mathcal L_{\widetilde D}} y)
 =\frac2\pi\arcsin\frac{G_{\widetilde D}(x,y)}{
       \sqrt{G_{\widetilde D}(x,x)G_{\widetilde D}(y,y)}}\lesssim G(x,y).
\end{equation}

\subsection{One-arm lower bound and capacity estimates}

We first give a weighted two-point identity.
\begin{lemma}\label{lem:two-point}
For every \(d\ge3\) and \(x,y\in\widetilde{\mathbb Z}^d\),
\begin{equation}\label{eq:weighted-two-point}
 \mathbb E[|\phi_x|\mathbf{1}_{\smash[t]{\{x\xleftrightarrow{\mathcal L} y\}}}]=\sqrt{\frac{2}{\pi G(y,y)}}G(x,y).
\end{equation}
\end{lemma}

\begin{proof}
By Lemma~\ref{thm:lupu-isomorphism}, we have
\begin{equation}\label{eq:conditional-cluster-signs}
 \begin{aligned}
 \mathbb E[\phi_x
     \operatorname{sgn}\phi_y
       \mid\widehat{\mathcal L}]
 ={}&\sqrt{2\widehat{\mathcal L}_x}
       \mathbf{1}_{\smash[t]{\{x\xleftrightarrow{\mathcal L}y\}}}.
 \end{aligned}
\end{equation}
Taking expectations in~\eqref{eq:conditional-cluster-signs}, we obtain
\begin{equation}
 \begin{aligned}
 \mathbb E[|\phi_x|\mathbf{1}_{\smash[t]{\{x\xleftrightarrow{\mathcal L}y\}}}]
 ={}&\mathbb E[\phi_x
                 \operatorname{sgn}\phi_y]
 =\sqrt{\frac{2}{\pi G(y,y)}}G(x,y).
 \end{aligned}
\end{equation}
This proves
\eqref{eq:weighted-two-point}.
\end{proof}

Recall that \(\widetilde{\mathcal C}(0)\) is the metric-graph cluster of the origin
in \(\mathcal L\), and that \(\mathcal C(0)=\widetilde{\mathcal C}(0)\cap\mathbb Z^d\)
is its discrete version.
We record the following truncated capacity moment estimates.

\begin{lemma}\label{lem:capacity-input}
For every \(d\ge3\) and \(T\ge1\), we have
\begin{equation}\label{eq:capacity-inputs}
 \begin{aligned}
 \mathbb E[\operatorname{cap}(\mathcal C(0))\wedge T]
 \asymp{}&\sqrt T,\qquad
 \mathbb E[|\phi_0|(\operatorname{cap}(\widetilde{\mathcal C}(0))\wedge T)]
 \lesssim\sqrt T,\qquad T\ge1.
 \end{aligned}
\end{equation}
\end{lemma}

\begin{proof}
Write \(K=\operatorname{cap}(\widetilde{\mathcal C}(0))\).
By~\cite[Proposition~5]{LW18} and the fact that
\(\widetilde{\mathcal C}(0)\) is a.s.\ bounded, 
we conclude that
\begin{equation}\label{eq:conditional-capacity}
 \mathbb E[e^{-uK}\mid |\phi_0|=a]
 =\exp\!\left(-G(0,0)^{-1}u-aG(0,0)^{-1}\sqrt{2u}\right).
\end{equation}
For \(q\in\{0,1\}\), multiplying~\eqref{eq:conditional-capacity} by \(a^q\),
integrating against the law of \(|\phi_0|\), and using \(1-e^{-v}\le v\), we obtain
\begin{equation}\label{eq:weighted-capacity-laplace}
 \begin{aligned}
 \mathbb E[|\phi_0|^q(1-e^{-uK})]
 &\le\frac{u}{G(0,0)}\mathbb E[|\phi_0|^q]
       +\frac{\sqrt{2u}}{G(0,0)}\mathbb E[|\phi_0|^{q+1}]\\
 &\lesssim\sqrt u,\qquad 0<u\le1.
 \end{aligned}
\end{equation}
Taking \(u=T^{-1}\) and using
\(s\wedge T\le T(1-e^{-s/T})/(1-e^{-1})\) for \(s\ge0\), we have
\(\mathbb E[|\phi_0|^q(K\wedge T)]\lesssim\sqrt T\).
Since \(\operatorname{cap}(\mathcal C(0))\le K\) by monotonicity,
the cases \(q=0\) and \(q=1\) give the upper bounds in~\eqref{eq:capacity-inputs}.

For the lower bound in the first estimate, restricting~\eqref{eq:conditional-capacity} to
\(1\le|\phi_0|\le2\) gives
\(\mathbb E[1-e^{-uK}]\gtrsim\sqrt u\) for \(0<u\le1\). Hence
\[
 \mathbb E[K\wedge T]
 \ge T\mathbb E[1-e^{-K/T}]\gtrsim\sqrt T,\qquad T\ge1.
\]
By~\eqref{eq:neighbor-capacity-comparison}, we have
\(K\le(2d+1)\operatorname{cap}(\mathcal C(0))\).
The lower bound in~\eqref{eq:capacity-inputs} follows.
\end{proof}

The following consequence of Lemma~\ref{lem:capacity-input} turns a
mean capacity bound into a one-arm lower bound.
Recall that \(\mathcal C_r(0)=\mathcal C(0)\cap B_r\) is the truncated cluster.
\begin{lemma}\label{lem:capacity-arm-transfer}
For every \(d\ge3\) and \(r\ge1\),
\begin{equation}\label{eq:capacity-arm-transfer}
 \theta_d(r)\gtrsim
 \mathbb E[\operatorname{cap}(\mathcal C_r(0))]^{-1}
 \gtrsim\mathbb E[|\phi_0|\operatorname{cap}(\mathcal C_r(0))]^{-1}.
\end{equation}
\end{lemma}

\begin{proof}
Put \(m=\mathbb E[\operatorname{cap}(\mathcal C_r(0))]\).
On the event \(\{0\xleftrightarrow{\mathcal L}\partial\widetilde B_r\}^c\),
we have \(\mathcal C(0)=\mathcal C_r(0)\).
The lower bound in Lemma~\ref{lem:capacity-input} therefore gives,
for every \(T\ge1\),
\[
 c_0\sqrt T
 \le\mathbb E[\operatorname{cap}(\mathcal C(0))\wedge T]
 \le m+T\theta_d(r).
\]
Since \(m\ge\operatorname{cap}(\{0\})>0\), we take \(T=Am^2\)
with a sufficiently large constant \(A\) to obtain
\(\theta_d(r)\ge(c_0\sqrt A-1)/(Am)\gtrsim m^{-1}\).
Since both \(|\phi_0|=\sqrt{2\widehat{\mathcal L}_0}\) and
\(\operatorname{cap}(\mathcal C_r(0))\) are increasing in the loop soup,
we apply the FKG inequality~\cite[Theorem~20.4]{LP18} to obtain
\begin{equation}\label{eq:capacity-fkg-comparison}
\mathbb E[|\phi_0|\operatorname{cap}(\mathcal C_r(0))] \ge 
              \mathbb E|\phi_0| \, \mathbb E[\operatorname{cap}(\mathcal C_r(0))].
\end{equation}
By \eqref{eq:phi0}, the second inequality follows.
\end{proof}

Using the capacity estimates above, we give an alternative proof of the
one-arm lower bounds below, which are sharp for \(d\ne6\).
For \(d=3,4,5\), the lower bound was first obtained
in~\cite[Theorem~5]{DW20} using exploration martingales.

\begin{lemma}[One-arm lower bound]\label{lem:one-arm-lower}
For every \(d\ge3\), there is \(c_d>0\) such that
\begin{equation}\label{eq:one-arm-lower}
 \theta_d(r)\ge c_d r^{-(\frac{d-2}{2}\wedge2)},\qquad r\ge1.
\end{equation}
\end{lemma}

\begin{proof}
Since
\(\operatorname{cap}(\mathcal C_r(0))
\le\operatorname{cap}(\mathcal C(0))\wedge\operatorname{cap}(B_r)\),
Lemma~\ref{lem:capacity-input} and~\eqref{eq:capacity-hitting-estimates} give
\[
 \mathbb E[\operatorname{cap}(\mathcal C_r(0))]
 \lesssim\operatorname{cap}(B_r)^{\frac12}
 \asymp r^{\frac{d-2}{2}}.
\]
On the other hand, since
\(\operatorname{cap}(\mathcal C_r(0))\le|\mathcal C_r(0)|\),
the two-point bound~\eqref{eq:two-point} and
\eqref{eq:green-estimates} give
\[
 \mathbb E[\operatorname{cap}(\mathcal C_r(0))]
 \le\mathbb E|\mathcal C_r(0)|
 =\sum_{x\in B_r}\mathbb P(0\xleftrightarrow{\mathcal L}x)
 \lesssim\sum_{x\in B_r}G(0,x)\lesssim r^2.
\]
Combining these bounds with the first inequality
in~\eqref{eq:capacity-arm-transfer}, we obtain
\[
 \theta_d(r)\gtrsim
 \mathbb E[\operatorname{cap}(\mathcal C_r(0))]^{-1}
 \gtrsim r^{-(\frac{d-2}{2}\wedge2)}.
\]
This proves~\eqref{eq:one-arm-lower}.
\end{proof}

\subsection{Weighted one-arm estimates}

Recall the weighted one-arm probabilities \(w_d(r)\) and \(w_d^{\le\delta r}(r)\)
from~\eqref{eq:wr} and~\eqref{eq:wdr}. The following lemma proves
\eqref{eq:large-loop-decomposition} by an argument analogous to
\cite[Proposition~4.5]{Gao26}.
\begin{lemma}\label{lem:large-loop-removal}
Write \(\beta=(d-2)/2\). For every \(d\ge3\), there exists
\(C=C(d)>0\) such that, for all \(0<\delta<1/8\) and \(r\ge1\),
\[
 w_d(r)\le w_d^{\le\delta r}(r)+C\delta^{-(2d-2)}r^{-\beta}.
\]
\end{lemma}

\begin{proof}
Write \(\mathcal E_r^{\le\delta r}
 =\{0\xleftrightarrow{\mathcal L^{\le\delta r}}\partial\widetilde B_r\}\).
On the event \((\mathcal E_r^{\le\delta r})^c\), the cluster \(K_0\)
of zero in \(\mathcal L^{\le\delta r}\) is contained in \(\widetilde B_r\).
Now if \(0\xleftrightarrow{\mathcal L}\partial\widetilde B_r\) occurs,
then there is a loop
of diameter greater than \(\delta r\) that intersects \(K_0\).
Cover \(\widetilde B_r\) by \(O(\delta^{-d})\) boxes of radius \(\delta r/32\),
and let \(K_i\) be the part of \(K_0\) in the \(i\)-th box.
By \cite[Lemma 2.7]{CS16}
and the fact that \(\operatorname{cap}(K_i)\le\operatorname{cap}(K_0)\), we obtain that
\begin{equation}\label{eq:large-loop-attachment}
 \begin{aligned}
 1-\exp\!\left\{-\frac12\mu^{\mathrm{loop}}
       (\gamma\cap K_0\ne\varnothing,\operatorname{diam}\gamma>\delta r)\right\}
 &\lesssim (\delta r)^{2-d}\sum_i\operatorname{cap}(K_i)\\
 &\lesssim\delta^{-(2d-2)}r^{2-d}\operatorname{cap}(K_0).
 \end{aligned}
\end{equation}
Since \(\sqrt{2\widehat{\mathcal L}^{\le\delta r}_{0}}\le |\phi_0|\) and \(\operatorname{cap}(K_0)\lesssim \operatorname{cap}(\widetilde{\mathcal C}(0))\wedge r^{d-2}\),
we obtain from \eqref{eq:large-loop-attachment} and the weighted capacity moment in Lemma~\ref{lem:capacity-input} that
\begin{equation}\label{eq:large-loop-connection-error}
 \begin{aligned}
 \mathbb E\!\left[\sqrt{2\widehat{\mathcal L}^{\le\delta r}_{0}}\mathbf{1}_{\smash[t]{\{0\xleftrightarrow{\mathcal L}\partial\widetilde B_r\}\setminus\mathcal E_r^{\le\delta r}}}\right]
 &\lesssim\delta^{-(2d-2)}r^{2-d}
       \mathbb E[|\phi_0|
                  (\operatorname{cap}(\widetilde{\mathcal C}(0))\wedge r^{d-2})]\\
 &\lesssim\delta^{-(2d-2)}r^{-\beta}.
 \end{aligned}
\end{equation}

Write \(\widehat{\mathcal L}^{>\delta r}_{0}:=\widehat{\mathcal L}_0-\widehat{\mathcal L}^{\le\delta r}_{0}\) for the occupation at \(0\) generated by large loops. 
By \cite[Section~2.4, Proposition~5]{LJ11} and~\cite[Proposition~4.6.2]{LL10}, we have
\begin{equation}\label{eq:large-loop-occupation}
 2\mathbb E\widehat{\mathcal L}^{>\delta r}_{0}
 \le G(0,0)-G_{\widetilde B_{\delta r/3}}(0,0)
   \lesssim (\delta r)^{2-d}.
\end{equation}
By the Cauchy--Schwarz inequality and~\eqref{eq:large-loop-occupation},
\begin{equation}\label{eq:large-loop-root-amplitude}
 \mathbb E\sqrt{2\widehat{\mathcal L}^{>\delta r}_{0}}
 \le\sqrt{2\mathbb E\widehat{\mathcal L}^{>\delta r}_{0}}
   \lesssim (\delta r)^{-\beta}.
\end{equation}
Using \(|\phi_0|\le \sqrt{2\widehat{\mathcal L}^{\le\delta r}_{0}}+\sqrt{2\widehat{\mathcal L}^{>\delta r}_{0}}\)
and substituting~\eqref{eq:large-loop-connection-error}
and~\eqref{eq:large-loop-root-amplitude}, we conclude that
\[
 \begin{aligned}
 w_d(r)
 ={}&\mathbb E[|\phi_0|\mathbf{1}_{\smash[t]{\{0\xleftrightarrow{\mathcal L}\partial\widetilde B_r\}}}]\\
 \le{}& w_d^{\le\delta r}(r)
    +\mathbb E\!\left[\sqrt{2\widehat{\mathcal L}^{\le\delta r}_{0}}\mathbf{1}_{\smash[t]{\{0\xleftrightarrow{\mathcal L}\partial\widetilde B_r\}\setminus\mathcal E_r^{\le\delta r}}}\right]
    +\mathbb E\sqrt{2\widehat{\mathcal L}^{>\delta r}_{0}}\\
 \le{}& w_d^{\le\delta r}(r)+C\delta^{-(2d-2)}r^{-\beta}.
 \end{aligned}
\]
We complete the proof of the lemma.
\end{proof}

Let \(X\) be continuous-time simple random walk with exponential holding times of mean one,
started at \(z\) and killed at rate \(k(x)\) at vertex \(x\) under \(\mathbb P_z^k\).
We write \(H_\Sigma\) for its first hitting time of \(\Sigma\)
and \(\zeta\) for its lifetime. The following inequality adapts the
argument of~\cite[Lemma~3.5]{CD24} to nonnegative vertex killing.

\begin{lemma}\label{lem:separator}
Let \(V\) be a finite set in $\mathbb Z^d$ and \(z\notin V\). Let
\(\Sigma\subset\mathbb Z^d\setminus(V\cup\{z\})\) be a deterministic set separating
\(z\) from \(V\). For nonnegative vertex killing
\(k\),
\begin{equation}\label{eq:separator}
 \mathbb E^k\!\left[\sqrt{2\widehat{\mathcal L}^{k}_{z}}\mathbf{1}_{\smash[t]{\{z\xleftrightarrow{\mathcal L^k}V\}}}\right]
 \le\sum_{u\in\Sigma}\mathbb P_z^k(X_{H_\Sigma}=u,H_\Sigma<\zeta)
       \mathbb E^k\!\left[\sqrt{2\widehat{\mathcal L}^{k}_{u}}\mathbf{1}_{\smash[t]{\{u\xleftrightarrow{\mathcal L^k}V\}}}\right].
\end{equation}
\end{lemma}

\begin{proof}
Let
\(\phi^k\) be the corresponding GFF on \(\widetilde{\mathbb Z}^d\), coupled
with the corresponding loop soup so that
\(|\phi_x^k|=\sqrt{2\widehat{\mathcal L}^{k}_{x}}\).
Let \(A\) be the closure of the union of \(V\) and all negative sign
clusters of \(\phi^k\) that intersect \(V\). By~\cite[Section~2]{DW20},
\(A\) is a local set for \(\phi^k\). Let \(h(z)=\mathbb E^k[\phi_z^k\mid A,\phi^k|_A]\) for \(z\in A^c\),
which is the harmonic extension of
\(\phi^k|_{\partial A}\) to \(A^c\) with killing \(k\).
Extend \(h\) by zero on \(A\setminus V\) and by
\(\phi^k\vee0\) on \(V\). 
Since \(\mathbb E^k\phi_z^k=0\), for \(z\notin V\) we have
\[
 \mathbb E^k h(z)
 =\mathbb E^k[\phi_z^k\mathbf{1}_{\smash[t]{\{z\notin A\}}}]
   =-\mathbb E^k[\phi_z^k\mathbf{1}_{\smash[t]{\{z\in A\}}}] =\mathbb E^k[|\phi_z^k|\mathbf{1}_{\smash[t]{\{z\xleftrightarrow{\mathcal L^k}V,
       \phi_z^k<0\}}}].
\]
By Lemma~\ref{thm:lupu-isomorphism} and the preceding identity, we obtain
\begin{equation}\label{eq:local-set-mean}
 \mathbb E^k h(z)=\frac12\mathbb E^k\!\left[
       \sqrt{2\widehat{\mathcal L}^{k}_{z}}\mathbf{1}_{\smash[t]{\{z\xleftrightarrow{\mathcal L^k}V\}}}\right].
\end{equation}
Since \(\phi^k=0\) on \(\partial A\setminus V\) and
\(\phi^k\ge0\) on \(V\cap\partial A\), we have \(h\ge0\).
By its harmonic representation, for \(z\notin A\) we have
\[
 h(z)=\mathbb E_z^k[h(Z_{H_\Sigma})\mathbf{1}_{\smash[t]{\{H_\Sigma<H_A,\ H_\Sigma<\zeta\}}}] \le\sum_{u\in\Sigma}
        \mathbb P_z^k(X_{H_\Sigma}=u,H_\Sigma<\zeta)h(u),
\]
where $Z$ is the Brownian motion on $\widetilde{\mathbb Z}^d$ started from $z$ with killing $k$.
For \(z\in A\setminus V\), the same inequality holds because \(h(z)=0\).
Taking expectations and applying~\eqref{eq:local-set-mean} at \(z\) and
at each \(u\in\Sigma\) proves~\eqref{eq:separator}.
\end{proof}

Generally speaking, Lemma~\ref{lem:separator} can be viewed as a BK-type (or Simon--Lieb-type)
inequality for the weighted one-arm probability, analogous to the
FK--Ising bound in~\cite[Proposition~5.12]{EGPS25}.

Now, for a nonnegative vertex killing rate \(k\), we define the weighted one-arm probability in the loop soup \(\mathcal L^k\) by
\begin{equation}\label{eq:wdk}
 w_d^k(r)=\mathbb E^k\!\left[\sqrt{2\widehat{\mathcal L}^{k}_{0}}\mathbf{1}_{\smash[t]{\{0\xleftrightarrow{\mathcal L^k}\partial\widetilde B_r\}}}\right].
\end{equation}
For \(d\ge3\), \(r\ge1\), and \(\kappa\ge1\), let
\(\boldsymbol{k}\) be the killing rate in \eqref{eq:protected-killing-rate}.
Then, we have the following contraction estimate under killing.

\begin{lemma}\label{lem:killed-contraction}
For every \(d\ge3\), there exist constants \(c>0\) and \(C\ge1\)
such that, for every \(\kappa\ge1\) and \(r\ge C\kappa\), we have
\(w_d^{\boldsymbol{k}}(r)\le Ce^{-c\kappa}w_d(r/4)\).
\end{lemma}

\begin{proof}
Increase \(C\) so that \(r\ge8\). Every connection in the metric graph from zero to
\(\partial\widetilde B_r\) visits the inner vertex boundary \(\partial B_r\).
For \(u\in\partial B_{r/2}\) and \(v\in\partial B_r\), we have
\(|v-u|\ge r/2-1\ge r/4\), including for noninteger \(r\).
Thus, by monotonicity and translation invariance, for every
\(u\in\partial B_{r/2}\) we have
\begin{equation}\label{eq:killed-terminal-connection}
 \mathbb E^{\boldsymbol{k}}\!\left[\sqrt{2\widehat{\mathcal L}^{\boldsymbol{k}}_{u}}\mathbf{1}_{\smash[t]{\{u\xleftrightarrow{\mathcal L^{\boldsymbol{k}}}\partial B_r\}}}\right]\le w_d(r/4).
\end{equation}
Applying Lemma~\ref{lem:separator} with \(V=\partial B_r\) and \(\Sigma=\partial B_{r/2}\) and using~\eqref{eq:separator} and~\eqref{eq:killed-terminal-connection},
we obtain
\[
 \begin{aligned}
 w_d^{\boldsymbol{k}}(r)
 \le{}&\mathbb E^{\boldsymbol{k}}\!\left[\sqrt{2\widehat{\mathcal L}^{\boldsymbol{k}}_{0}}\mathbf{1}_{\smash[t]{\{0\xleftrightarrow{\mathcal L^{\boldsymbol{k}}}\partial B_r\}}}\right]\\
 \le{}&\sum_{u\in\partial B_{r/2}}
       \mathbb P_0^{\boldsymbol{k}}(X_{H_{\partial B_{r/2}}}=u,H_{\partial B_{r/2}}<\zeta)
       \mathbb E^{\boldsymbol{k}}\!\left[\sqrt{2\widehat{\mathcal L}^{\boldsymbol{k}}_{u}}\mathbf{1}_{\smash[t]{\{u\xleftrightarrow{\mathcal L^{\boldsymbol{k}}}\partial B_r\}}}\right]\\
 \le{}&\mathbb P_0^{\boldsymbol{k}}(H_{\partial B_{r/2}}<\zeta)w_d(r/4)
 \le Ce^{-c\kappa}w_d(r/4),
 \end{aligned}
\]
where the last inequality uses
\(\mathbb P_0^{\boldsymbol{k}}(H_{\partial B_{r/2}}<\zeta)\le Ce^{-c\kappa}\),
which follows from standard estimates in
\cite[Section~1.2 and Proposition~2.1.2]{LL10}.
\end{proof}

\begin{lemma}\label{lem:small-ball}
For every \(d\ge3\), there exists a constant \(C_d\ge1\)
such that, for every \(s\ge1\) and every \(x\in\mathbb Z^d\) with \(|x|\ge C_ds\),
\begin{equation}\label{eq:small-ball}
 \mathbb E[|\phi_0|\mathbf{1}_{\smash[t]{\{0\xleftrightarrow{\mathcal L} B(x,s)\}}}]
 \lesssim \left(\frac{s}{|x|}\right)^{d-2}w_d(s).
\end{equation}
\end{lemma}

\begin{proof}
We will apply Lemma~\ref{lem:separator} with \(k=0\), \(V=B(x,s)\), and
\(\Sigma=\partial B(x,3s)\), with
integer rounding if needed. If \(|x|\ge C_ds\) for some large $C_d$, then $\Sigma$ 
separates zero from the target $V$. The standard hitting estimate~\cite[Propositions~6.5.1--6.5.2]{LL10} gives
\begin{equation}\label{eq:small-ball-hitting}
 \mathbb P_0(H_{\partial B(x,3s)}<\infty)
 \lesssim \left(\frac{s}{|x|}\right)^{d-2}.
\end{equation}
Every connection from \(u\in\partial B(x,3s)\) to \(B(x,s)\) reaches distance
at least \(s\) from \(u\), and hence the weighted one-arm probability is at
most \(w_d(s)\) by translation invariance. Substituting this bound
and~\eqref{eq:small-ball-hitting} into~\eqref{eq:separator} concludes the proof.
\end{proof}

\section{One-arm upper bounds for \texorpdfstring{\(d\ne6\)}{d not equal to 6}}
\label{sec:exploration}

We first prove the upper bound in dimension three.
This is the simplest case and serves as a good example to 
demonstrate the main idea.
We then use the exploration entropy method of~\cite{DM23PB,DM23}
to refine the estimates and give a unified
proof of Theorem~\ref{thm:main} for \(d\ne6\).

\subsection{A direct proof in dimension three}
\label{subsec:three-dimensional}

We prove the one-arm upper bound for $d=3$ as a warm-up case. 
We first recall the relative entropy
that will be used in this section.
For probability measures \(P,Q\), define the relative entropy by
\(\mathrm{H}(P\Vert Q)=\int\log(dP/dQ)\,dP\) when
\(P\ll Q\), and \(+\infty\) otherwise.
In the case of Poisson point processes, we have the following result.
\Needspace{11\baselineskip}
\begin{lemma}[{\cite[Theorem~6]{Les24}}]
\label{lem:poisson-relative-entropy}
Let \(\mu\) be a sigma-finite measure on a measurable space \(E\),
and let \(f:E\to[0,\infty)\) be measurable.
Writing \(\operatorname{PPP}(\nu)\) for the law of a Poisson point
process with intensity measure \(\nu\), we have
\begin{equation}
 \mathrm{H}\bigl(\operatorname{PPP}(\mu)\Vert
                         \operatorname{PPP}(f\mu)\bigr)
 =\int_E\bigl(f-1-\log f\bigr)d\mu.
\end{equation}
The integrand is interpreted as \(+\infty\) where \(f=0\).
The equality holds in \([0,\infty]\).
\end{lemma}

Next, we prove the one-arm upper bound for \(d=3\) within three steps. 
The proofs for \(d\ge 4\) are dedicated to refining the crude entropy bound in the 
first step.

\begin{proof}[Proof of the upper bound in Theorem~\ref{thm:main} for \(d=3\)]
Fix \(0<\delta<1/32\), \(\kappa\ge1\), and
\(r\ge C\max\{\delta^{-1},\kappa\}\), with a sufficiently large
absolute constant \(C\). Use the killing
\(\boldsymbol{k}\) from~\eqref{eq:protected-killing-rate}, and set
$
 \Gamma^{\le\delta r}=\{\gamma\subset\widetilde B_{2r}:
                         \operatorname{diam}\gamma\le\delta r\}
$.
Let \(\widehat P,\widehat P^{\boldsymbol{k}}\) be the
laws of the loop soups with intensities
\(\frac12\mu_{\widetilde B_{2r}}^{\mathrm{loop}}|_{\Gamma^{\le\delta r}}\)
and \(\frac12e^{-\boldsymbol J}
\mu_{\widetilde B_{2r}}^{\mathrm{loop}}|_{\Gamma^{\le\delta r}}\), respectively.
Here \(\boldsymbol J(\gamma)=\sum_x\boldsymbol{k}(x)\ell_x(\gamma)
\le\kappa^2r^{-2}\mathcal V(\gamma)\).

We further split \(\Gamma^{\le\delta r}\) into loops that intersect
\(\widetilde B_{r/16}\) and loops that avoid it, denoted by 
\(\Gamma^{\le\delta r}_{\mathrm{in}}\) and \(\Gamma^{\le\delta r}_{\mathrm{out}}\), respectively.
Then, \(\boldsymbol J(\gamma)=0\) for all $\gamma\in \Gamma^{\le\delta r}_{\mathrm{in}}$
since these loops have radius \((1/16+\delta)r<r/8\).
Let \(\widehat P_{\mathrm{out}},\widehat P_{\mathrm{out}}^{\boldsymbol{k}}\)
be the respective laws of \(\widehat P,\widehat P^{\boldsymbol{k}}\) 
restricted on \(\Gamma^{\le\delta r}_{\mathrm{out}}\).

\medskip\noindent\textbf{Step 1: The entropy bound.}
Lemma~\ref{lem:poisson-relative-entropy},
\(t-1+e^{-t}\le t^2/2\), and~\eqref{eq:loop-second-occupation} give
\begin{equation}\label{eq:three-static-poisson-entropy}
 \begin{aligned}
 \mathrm H(\widehat P_{\mathrm{out}}\Vert
                   \widehat P_{\mathrm{out}}^{\boldsymbol{k}})
 &\le \frac14\int_{\Gamma^{\le\delta r}}\boldsymbol J(\gamma)^2
                   \,\mu_{\widetilde B_{2r}}^{\mathrm{loop}}(d\gamma)\\
 &\le C\kappa^4r^{-4}
       \sum_{x\in B_{2r}}\sum_{|y-x|\le\delta r}G(x,y)^2\\
 &\le C\kappa^4r^{-1}(\delta r+1)
 \le C\kappa^4\delta.
 \end{aligned}
\end{equation}
The last line uses \(d=3\), \eqref{eq:green-estimates}, and
\(\delta r\ge1\). 
Now, Pinsker's inequality~\cite[Remark~2.4]{DM23} consequently bounds the total variation 
distance by
\begin{equation}\label{eq:three-static-total-variation}
 \bigl\|\widehat P_{\mathrm{out}}
       -\widehat P_{\mathrm{out}}^{\boldsymbol{k}}\bigr\|_{\mathrm{TV}}
 \le \sqrt{\frac12\mathrm H\bigl(\widehat P_{\mathrm{out}}\Vert
                   \widehat P_{\mathrm{out}}^{\boldsymbol{k}}\bigr)}
 \le C\kappa^2\sqrt\delta.
\end{equation}

\medskip\noindent\textbf{Step 2: The weighted comparison.}
If \(0\not\xleftrightarrow{}\partial\widetilde B_{r/16}\), the
conditional probability of
\(0\xleftrightarrow{}\partial\widetilde B_r\) is zero under
both laws on \(\Gamma^{\le\delta r}_{\mathrm{out}}\). Otherwise, 
their difference is bounded by the total
variation distance. Since the event
\(\{0\xleftrightarrow{}\partial\widetilde B_{r/16}\}\) and the weight
\(\sqrt{2\widehat{\mathcal L}^{\le\delta r}_0}\) are determined by the
configuration on \(\Gamma^{\le\delta r}_{\mathrm{in}}\),
independent of loops in \(\Gamma^{\le\delta r}_{\mathrm{out}}\),
we obtain that 
\begin{equation}\label{eq:three-static-comparison}
 \begin{aligned}
 w_3^{\le\delta r}(r)
 &=\mathbb E_{\widehat P}\Bigl[\sqrt{2\widehat{\mathcal L}^{\le\delta r}_0}
       \,\mathbf1_{\{0\xleftrightarrow{}\partial\widetilde B_r\}}\Bigr]\\
 &\le \mathbb E_{\widehat P^{\boldsymbol{k}}}\Bigl[
       \sqrt{2\widehat{\mathcal L}^{\le\delta r}_0}
       \,\mathbf1_{\{0\xleftrightarrow{}\partial\widetilde B_r\}}\Bigr]
       +\bigl\|\widehat P_{\mathrm{out}}
       -\widehat P_{\mathrm{out}}^{\boldsymbol{k}}\bigr\|_{\mathrm{TV}}
       \,\mathbb E_{\widehat P}\Bigl[
       \sqrt{2\widehat{\mathcal L}^{\le\delta r}_0}
       \,\mathbf1_{\{0\xleftrightarrow{}\partial\widetilde B_{r/16}\}}\Bigr]\\
 &\le w_3^{\boldsymbol{k}}(r)
       +C\kappa^2\sqrt\delta\,w_3(r/16).
 \end{aligned}
\end{equation}
The last inequality follows from~\eqref{eq:three-static-total-variation}
and the definition \eqref{eq:wdr}.
Combining Lemmas~\ref{lem:large-loop-removal}
and~\ref{lem:killed-contraction}, we get
\begin{equation}\label{eq:three-static-renormalisation}
 w_3(r)\le C_\delta r^{-1/2}
       +Ce^{-c\kappa}w_3(r/4)
       +C\kappa^2\sqrt\delta\,w_3(r/16).
\end{equation}

\medskip\noindent\textbf{Step 3: Closing the induction.}
Choose \(\kappa\) large and then \(\delta\) small such that
$
 2Ce^{-c\kappa}+4C\kappa^2\sqrt\delta\le\frac12
$.
Take a dyadic radius \(r_0\) satisfying the preceding scale
conditions, and choose \(M\ge2C_\delta\) so that
\(w_3(s)\le Ms^{-1/2}\) for \(1\le s\le r_0\).
Induction over dyadic radii in~\eqref{eq:three-static-renormalisation}
then gives \(w_3(r)\le Mr^{-1/2}\).
We can use monotonicity to handle arbitrary radii after increasing the constant.
Finally, the FKG inequality gives
\(w_3(r)\ge\mathbb E|\phi_0|\,\theta_3(r)\), proving the desired
upper bound. We finish the proof for $d=3$.
\end{proof}

\subsection{The exploration entropy bound}\label{subsec:entropy}

For an integer \(r\ge1\), we work with the critical loop soup
\(\mathcal L_{\widetilde B_{2r}}\)
in the box \(\widetilde B_{2r}\).
Let \(\zeta(\gamma)\) be the duration of a loop \(\gamma\).
For \(\varepsilon>0\), we retain the loops in
\begin{equation}\label{eq:finite-loop-class}
 \Gamma_\varepsilon^{\le\delta r}=
 \{\gamma:\operatorname{diam}\gamma\le\delta r,\ 
                              \zeta(\gamma)\ge\varepsilon\},
\end{equation}
which has finite \(\mu_{\widetilde B_{2r}}^{\mathrm{loop}}\)-measure.

Let \(P_\varepsilon\) be the law of the loop soup with intensity
\(\frac12\mu_{\widetilde B_{2r}}^{\mathrm{loop}}\) restricted to
\(\Gamma_\varepsilon^{\le\delta r}\). Let \(P_\varepsilon^{\boldsymbol{k}}\) be the law of the loop soup
with the additional killing rate \(\boldsymbol{k}\), also restricted to
\(\Gamma_\varepsilon^{\le\delta r}\), so that the
intensity of each retained loop \(\gamma\) is multiplied by
\(e^{-\boldsymbol{J}(\gamma)}\), where
\(\boldsymbol{J}(\gamma)=\sum_{x\in\mathbb Z^d}\boldsymbol{k}(x)\ell_x(\gamma)\) (see \eqref{eq:kill}).
Under either law, write \(\mathcal L^{\le\delta r,\varepsilon}\)
for the loop soup and \(\widehat{\mathcal L}^{\le\delta r,\varepsilon}_{x}\)
for its occupation at \(x\).

\begin{exploration}\label{exp:root-cluster}
We explore the cluster of \(0\) in \(\mathcal L^{\le\delta r,\varepsilon}\)
by successively revealing loops that intersect the discovered trace.
The exploration stops when this trace meets \(\partial\widetilde B_r\)
or no new loop is found. 
\begin{enumerate}[label=\arabic*.,leftmargin=*,itemsep=3pt,topsep=4pt]
\item First reveal the loop soup on the set of loops in
\(\Gamma_\varepsilon^{\le\delta r}\) that visit zero, and call this initial
configuration \(Z_0\). Let \(K_0\) be the union of zero and the
traces of these loops. If \(K_0\) meets \(\partial\widetilde B_r\), stop.

\item At step \(j\ge1\), let \(U_j\) be the set of loops in
\(\Gamma_\varepsilon^{\le\delta r}\) that intersect \(K_{j-1}\) and lie outside
all regions of loop space on which the loop soup has already been revealed.
Reveal the loop soup restricted to \(U_j\), and denote the resulting
configuration by \(Z_j\). Add the traces of these newly discovered
loops to \(K_{j-1}\) to obtain \(K_j\).

\item Stop if \(K_j\) meets \(\partial\widetilde B_r\) or if no new loop was found;
otherwise proceed to the next step. After stopping, take all
subsequent \(U_j\) and \(Z_j\) to be empty and keep the discovered
set unchanged.
\end{enumerate}
Here \(U_j\) is a region of loop space selected from the preceding
record, whereas \(Z_j\) records the loops of the loop soup in that region.
In particular, revealing the loop soup on \(U_j\) also records the absence
of any other loops there. 
The regions \(U_j\) are disjoint from one
another and from the initial region of loops visiting zero.
\end{exploration}

For the entropy bound, we
define \(\widehat P_\varepsilon\) and \(\widehat P_\varepsilon^{\boldsymbol{k}}\) as the
laws of the exploration record obtained by reweighting
\(P_\varepsilon\) and \(P_\varepsilon^{\boldsymbol{k}}\) by the root occupation weight
\(\sqrt{2\widehat{\mathcal L}^{\le\delta r,\varepsilon}_{0}}\), respectively. 
Explicitly, for every measurable
set \(A\) of records $(Z_j)_{j\ge0}$,
\begin{equation}\label{eq:weighted-exploration-laws}
 \widehat P_\varepsilon(A)
 =\frac{\mathbb E_{P_\varepsilon}\!\left[
 \sqrt{2\widehat{\mathcal L}^{\le\delta r,\varepsilon}_{0}}
 \mathbf 1_{\{(Z_j)_{j\ge0}\in A\}}\right]}
 {\mathbb E_{P_\varepsilon}\!\left[
 \sqrt{2\widehat{\mathcal L}^{\le\delta r,\varepsilon}_{0}}\right]},\qquad
 \widehat P_\varepsilon^{\boldsymbol{k}}(A)
 =\frac{\mathbb E_{P_\varepsilon^{\boldsymbol{k}}}\!\left[
 \sqrt{2\widehat{\mathcal L}^{\le\delta r,\varepsilon}_{0}}
 \mathbf 1_{\{(Z_j)_{j\ge0}\in A\}}\right]}
 {\mathbb E_{P_\varepsilon^{\boldsymbol{k}}}\!\left[
 \sqrt{2\widehat{\mathcal L}^{\le\delta r,\varepsilon}_{0}}\right]}.
\end{equation}
Let \(\operatorname{Rev}(\gamma)=\{\gamma\in\bigcup_{j\ge1}U_j\}\) 
be the event that the candidate loop $\gamma$ is revealed.
The corresponding weighted revealment probability is 
\begin{equation}\label{eq:weighted-revealment}
 R_\varepsilon(\gamma)
 =\mathbb E_{P_\varepsilon}\!\left[
 \sqrt{2\widehat{\mathcal L}^{\le\delta r,\varepsilon}_{0}}
 \mathbf{1}_{\smash[t]{\operatorname{Rev}(\gamma)}}\right].
\end{equation}
Set \(R_\varepsilon=0\) outside \(\Gamma_\varepsilon^{\le\delta r}\).
Recall that we use 
\(\mathcal V(\gamma)=\sum_{x\in\mathbb Z^d}\ell_x(\gamma)\) to denote 
the total vertex local time of \(\gamma\).

\begin{proposition}[Entropic bound]\label{lem:query}
There exist constants \(L_0\ge1\) and \(\varepsilon_0>0\),
depending only on the dimension, such that,
whenever \(0<\delta<\frac18\), \(\delta r\ge L_0\), and
\(0<\varepsilon\le\varepsilon_0\),
Algorithm~\ref{exp:root-cluster} in \(\widetilde B_{2r}\) satisfies
\begin{equation}\label{eq:exploration-entropy}
 \begin{aligned}
 \Delta_\varepsilon(r)
 :={}&\mathbb E_{P_\varepsilon}\!\left[
       \sqrt{2\widehat{\mathcal L}^{\le\delta r,\varepsilon}_{0}}\right]
       \mathrm{H}(\widehat P_\varepsilon\Vert\widehat P_\varepsilon^{\boldsymbol{k}})\\
 \lesssim{}& \kappa^4r^{-4}
 \int_{\substack{\operatorname{diam}\gamma\le\delta r\\\boldsymbol{J}(\gamma)>0}}
 \mathcal V(\gamma)^2R_\varepsilon(\gamma)\,
                         \mu_{\widetilde B_{2r}}^{\mathrm{loop}}(d\gamma).
 \end{aligned}
\end{equation}
Moreover, for every event \(A\) determined by Algorithm~\ref{exp:root-cluster},
\begin{equation}\label{eq:entropy-event}
 \mathbb E_{P_\varepsilon}[\sqrt{2\widehat{\mathcal L}^{\le\delta r,\varepsilon}_{0}}\mathbf{1}_{\smash[t]{A}}]
 \le2\mathbb E_{P_\varepsilon^{\boldsymbol{k}}}[\sqrt{2\widehat{\mathcal L}^{\le\delta r,\varepsilon}_{0}}\mathbf{1}_{\smash[t]{A}}]+2\Delta_\varepsilon(r).
\end{equation}
\end{proposition}

\begin{proof}
Let \(U_0=\{\gamma\in\Gamma_\varepsilon^{\le\delta r}:0\in\gamma\}\)
be the initial query region. Every \(\gamma\in U_0\) contains zero
and has diameter at most \(\delta r\). Since \(\delta<\frac18\), its
trace is contained in \(\{z:|z|<r/8\}\), where the vertex
killing rate \(\boldsymbol{k}\) vanishes by
\eqref{eq:protected-killing-rate}. Thus \(\boldsymbol{J}=0\) on \(U_0\), and
the loop intensities under \(P_\varepsilon\) and \(P_\varepsilon^{\boldsymbol{k}}\)
agree there. The initial record \(Z_0\) therefore has the same
Poisson law under both measures,
\begin{equation}
 \operatorname{Law}_{P_\varepsilon}(Z_0)
 =\operatorname{Law}_{P_\varepsilon^{\boldsymbol{k}}}(Z_0)
 =\operatorname{PPP}\!\left(
       \frac12\mathbf 1_{U_0}\mu_{\widetilde B_{2r}}^{\mathrm{loop}}\right).
\end{equation}
The root occupation weight is determined by \(Z_0\), so it also
has the same law under both measures. We may therefore write
\begin{equation}\label{eq:query-weight-notation}
 W_\varepsilon=\sqrt{2\widehat{\mathcal L}^{\le\delta r,\varepsilon}_{0}},
 \qquad m_\varepsilon=\mathbb E_{P_\varepsilon}W_\varepsilon
                   =\mathbb E_{P_\varepsilon^{\boldsymbol{k}}}W_\varepsilon,
 \qquad \mathcal F_j=\sigma(Z_0,\ldots,Z_j).
\end{equation}

For suitable \(L_0\ge1\) and \(\varepsilon_0>0\),
\(m_\varepsilon\) is uniformly bounded above and away from zero
under the stated assumptions, so the weighted laws in \eqref{eq:weighted-exploration-laws} are well defined.
The weight \(W_\varepsilon\) is \(\mathcal F_0\)-measurable.
Thus, for every \(A_0\in\mathcal F_0\),
\begin{equation}
 \widehat P_\varepsilon(A_0)
 =\frac{\mathbb E_{P_\varepsilon}[W_\varepsilon\mathbf{1}_{\smash[t]{A_0}}]}{m_\varepsilon}
 =\frac{\mathbb E_{P_\varepsilon^{\boldsymbol{k}}}[W_\varepsilon\mathbf{1}_{\smash[t]{A_0}}]}{m_\varepsilon}
 =\widehat P_\varepsilon^{\boldsymbol{k}}(A_0).
\end{equation}
In particular,
\(\mathrm{H}(\widehat P_\varepsilon|_{\mathcal F_0}
\Vert\widehat P_\varepsilon^{\boldsymbol{k}}|_{\mathcal F_0})=0\).
Similarly, for \(j\ge1\) and every bounded
measurable function \(f\), 
\begin{equation}\label{eq:query-tilt-cancellation}
 \mathbb E_{\widehat P_\varepsilon}[f(Z_j)\mid\mathcal F_{j-1}]
 =\mathbb E_{P_\varepsilon}[f(Z_j)\mid\mathcal F_{j-1}].
\end{equation}
The identical calculation applies
to \(\widehat P_\varepsilon^{\boldsymbol{k}}\) and \(P_\varepsilon^{\boldsymbol{k}}\).

Conditional on the past record \((Z_0,\ldots,Z_{j-1})\), \(U_j\) is a fixed subset of
the unqueried region
\(\Gamma_\varepsilon^{\le\delta r}\setminus\bigcup_{i=0}^{j-1}U_i\).
By the restriction property and~\eqref{eq:query-tilt-cancellation},
for the same past record under both laws, we obtain
\begin{equation}\label{eq:conditional-query-intensities}
 \begin{aligned}
 \operatorname{Law}_{\widehat P_\varepsilon}(Z_j\mid\mathcal F_{j-1})
 ={}&\operatorname{PPP}(\lambda_j),
 &\lambda_j(d\gamma)
 ={}&\frac12\mathbf 1_{U_j}(\gamma)
                          \mu_{\widetilde B_{2r}}^{\mathrm{loop}}(d\gamma),\\
 \operatorname{Law}_{\widehat P_\varepsilon^{\boldsymbol{k}}}(Z_j\mid\mathcal F_{j-1})
 ={}&\operatorname{PPP}(\rho_j),
 &\rho_j(d\gamma)={}&e^{-\boldsymbol{J}(\gamma)}\lambda_j(d\gamma).
 \end{aligned}
\end{equation}
The conditional entropy is therefore, by
Lemma~\ref{lem:poisson-relative-entropy} with \(\mu=\lambda_j\)
and \(f=e^{-\boldsymbol{J}}\),
\begin{equation}
 h_j:=\mathrm{H}\bigl(
           \operatorname{PPP}(\lambda_j)\Vert\operatorname{PPP}(\rho_j)\bigr)
 =\frac12\int_{U_j}
       \bigl(\boldsymbol{J}(\gamma)-1+e^{-\boldsymbol{J}(\gamma)}\bigr)
                          \mu_{\widetilde B_{2r}}^{\mathrm{loop}}(d\gamma).
\end{equation}
Using \(t-1+e^{-t}\le t^2/2\) for \(t\ge0\) and the chain rule of
relative entropy~\cite[(2.2) and Lemma~2.7]{DM23PB}, we obtain
\begin{equation}
 \mathrm{H}(\widehat P_\varepsilon|_{\mathcal F_N}
                    \Vert\widehat P_\varepsilon^{\boldsymbol{k}}|_{\mathcal F_N})
 =\sum_{j=1}^N\mathbb E_{\widehat P_\varepsilon}[h_j]
 \le\frac{1}{4m_\varepsilon}\mathbb E_{P_\varepsilon}\!\left[
       W_\varepsilon\sum_{j=1}^N\int_{U_j}\boldsymbol{J}(\gamma)^2
                          \mu_{\widetilde B_{2r}}^{\mathrm{loop}}(d\gamma)\right].
\end{equation}
As \(N\to\infty\), the relative entropies on the left increase
to \(\mathrm{H}(\widehat P_\varepsilon\Vert\widehat P_\varepsilon^{\boldsymbol{k}})\)
by~\cite[Theorem~21]{vEH14}. Multiplying by \(m_\varepsilon\),
and using monotone convergence and the disjointness of the query regions, we obtain
\begin{equation}\label{eq:query-entropy}
 \begin{aligned}
 \Delta_\varepsilon(r)
 \le{}&\frac14\mathbb E_{P_\varepsilon}\!\left[
       W_\varepsilon\int_{\Gamma_\varepsilon^{\le\delta r}}\boldsymbol{J}(\gamma)^2
       \mathbf{1}_{\smash[t]{\operatorname{Rev}(\gamma)}}
                          \,\mu_{\widetilde B_{2r}}^{\mathrm{loop}}(d\gamma)\right]\\
 ={}&\frac14\int_{\Gamma_\varepsilon^{\le\delta r}}\boldsymbol{J}(\gamma)^2
       \mathbb E_{P_\varepsilon}\!\left[
          W_\varepsilon\mathbf{1}_{\smash[t]{\operatorname{Rev}(\gamma)}}\right]
                          \,\mu_{\widetilde B_{2r}}^{\mathrm{loop}}(d\gamma)\\
 ={}&\frac14\int_{\substack{\gamma\in\Gamma_\varepsilon^{\le\delta r}\\\boldsymbol{J}(\gamma)>0}}
       \boldsymbol{J}(\gamma)^2R_\varepsilon(\gamma)
                          \,\mu_{\widetilde B_{2r}}^{\mathrm{loop}}(d\gamma).
 \end{aligned}
\end{equation}
Now by \eqref{eq:protected-killing-rate}, we have
\[
 \boldsymbol{J}(\gamma)
 =\kappa^2r^{-2}\sum_{|x|\ge r/8}\ell_x(\gamma)
 \le\kappa^2r^{-2}\sum_{x\in\mathbb Z^d}\ell_x(\gamma)
 =\kappa^2r^{-2}\mathcal V(\gamma).
\]
Substituting this bound into \eqref{eq:query-entropy} proves
\eqref{eq:exploration-entropy}.

Finally, set \(p=\widehat P_\varepsilon(A)\),
\(q=\widehat P_\varepsilon^{\boldsymbol{k}}(A)\), and
\(D=\mathrm{H}(\widehat P_\varepsilon\Vert\widehat P_\varepsilon^{\boldsymbol{k}})\).
If \(p>q\), it follows from~\cite[Lemma~2.3]{DM23} that
\begin{equation}
 p-q\le\sqrt{2pD}\le\frac p2+D,
 \qquad p\le2q+2D.
\end{equation}
The last bound also holds if \(p\le q\). Multiplying it by
\(m_\varepsilon\) and using
\(m_\varepsilon D=\Delta_\varepsilon(r)\) proves \eqref{eq:entropy-event}.
\end{proof}

\Needspace{11\baselineskip}
We next remove the duration cutoff \(\varepsilon\) by monotone convergence.

\begin{lemma}
\label{lem:query-transfer}
Let \(0<\delta<\frac18\), \(\delta r\ge L_0\), and
\(0<\varepsilon\le\varepsilon_0\), where \(L_0\) and
\(\varepsilon_0\) are as in Proposition~\ref{lem:query}.
Fix \(x\in\mathbb Z^d\) and an integer radius \(\rho\ge1\).
If a candidate loop \(\gamma\in\Gamma_\varepsilon^{\le\delta r}\) satisfies
\(\sup_{y\in\gamma}|y-x|\le\rho\), then
\begin{equation}\label{eq:query-probability}
 R_\varepsilon(\gamma)
 \le\mathbb E[|\phi_0|\mathbf{1}_{\smash[t]{\{0\xleftrightarrow{\mathcal L}B(x,2\rho)\}}}].
\end{equation}
If \(\Delta_\varepsilon(r)\le \Delta(r)\) uniformly for sufficiently small
\(\varepsilon>0\), then
\[
 w_d^{\le\delta r}(r)\le2w_d^{\boldsymbol{k}}(r)+2\Delta(r).
\]
\end{lemma}

\begin{proof}
We use the weight \(W_\varepsilon=\sqrt{2\widehat{\mathcal L}^{\le\delta r,\varepsilon}_{0}}\) in
\eqref{eq:query-weight-notation}, and fix a candidate loop
\(\gamma\) satisfying \(\sup_{y\in\gamma}|y-x|\le\rho\).
By Algorithm~\ref{exp:root-cluster}, we have \(\operatorname{Rev}(\gamma)
 \subseteq\{0\xleftrightarrow{
           \mathcal L^{\le\delta r,\varepsilon}}B(x,2\rho)\}\). Therefore,
by monotonicity, we have
\[
 R_\varepsilon(\gamma)
 =\mathbb E_{P_\varepsilon}\!\left[
       W_\varepsilon\mathbf{1}_{\smash[t]{\operatorname{Rev}(\gamma)}}\right]
 \le\mathbb E_{P_\varepsilon}\!\left[
       W_\varepsilon\mathbf{1}_{\smash[t]{\{0\xleftrightarrow{\mathcal L^{\le\delta r,\varepsilon}}B(x,2\rho)\}}}\right]
 \le\mathbb E[|\phi_0|\mathbf{1}_{\smash[t]{\{0\xleftrightarrow{\mathcal L}B(x,2\rho)\}}}].
\]
This proves \eqref{eq:query-probability}.
Applying \eqref{eq:entropy-event} to the event $\mathcal E_r^{\le\delta r,\varepsilon}
 :=\{0\xleftrightarrow{
                 \mathcal L^{\le\delta r,\varepsilon}}\partial\widetilde B_r\}$,
and using
\(\Delta_\varepsilon(r)\le \Delta(r)\), we obtain
\begin{equation}\label{eq:query-cutoff-event-comparison}
 \mathbb E_{P_\varepsilon}
       [W_\varepsilon\mathbf{1}_{\smash[t]{\mathcal E_r^{\le\delta r,\varepsilon}}}]
 \le2\mathbb E_{P_\varepsilon^{\boldsymbol{k}}}
       [W_\varepsilon\mathbf{1}_{\smash[t]{\mathcal E_r^{\le\delta r,\varepsilon}}}]
       +2\Delta_\varepsilon(r)
 \le2w_d^{\boldsymbol{k}}(r)+2\Delta(r).
\end{equation}
As
\(\varepsilon\downarrow0\), we have
\(W_\varepsilon\uparrow\sqrt{2\widehat{\mathcal L}^{\le\delta r}_{0}}\)
and \(\mathbf 1_{\mathcal E_r^{\le\delta r,\varepsilon}}
\uparrow\mathbf{1}_{\smash[t]{\{0\xleftrightarrow{
                 \mathcal L^{\le\delta r}}\partial\widetilde B_r\}}}\).
Monotone convergence in \eqref{eq:query-cutoff-event-comparison}
now gives
\begin{equation}
 w_d^{\le\delta r}(r)
 =\lim_{\varepsilon\downarrow0}
       \mathbb E_{P_\varepsilon}
       [W_\varepsilon\mathbf{1}_{\smash[t]{\mathcal E_r^{\le\delta r,\varepsilon}}}]
 \le2w_d^{\boldsymbol{k}}(r)+2\Delta(r).
\end{equation}
This completes the proof.
\end{proof}

\subsection{Upper bounds in dimensions three to five}

We first estimate the second moment of the total vertex local time under the loop measure.
\begin{lemma}\label{lem:layer}
For \(d\ge3\) and \(1\le s\le r\),
\begin{equation}\label{eq:diameter-moment}
 \int_{s<\operatorname{diam}\gamma\le2s}
       \mathcal V(\gamma)^2\,\mu_{\widetilde B_{2r}}^{\mathrm{loop}}(d\gamma)
 \lesssim r^d s^{4-d}.
\end{equation}
\end{lemma}

\begin{proof}
We use the occupation identity~\eqref{eq:loop-second-occupation},
recalled in Section~\ref{subsec:basic}.
For \(x\in B_{2r}\) and \(|x-y|\le s/6\),
the restriction property and~\cite[Section~2.2]{LJ11} give
\begin{equation}\label{eq:diameter-green-exit}
 \begin{aligned}
 \mu^{\mathrm{loop}}[\ell_x\ell_y;\operatorname{diam}\gamma>s]
 \le{}&\mu^{\mathrm{loop}}[\ell_x\ell_y;
\gamma\cap\partial\widetilde B(x,s/3)\ne\varnothing]\\
 ={}&G(x,y)^2-G_{\widetilde B(x,s/3)}(x,y)^2
 \lesssim s^{2-d}G(x,y),
 \end{aligned}
\end{equation}
where the last inequality uses the Green exit decomposition and~\cite[Proposition~4.6.2 and Theorem~4.3.1]{LL10}.
For fixed \(x\in B_{2r}\), the Green sum in
\eqref{eq:green-estimates} gives
\(\sum_{|y-x|\le s/6}s^{2-d}G(x,y)\lesssim s^{4-d}\).
For the remaining pairs, uniformly in \(x\in B_{2r}\), we have
\(\sum_{s/6<|y-x|\le2s}G(x,y)^2\lesssim s^d s^{4-2d}=s^{4-d}\).
We now expand
\(\mathcal V^2=\sum_{x,y}\ell_x\ell_y\). Bounding $\mu_{\widetilde B_{2r}}^{\mathrm{loop}}$ by 
$\mu^{\mathrm{loop}}$ and applying the preceding two bounds for each \(x\in B_{2r}\), we obtain
\[
 \begin{aligned}
 \int_{s<\operatorname{diam}\gamma\le2s}
       \mathcal V(\gamma)^2\,\mu_{\widetilde B_{2r}}^{\mathrm{loop}}(d\gamma)
 &\le\sum_{x\in B_{2r}}
       \sum_{|y-x|\le2s}
       \mu^{\mathrm{loop}}[\ell_x\ell_y;
                            s<\operatorname{diam}\gamma\le2s]\\
 &\lesssim\sum_{x\in B_{2r}}s^{4-d}
       \lesssim r^d s^{4-d}.
 \end{aligned}
\]
This proves \eqref{eq:diameter-moment}.
\end{proof}

\begin{lemma}[Renormalization inequality]\label{lem:below-six-renormalisation}
Let \(d\in\{3,4,5\}\) and \(\beta=(d-2)/2\).
There exist constants \(c,C>0\) and \(\delta_0\in(0,\frac1{16})\),
depending only on \(d\), such that the
following holds. Let \(\kappa\ge1\), \(0<\delta\le\delta_0\), and
let \(r\ge C\kappa\) be an integer with \(\delta r\ge L_0\),
where \(L_0\) is as in Proposition~\ref{lem:query}.
Then
\begin{equation}\label{eq:below-six-renormalisation}
 w_d(r)\le C_\delta r^{-\beta}
       +Ce^{-c\kappa}w_d(r/4)
       +C\kappa^4r^{-2}
       \sum_{j=0}^{\lfloor\log_2(\delta r)\rfloor}
              2^{2j}w_d(2^{j+1}).
\end{equation}
Here \(C_\delta\) depends only on \(d\) and \(\delta\).
\end{lemma}

\begin{proof}
Fix \(0<\varepsilon\le\varepsilon_0\), with \(\varepsilon_0\) as in
Proposition~\ref{lem:query}.
Let \(C_d\ge1\) be the constant in~Lemma~\ref{lem:small-ball},
and decrease \(\delta_0\), if necessary, so that
\(\delta_0\le(64C_d)^{-1}\).
Every candidate loop contributing to~\eqref{eq:exploration-entropy}
satisfies \(\operatorname{diam}\gamma\le\delta r\) and
\(\boldsymbol{J}(\gamma)>0\). For each such \(\gamma\),
choose a vertex \(z\in\gamma\) with \(\boldsymbol{k}(z)>0\)
and \(|z|\ge r/8\).
Since \(\operatorname{diam}\gamma\le\delta r<r/16\) and
\(r/8-\delta r>r/16\),
 the vertices of all contributing loops lie in the annulus
\(B_{2r}\setminus B_{r/16}\). We split the contributing loops according to their diameters.

\smallskip\noindent\textbf{Case 1.} \(\operatorname{diam}\gamma\le2\).
Then \(\sup_{y\in\gamma}|y-z|\le2\). Since \(\delta r\ge L_0\ge1\),
by our choice of \(\delta_0\), we have
\(4C_d\le4C_d\delta_0r\le r/16<|z|\).
Applying~\eqref{eq:query-probability} and~\eqref{eq:small-ball}, we obtain
\[
 R_\varepsilon(\gamma)\lesssim r^{2-d}w_d(4)
 \le r^{2-d}w_d(2),
\]
since \(w_d\) is nonincreasing. Moreover, by~\eqref{eq:loop-second-occupation}, we have
\[
 \int_{\operatorname{diam}\gamma\le2}\mathcal V(\gamma)^2
             \,\mu_{\widetilde B_{2r}}^{\mathrm{loop}}(d\gamma)
 \le\sum_{x\in B_{2r}}\sum_{|y-x|\le2}G(x,y)^2
 \lesssim r^d.
\]
Thus the contribution of this diameter range to the right-hand side
of~\eqref{eq:exploration-entropy} is
\(O(\kappa^4r^{-2}w_d(2))\).

\smallskip\noindent\textbf{Case 2.} \(2<\operatorname{diam}\gamma\le\delta r\).
Choose the unique integer \(j\ge1\) such that
\(2^j<\operatorname{diam}\gamma\le2^{j+1}\); in particular, \(2^j\le\delta r\).
Then \(\sup_{y\in\gamma}|y-z|\le2^{j+1}\), and
$
 C_d2^{j+2}\le4C_d\delta_0r\le r/16<|z|
$.
Applying~\eqref{eq:query-probability} and~\eqref{eq:small-ball}, and then
using \(w_d(2^{j+2})\le w_d(2^{j+1})\), we obtain
\begin{equation}\label{eq:dyadic-query}
 R_\varepsilon(\gamma)
 \lesssim (2^j/r)^{d-2}w_d(2^{j+1}).
\end{equation}
For this layer, we obtain from Lemma~\ref{lem:layer}
\[
 \int_{2^j<\operatorname{diam}\gamma\le2^{j+1}}\mathcal V(\gamma)^2
             \,\mu_{\widetilde B_{2r}}^{\mathrm{loop}}(d\gamma)
 \lesssim r^d2^{(4-d)j}.
\]
Hence its contribution to the right-hand side
of~\eqref{eq:exploration-entropy} is
\[
 O\!\left(\kappa^4r^{-4}\cdot r^d2^{(4-d)j}(2^j/r)^{d-2}w_d(2^{j+1})\right)
 =O\!\left(\kappa^4r^{-2}2^{2j}w_d(2^{j+1})\right).
\]
\smallskip
Combining the two cases in~\eqref{eq:exploration-entropy}, we obtain
\[
 \begin{aligned}
 \Delta_\varepsilon(r)
 \lesssim{}& \kappa^4r^{-2}\Bigl(
       w_d(2)+\sum_{j=1}^{\lfloor\log_2(\delta r)\rfloor}
       2^{2j}w_d(2^{j+1})\Bigr)\\
 ={}&\kappa^4r^{-2}
       \sum_{j=0}^{\lfloor\log_2(\delta r)\rfloor}
       2^{2j}w_d(2^{j+1}).
 \end{aligned}
\]
This bound is uniform for \(0<\varepsilon\le\varepsilon_0\), hence
Lemma~\ref{lem:query-transfer} applies. Combining it with
Lemmas~\ref{lem:large-loop-removal} and~\ref{lem:killed-contraction} proves
\eqref{eq:below-six-renormalisation}.
\end{proof}

\begin{proof}[Proof of Theorem~\ref{thm:main} for \(d=3,4,5\)]
Fix \(d\in\{3,4,5\}\) and set \(\beta=(d-2)/2\).
It suffices to show that there exists \(M>0\) such that, for every
\(k\ge0\),
\begin{equation}\label{eq:induction-hypothesis}
 w_d(2^k)\le M2^{-\beta k}.
\end{equation}
With \(C\) as in \eqref{eq:below-six-renormalisation}, set
\(C_0:=C\max\{4^\beta,(1-2^{\beta-2})^{-1}\}\).
Next choose \(\kappa\ge1\) such that \(C_0e^{-c\kappa}\le\frac14\), and then choose
\(0<\delta\le\delta_0\) such that
\(C_0\kappa^4\delta^{2-\beta}\le\frac14\).
Choose an integer \(k_0\ge2\) such that the scale conditions of
Lemma~\ref{lem:below-six-renormalisation} hold whenever \(r\ge2^{k_0}\).
Set \(M:=2C_\delta\vee\bigl(2^{\beta k_0}\sqrt{2G(0,0)/\pi}\bigr)\),
where \(C_\delta\) is the constant in~\eqref{eq:below-six-renormalisation}.
We prove~\eqref{eq:induction-hypothesis} by induction on \(k\ge0\).
Since \(w_d(2^k)\le\mathbb E|\phi_0|=\sqrt{2G(0,0)/\pi}\),
\eqref{eq:induction-hypothesis} holds for \(k\le k_0\) by definition of \(M\).
Now, fix \(k>k_0\) and assume that
\[
 w_d(2^j)\le M2^{-\beta j},\qquad j=0,\ldots,k-1.
\]
Applying Lemma~\ref{lem:below-six-renormalisation} with \(r=2^k\), we obtain
\[
 w_d(r)\le C_\delta r^{-\beta}+Ce^{-c\kappa}w_d(r/4)
 +C\kappa^4r^{-2}
 \sum_{j=0}^{\lfloor\log_2(\delta r)\rfloor}2^{2j}w_d(2^{j+1}).
\]
Here \(r/4=2^{k-2}\), and every index \(j\) in the sum satisfies
\(2^{j+1}\le2\delta r\le r/2=2^{k-1}\). Thus, by the induction
hypothesis, we have
\(w_d(r/4)\le M(r/4)^{-\beta}=4^\beta Mr^{-\beta}\).
Writing \(J=\lfloor\log_2(\delta r)\rfloor\), we also have
\[
 \sum_{j=0}^{J}2^{2j}w_d(2^{j+1})
 \le M\sum_{j=0}^{J}2^{(2-\beta)j}
 \le \frac{M}{1-2^{\beta-2}}(\delta r)^{2-\beta},
\]
where the last inequality uses \(2-\beta>0\) and \(2^J\le\delta r\).
Substituting these two bounds into \eqref{eq:below-six-renormalisation}
and using the definition of \(C_0\), we conclude that
\begin{equation}\label{eq:subcritical-recursion}
 \begin{aligned}
 r^\beta w_d(r)
 &\le C_\delta+CM\left(4^\beta e^{-c\kappa}
             +\frac{\kappa^4\delta^{2-\beta}}{1-2^{\beta-2}}\right)\\
 &\le C_\delta+C_0M\big(e^{-c\kappa}
                       +\kappa^4\delta^{2-\beta}\big)
 \le C_\delta+M/2\le M.
 \end{aligned}
\end{equation}
This proves~\eqref{eq:induction-hypothesis} at \(k\), completing the
induction. For arbitrary \(r\ge1\), choose \(k\ge0\) with
\(2^k\le r<2^{k+1}\). By monotonicity,
\(w_d(r)\le w_d(2^k)\le2^\beta Mr^{-\beta}\), as required.

Together with~\eqref{eq:weighted-to-one-arm}, the preceding estimate on \(w_d\) proves the
upper bound in~\eqref{eq:known-one-arm-estimates} for \(d=3,4,5\);
the lower bound follows from Lemma~\ref{lem:one-arm-lower}.
\end{proof}

We next bound the entropy cost using the triangle diagram in dimensions \(d\ge6\).

\subsection{The triangle diagram and dimensions above six}
\label{sec:above-six}

For a finite discrete set \(D\subset\mathbb Z^d\), regard
\(G_D=(G_D(x,y))_{x,y\in D}\)
as the Green matrix, with \((G_D^m)(x,y)\) denoting the \((x,y)\)-entry of
its \(m\)-th matrix power. Define the triangle diagrams
\[
 \nabla_D
 =\sup_{x\in D}\sum_{y,z\in D}
       G_D(x,y)G_D(y,z)G_D(z,x),\qquad \nabla_d
 =\sum_{y,z\in\mathbb Z^d}G(0,y)G(y,z)G(z,0).
\]
The latter is finite for \(d>6\); see e.g.~\cite[Section~5.4]{HH17}.
We first express a loop-occupation moment in terms of the Green matrix,
then bound it by the triangle diagram, including its logarithmic growth in dimension six.

\begin{lemma}
\label{lem:triangle-occupation}
Let \(\widetilde D\) be a bounded open box, and put \(D=\widetilde D\cap\mathbb Z^d\).
For every \(x\in D\),
\begin{equation}\label{eq:triangle-occupation-identity}
 \int \mathcal V(\gamma)^2\mathbf{1}_{\smash[t]{\{x\in\gamma\}}}
                   \,\mu_{\widetilde D}^{\mathrm{loop}}(d\gamma) =
 \frac{2(G_D^3)(x,x)}{G_D(x,x)}
 -\left(\frac{(G_D^2)(x,x)}{G_D(x,x)}\right)^2
 \le 2\nabla_D.
\end{equation}
For every integer \(r\ge2\) and finite discrete set
\(D\subset B_{2r}\),
\begin{equation}\label{eq:high-finite-triangle}
 \nabla_D\le\nabla_d<\infty\quad(d>6),\qquad
 \nabla_D\lesssim\log r\quad(d=6).
\end{equation}
\end{lemma}

\begin{proof}
Writing \(G_{D,t}=(G_D^{-1}+tI)^{-1}\),
the restriction property and~\cite[Section~2, Proposition~7]{LJ11} give
\begin{equation}\label{eq:high-restricted-loop-moment}
 F_x(t):=\mu_{\widetilde D}^{\mathrm{loop}}\!\left[
       (e^{-t\mathcal V}-1)\mathbf{1}_{\smash[t]{\{x\in\gamma\}}}\right]
 =\log\frac{G_{D,t}(x,x)}{G_D(x,x)}.
\end{equation}
Since \(G_{D,t}'=-G_{D,t}^2\) and
\(G_{D,t}''=2G_{D,t}^3\), differentiating twice for
\(t>0\) in \eqref{eq:high-restricted-loop-moment} and then letting \(t\downarrow0\), by monotone convergence, we have
\[
 \int\mathcal V(\gamma)^2\mathbf{1}_{\smash[t]{\{x\in\gamma\}}}
             \,\mu_{\widetilde D}^{\mathrm{loop}}(d\gamma)
 =F_x''(0)
 =\frac{2(G_D^3)(x,x)}{G_D(x,x)}
  -\left(\frac{(G_D^2)(x,x)}{G_D(x,x)}\right)^2.
\]
Since
\(G_D(x,x)\ge1\) and \((G_D^3)(x,x)\le\nabla_D\), this proves \eqref{eq:triangle-occupation-identity}.

For \(d>6\), by \(G_D\le G\), we have
\(\nabla_D\le\nabla_d\).
Write \(\langle z\rangle=1+|z|\). By \(G_D\le G\),
the Green-function bound \(G(x,y)\lesssim \langle x-y\rangle^{2-d}\)
from~\cite[Theorem~4.3.1]{LL10}, and the standard convolution estimate
\begin{equation}\label{eq:high-green-convolution}
 \sum_{y\in\mathbb Z^d}G(x,y)G(y,z)
 \lesssim \langle z-x\rangle^{4-d},\qquad d>4,
\end{equation}
we have, for \(d=6\), \(D\subset B_{2r}\), and \(x\in D\),
\begin{equation}\label{eq:high-triangle-spatial-sum}
 (G_D^3)(x,x)
 \lesssim \sum_{z\in D}\langle z-x\rangle^{-6}
 \lesssim \log r.
\end{equation}
This completes the proof of \eqref{eq:high-finite-triangle}.
\end{proof}

The preceding loop-moment estimate bounds the entropy error for \(d\ge6\).
\begin{lemma}
\label{lem:high-triangle-entropy}
For every \(d\ge6\), there exists \(\delta_0(d)>0\) such that the
following holds. 
Let \(L_0\) and \(\varepsilon_0\) be as in Proposition~\ref{lem:query}.
For \(\kappa\ge1\), \(0<\delta\le\delta_0(d)\),
an integer \(r\ge C_d\kappa\) with \(\delta r\ge L_0\),
and \(0<\varepsilon\le\varepsilon_0\),
\begin{equation}\label{eq:high-triangle-entropy}
 \Delta_\varepsilon(r)\lesssim \kappa^4r^{-2}\nabla_{B_{2r}}.
\end{equation}
\end{lemma}

\begin{proof}
Fix a candidate loop \(\gamma\in\Gamma_\varepsilon^{\le\delta r}\) with
\(\boldsymbol{J}(\gamma)>0\). If it is revealed, then some vertex \(x\in\gamma\cap\mathbb Z^d\)
has a nearest neighbor \(y\) connected to \(0\) by \(\mathcal L\).
By the definition of \(R_\varepsilon\)
and the weighted two-point identity~\eqref{eq:weighted-two-point}, we obtain
\begin{equation}\label{eq:high-vertex-revealment}
 R_\varepsilon(\gamma)
 \lesssim \sum_{x\in\gamma\cap\mathbb Z^d}
       \sum_{y\sim x}
       \mathbb E[|\phi_0|\mathbf{1}_{\smash[t]{\{0\xleftrightarrow{\mathcal L}y\}}}] \lesssim \sum_{x\in\gamma\cap\mathbb Z^d}
       G(0,x).
\end{equation}
Here we used the fact that \(0\notin\gamma\), and
\(\sum_{y\sim x}G(0,y)=2d\,G(0,x)\).
Substitute~\eqref{eq:high-vertex-revealment}
into~\eqref{eq:exploration-entropy}. Dropping the diameter and duration restrictions,
and using Lemma~\ref{lem:triangle-occupation} and domain
monotonicity, we obtain
\[
 \begin{aligned}
 \Delta_\varepsilon(r)
 \lesssim{}& \kappa^4r^{-4}\sum_{x\in B_{2r}}
       G(0,x)
       \int\mathcal V(\gamma)^2\mathbf{1}_{\smash[t]{\{x\in\gamma\}}}
                         \,\mu_{\widetilde B_{2r}}^{\mathrm{loop}}(d\gamma)\\
 \lesssim{}& \kappa^4r^{-4}\nabla_{B_{2r}}
       \sum_{x\in B_{2r}}G(0,x)
 \lesssim \kappa^4r^{-2}\nabla_{B_{2r}}.
 \end{aligned}
\]
The last line follows from~\eqref{eq:green-estimates} and completes the proof of \eqref{eq:high-triangle-entropy}.
\end{proof}

\Needspace{4\baselineskip}
We record the following simple result, which is useful in the induction.
\begin{lemma}\label{lem:contraction-iteration}
Let \(w:[1,\infty)\to[0,\infty)\) be bounded and nonincreasing,
and let \(h:[1,\infty)\to[1,\infty)\) be nondecreasing.
Suppose that there exist an integer \(r_0\ge1\), a constant \(C<\infty\), and
\(0<\rho<1\) such that
\[
 r^2w(r)\le\rho(r/4)^2w(r/4)+Ch(r)
\]
for every integer \(r\ge4r_0\). Then there is a constant
\(C'<\infty\) such that \(w(r)\le C'r^{-2}h(r)\) for every
\(r\ge1\).
\end{lemma}

\begin{proof}
Set \(r_n=4^nr_0\). Iterating the assumed inequality and using the monotonicity
of \(h\), we obtain
\[
 r_n^2w(r_n)
 \le\rho^n r_0^2w(r_0)
       +C\sum_{j=0}^{n-1}\rho^jh(r_{n-j}) \le\left(r_0^2w(r_0)+\frac{C}{1-\rho}\right)h(r_n).
\]
If \(r_n\le r<r_{n+1}\), then by monotonicity of \(w\) and \(h\),
we have \(w(r)\le w(r_n)\le16C'r^{-2}h(r)\).
Boundedness of \(w\) and \(h\ge1\) cover \(1\le r<r_0\)
after increasing the constant. This completes the proof.
\end{proof}

\begin{proof}[Proof of Theorem~\ref{thm:main} for \(d>6\)]
Fix \(d>6\). For the upper bound, we combine
Lemmas~\ref{lem:large-loop-removal}, \ref{lem:killed-contraction},
and~\ref{lem:query-transfer} with~\eqref{eq:high-triangle-entropy}
to obtain the renormalization inequality
\begin{equation}\label{eq:high-renormalisation}
 w_d(r)\le Ce^{-c\kappa}w_d(r/4)
       +C_\delta r^{-(d-2)/2}
       +C\kappa^4r^{-2}\nabla_d.
\end{equation}
Choose \(\kappa\ge1\) so that \(16Ce^{-c\kappa}\le\frac12\).
Set \(\delta=\delta_0(d)\wedge1/16\) and
\(r_*=\lceil\max\{4,C_d\kappa,L_0/\delta\}\rceil\).
Here \(\delta_0(d)\), \(L_0\), and \(C_d\) are as in
Lemma~\ref{lem:high-triangle-entropy}, with \(C_d\) enlarged, if necessary,
so that the radius condition in Lemma~\ref{lem:killed-contraction}
also holds.
Since \((d-2)/2>2\), multiplying
\eqref{eq:high-renormalisation} by \(r^2\), we obtain
\[
 r^2w_d(r)\le\tfrac12(r/4)^2w_d(r/4)+C,
 \qquad r\ge r_*.
\]
The function \(w_d\) is nonincreasing and bounded by
\(\mathbb E|\phi_0|\). Applying Lemma~\ref{lem:contraction-iteration}
with \(h=1\) and \(\rho=\frac12\), we obtain \(w_d(r)\lesssim r^{-2}\).
The comparison~\eqref{eq:weighted-to-one-arm} proves the upper bound for \(\theta_d(r)\).

The lower bound follows from Lemma~\ref{lem:one-arm-lower}.
\end{proof}

\begin{remark}[The triangle diagram in dimension six]
\label{rem:triangle-six}
The same argument already gives
\begin{equation}\label{eq:six-elementary-triangle-upper}
 \theta_6(r)\lesssim r^{-2}\log r,\qquad r\ge2.
\end{equation}
Indeed, by~\eqref{eq:high-finite-triangle} and
\eqref{eq:high-triangle-entropy}, we have
\(\Delta_\varepsilon(r)\lesssim \kappa^4r^{-2}\log r\).
The large-loop error is now \(O_\delta(r^{-2})\), so with the
same choices of \(\kappa\) and \(\delta\), we obtain
\[
 r^2w_6(r)\le\tfrac12(r/4)^2w_6(r/4)+C\log r.
\]
By Lemma~\ref{lem:contraction-iteration}, with \(h(r)=\log(er)\)
and \(\rho=\frac12\), we obtain \(w_6(r)\lesssim r^{-2}\log r\) for
\(r\ge2\), and
\eqref{eq:six-elementary-triangle-upper} follows. 
\end{remark}

Thus, finiteness of the
full-space triangle diagram closes the renormalization inequality
above dimension six, whereas its finite-volume growth produces a
logarithmic loss at dimension six, as in~\cite{EGPS25}.
In the next section, we use the additional capacity estimates below to sharpen the logarithmic factor
and prove the sharp bound in~\eqref{eq:main-six}.

\section{One-arm bounds in dimension six}\label{sec:capacity}

Throughout the remainder of the paper, we fix \(d=6\).
The first subsection proves the lower bound in~\eqref{eq:main-six}.
The remaining subsections establish the upper bound using the entropy
estimate from Section~\ref{sec:exploration}. Both arguments use the
weighted escape estimate in Proposition~\ref{s6:escape-bound}, whose
proof is deferred to Section~\ref{sec:escape-proof}.

\subsection{The one-arm lower bound}\label{sec:escape}

We first recall the following isomorphisms, which are fundamental in our analysis.
For an open domain \(\widetilde D\subset\widetilde{\mathbb Z}^d\),
write \(\widetilde p_{\widetilde D}\) for the transition density of the Brownian motion
\(Z\) killed on leaving \(\widetilde D\), and
\(\widetilde{\mathbb P}_{\widetilde D,t}^{x,y}\) for its bridge law
from \(x\) to \(y\) in time \(t>0\).
For \(x,y\in\widetilde D\), define the path measure
\cite[Definition~3.1]{Lup18} by
\begin{equation}\label{eq:integrated-bridge-measure}
 \mu_{\widetilde D}^{x,y}(d\gamma)
   =\int_0^\infty \widetilde p_{\widetilde D}(t;x,y)\,
                       \widetilde{\mathbb P}_{\widetilde D,t}^{x,y}(d\gamma)\,dt.
\end{equation}
Its total mass is \(G_{\widetilde D}(x,y)\).
We write
\(\ell(\gamma)\) for the occupation field of $\gamma$.
For every bounded nonnegative measurable functional \(F\) of the
full occupation field on \(\widetilde D\),
\begin{equation}\label{eq:bridge}
 \mathbb E_{\widetilde D}\!\left[|\phi_{\widetilde D,x}\phi_{\widetilde D,y}|
       \mathbf{1}_{\smash[t]{\{x\xleftrightarrow{\mathcal L_{\widetilde D}}y\}}}
       F(\widehat{\mathcal L})\right]
 =\int\mu_{\widetilde D}^{x,y}(d\gamma)\,
       \mathbb E_{\widetilde D}[F(\widehat{\mathcal L}+\ell(\gamma))].
\end{equation}
This identity follows from the BFS--Dynkin isomorphism~\cite{BFS82,Dyn84a,Dyn84b,LJ11}
and the coupling in
Lemma~\ref{thm:lupu-isomorphism}. A discrete formulation is
given in \cite[Theorem~1.4]{Lup22}, and the metric-graph formulation can be found in~\cite[Remark~1.13]{JLQ23fields}.
Applying \eqref{eq:bridge} to
\(F(\widehat{\mathcal L})
  ((2\widehat{\mathcal L}_x)^{-\frac12}\wedge M)\)
and letting \(M\uparrow\infty\),
we obtain by monotone convergence 
\begin{equation}\label{eq:one-sided-bridge}
 \mathbb E_{\widetilde D}[|\phi_{\widetilde D,y}| F(\widehat{\mathcal L});
                x\xleftrightarrow{\mathcal L_{\widetilde D}}y]
 =\int\mu_{\widetilde D}^{x,y}(d\gamma)\,
   \mathbb E_{\widetilde D}\!\left[
     F(\widehat{\mathcal L}+\ell(\gamma))
     \bigl[2(\widehat{\mathcal L}_x+\ell_x(\gamma))\bigr]^{-\frac12}
   \right].
\end{equation}
This is the one-point weighted version of \eqref{eq:bridge}.
Both \eqref{eq:bridge} and \eqref{eq:one-sided-bridge} extend to nonnegative measurable
\(F\) and to full space.

\medskip

We now define the occupation field and the cluster containing \(0\) after adding
an independent stopped Brownian path to the loop soup.
For an integer
\(R\ge1\), let \(\mathcal L_{\widetilde B_R}\) be the critical
loop soup in \(\widetilde B_R\). Independently, let \(Z\) be the 
Brownian motion on \(\widetilde{\mathbb Z}^6\) started from \(0\), stopped at
\(\zeta_{\widetilde B_R}=\inf\{t\ge0:Z_t\in\partial\widetilde B_R\}\).
Write \(\widehat{\mathcal L}_{\widetilde B_R,z}\) for the loop soup occupation field
at \(z\), and \(L_t^z(Z)\) for the local time of \(Z\) at \(z\)
up to time \(t\). The total occupation field of the loop soup and the path is
\begin{equation}\label{eq:six-enhanced-occupation}
 \Gamma_R(z)=\widehat{\mathcal L}_{\widetilde B_R,z}
                  +L^z_{\zeta_{\widetilde B_R}}(Z),\qquad z\in\widetilde{\mathbb Z}^6.
\end{equation}
Let \(\mathcal K_R\) be the discrete set of the component of
\(\{\Gamma_R>0\}\) containing zero,
\begin{equation}\label{eq:six-enhanced-root-cluster}
 \mathcal K_R=\{x\in B_R:0\xleftrightarrow{\{\Gamma_R>0\}}x\}.
\end{equation}
Equivalently, take the union of the Brownian path
\(Z([0,\zeta_{\widetilde B_R}))\) and all metric-graph clusters of
\(\mathcal L_{\widetilde B_R}\) that intersect this path;
then \(\mathcal K_R\) is the intersection of this union with \(\mathbb Z^6\).

For a simple random walk \(Y\) from \(0\), independent of the loop soup
and \(Z\), the equilibrium mass at \(0\) 
\begin{equation}\label{eq:six-equilibrium-weight}
 e_{\mathcal K_R}(0)=\mathbb P(Y_n\notin \mathcal K_R\text{ for all }n\ge1\mid \mathcal K_R)
\end{equation}
is the escape probability conditional on \(\mathcal K_R\).
In accordance with \eqref{eq:one-sided-bridge}, define the weighted
escape probability by
\begin{equation}\label{eq:six-b-definition}
 b(R)=\mathbb E\!\left[
       e_{\mathcal K_R}(0)(2\Gamma_R(0))^{-\frac12}\right],
\end{equation}
where the expectation is taken over the independent loop soup \(\mathcal L_{\widetilde B_R}\) and the Brownian motion \(Z\).
For a noninteger argument \(R\ge1\), we use
\(b(R)=b(\lfloor R\rfloor)\), and apply the same convention to
\(\mathcal K_R\) and \(\Gamma_R\).

We first record a uniform bound on \(b(R)\).
By Lemma~\ref{thm:lupu-isomorphism}, we see that
\(\widehat{\mathcal L}_{\widetilde B_R,0}\sim\operatorname{Gamma}(\frac12,G_{\widetilde B_R}(0,0))\),
where the second parameter is the scale. Independently,
\(L^0_{\zeta_{\widetilde B_R}}(Z)\) is exponential with mean \(G_{\widetilde B_R}(0,0)\). 
Consequently,
\(\Gamma_R(0)\sim\operatorname{Gamma}(\frac32,G_{\widetilde B_R}(0,0))\).
Since \(G_{\widetilde B_R}(0,0)\ge1\), we have
\begin{equation}\label{eq:six-inverse-amplitude-moment}
 \mathbb E[(2\Gamma_R(0))^{-1}]=G_{\widetilde B_R}(0,0)^{-1}\le1.
\end{equation}
Cauchy--Schwarz and \(0\le e_{\mathcal K_R}(0)\le1\) then give
\begin{equation}\label{eq:six-b-uniform-bound}
 b(R)^2\le\mathbb E[e_{\mathcal K_R}(0)^2]
       \le\mathbb E[e_{\mathcal K_R}(0)]\le1.
\end{equation}
Note that both
\(\mathbb E[e_{\mathcal K_R}(0)]\) and \(b(R)\)
are nonincreasing in \(R\).

The one-sided bridge identity~\eqref{eq:one-sided-bridge} now bounds a
weighted equilibrium mass in terms of \(b\).

\begin{lemma}\label{s6:local-bridge}
Let \(W\subset\mathbb Z^6\) be deterministic, let \(\rho\ge1\) be an
integer, and suppose that \(B(x,\rho)\subset W\) and
\(|x|\ge8\rho\). Then
\begin{equation}\label{eq:six-local-bridge-bound}
 \mathbb E\!\left[|\phi_0|e_{\mathcal C(0)\cap W}(x);
                   0\xleftrightarrow{\mathcal L}x\right]
 \lesssim G(0,x)b(\rho).
\end{equation}
\end{lemma}

\begin{proof}
For a nonnegative function \(f\) on
\(\widetilde{\mathbb Z}^6\), let \(\mathcal C_f(x)\) be the discrete
set of the component of \(\{f>0\}\) containing \(x\), as in \eqref{eq:six-enhanced-root-cluster}, with
\(\mathcal C_f(x)=\{x\}\) if \(f(x)=0\). Take
\(F(f)=e_{\mathcal C_f(x)\cap W}(x)\).
On
\(\{0\xleftrightarrow{\mathcal L}x\}\), we have
\(\mathcal C_{\widehat{\mathcal L}}(x)=\mathcal C(0)\), and
\(F(\widehat{\mathcal L})=e_{\mathcal C(0)\cap W}(x)\).
For \(\gamma\) sampled from the path measure \(\mu_{\widetilde{\mathbb Z}^6}^{x,0}\), set
\(\Gamma^\gamma=\widehat{\mathcal L}+\ell(\gamma)\) and
\(K^\gamma=\mathcal C_{\Gamma^\gamma}(x)\).
Applying \eqref{eq:one-sided-bridge} in full space with
\(y=0\) and the functional \(F\), we have
\begin{equation}\label{eq:six-local-bridge-insertion}
 \mathcal I:=\mathbb E\!\left[|\phi_0|e_{\mathcal C(0)\cap W}(x);
                     0\xleftrightarrow{\mathcal L}x\right]
 =\int\mu_{\widetilde{\mathbb Z}^6}^{x,0}(d\gamma)\,
       \mathbb E\!\left[
       \frac{e_{K^\gamma\cap W}(x)}{\sqrt{2\Gamma^\gamma(x)}}\right].
\end{equation}

Let \(\sigma\) denote the first hitting time of \(\partial\widetilde B(x,\rho)\).
For every
nonnegative functional \(h\), by the Markov property, we have
\begin{equation}\label{eq:six-local-bridge-stopping}
 \int\mu_{\widetilde{\mathbb Z}^6}^{x,0}(d\gamma)\,h(\gamma_{[0,\sigma]})
 =\mathbb E_x\!\left[h(Z_{[0,\sigma]})G(Z_\sigma,0)\right],
\end{equation}
where \(Z\) on the right is ordinary Brownian motion from \(x\).
The assumption \(|x|\ge8\rho\) and \eqref{eq:green-estimates} imply
\(G(Z_\sigma,0)\lesssim G(x,0)\). 
Restrict the loop soup to \(\widetilde B(x,\rho)\) and retain only the stopped
path \(\gamma_{[0,\sigma]}\). The resulting cluster is contained in \(K^\gamma\cap W\),
and its occupation at \(x\) is at most \(\Gamma^\gamma(x)\).
Both the equilibrium mass and the inverse square-root weight therefore
increase under this restriction. Applying~\eqref{eq:six-local-bridge-stopping}
and translation invariance to the restricted loop soup and path, we obtain
\(\mathcal I\lesssim G(0,x)b(\rho)\).
This completes the proof.
\end{proof}

For \(r\ge1\), recall the truncated vertex cluster of \(0\)
and define its weighted mean capacity by
\begin{equation}\label{eq:six-local-capacity-target}
 \mathcal C_r(0):=\mathcal C(0)\cap B_r,\qquad
 M(r):=\mathbb E[|\phi_0|\operatorname{cap}(\mathcal C_r(0))].
\end{equation}
\begin{proposition}\label{s6:mean-capacity-ball}
For all sufficiently large \(r\),
\begin{equation}\label{eq:six-mean-capacity-bound}
 M(r)\lesssim r^2b(\sqrt r)+r^{\frac32}.
\end{equation}
\end{proposition}

\begin{proof}
To estimate \(M(r)\), set \(\rho=\lfloor\sqrt r\rfloor\).
By \eqref{eq:capacity}, we can write 
\[
 \operatorname{cap}(\mathcal C_r(0))
 =\sum_{x\in B_r}e_{\mathcal C_r(0)}(x)
                    \mathbf{1}_{\smash[t]{\{0\xleftrightarrow{\mathcal L}x\}}}.
\]
For every \(x\in B_{r-2\rho}\setminus B_{8\rho}\), we have
\(B(x,\rho)\subset B_r\) and \(|x|\ge8\rho\).
By \eqref{eq:six-local-bridge-bound}
and the Green estimate \eqref{eq:green-estimates}, the contribution
of these vertices to \(M(r)\) satisfies
\begin{equation}\label{eq:six-mean-capacity-interior}
 \sum_{x\in B_{r-2\rho}\setminus B_{8\rho}}
 \mathbb E\!\left[|\phi_0|e_{\mathcal C_r(0)}(x);
                  0\xleftrightarrow{\mathcal L}x\right]
 \lesssim b(\rho)\sum_{x\in B_r}G(0,x)
 \lesssim r^2b(\rho).
\end{equation}
The remaining vertices lie in a
neighborhood of the root or in the boundary layer of \(B_r\).
Using \(e\le1\) and the weighted two-point identity
\eqref{eq:weighted-two-point}, we have
\[
 \sum_{x\in B_{8\rho}}G(0,x)\lesssim\rho^2,
 \qquad
 \sum_{x\in B_r\setminus B_{r-2\rho}}G(0,x)
 \lesssim r^5\rho\,r^{-4}\lesssim r\rho.
\]
Consequently,
\[
 M(r)\lesssim r^2b(\rho)+\rho^2+r\rho
       \lesssim r^2b(\rho)+r^{\frac32}.
\]
This proves \eqref{eq:six-mean-capacity-bound}.
\end{proof}

\Needspace{4\baselineskip}
The following estimate for \(b(R)\), proved in
Section~\ref{sec:escape-proof}, will be used throughout this section.

\begin{proposition}\label{s6:escape-bound}
For every integer \(R\ge2\),
\(b(R)\lesssim(\log R)^{-1/2}\).
\end{proposition}

\Needspace{3\baselineskip}
We can now deduce the one-arm lower bound.

\begin{proof}[Proof of the lower bound in \eqref{eq:main-six}]
Combining Propositions~\ref{s6:mean-capacity-ball} and~\ref{s6:escape-bound}, we get
\begin{equation}\label{eq:Mr}
  M(r)=\mathbb E[|\phi_0|\operatorname{cap}(\mathcal C_r(0))]\lesssim r^2(\log r)^{-1/2}.
\end{equation}
By Lemma~\ref{lem:capacity-arm-transfer}, we then have
\(\theta_6(r)\gtrsim M(r)^{-1}
\gtrsim r^{-2} (\log r)^{1/2}\).
\end{proof}

\subsection{Loop-duration estimates}

Let \(X\) be a continuous-time simple random walk on \(\mathbb Z^6\) with jump rate one.
For a finite set \(V\subset\mathbb Z^6\), let \(p_t^V\) be the transition kernel of \(X\) killed on leaving \(V\)
and extended by zero outside \(V\).
Write \(H_V=\inf\{u\ge0:X_u\in V\}\).
Recall that \(\mathcal V(\gamma)=\sum_x\ell_x(\gamma)\) is the total
vertex local time of $\gamma$. For an integer \(r\ge2\), define
\[
 \mathcal I_{2r}(V)=\int \mathcal V(\gamma)^2
       \mathbf 1_{\{\gamma\cap V\ne\varnothing\}}
       \,\mu_{\widetilde B_{2r}}^{\mathrm{loop}}(d\gamma).
\]

\begin{lemma}\label{s6:duration-trace}
For every integer \(r\ge2\) and finite \(V\subset\mathbb Z^6\),
\[
 \mathcal I_{2r}(V)\le\int_0^\infty tT_V(t)\,dt,
 \qquad
 T_V(t):=\sum_{x\in B_{2r}}\bigl[
 p_t^{B_{2r}}(x,x)-p_t^{B_{2r}\setminus V}(x,x)\bigr].
\]
Moreover, there is a constant \(c>0\) such that, for every \(t>0\),
\begin{equation}\label{eq:six-loop-hitting-heat-bound}
 T_V(t)\lesssim(1+t)^{-3}e^{-ct/r^2}
                  \sum_x\mathbb P_x(H_V\le t/2).
\end{equation}
\end{lemma}

\begin{proof}
By \cite[Section~2]{Lup16} and
\cite[Section~7.2]{LJ11}, we have
\[
 \mathcal I_{2r}(V)
 \le\int_0^\infty\frac{t^2}{t}
       \sum_{x\in B_{2r}}
       \mathbb P_x(X_t=x,\ t<H_{B_{2r}^c},\ H_V\le t)\,dt.
\]
For each \(x\in B_{2r}\), subtracting the paths that avoid
\(V\), we obtain
\[
 \mathbb P_x(X_t=x,\ t<H_{B_{2r}^c},\ H_V\le t)
 =p_t^{B_{2r}}(x,x)-p_t^{B_{2r}\setminus V}(x,x).
\]
Summing over \(x\) and integrating proves the first bound.
In particular, \(T_V(t)\) is the total mass of closed paths of
duration \(t\) that stay in \(B_{2r}\) and visit \(V\).

To prove \eqref{eq:six-loop-hitting-heat-bound}, we use the heat-kernel
and box eigenvalue bounds in \cite[Sections~2.1 and 8.2]{LL10} to get
\[
 p_t^{B_{2r}}(x,y)\lesssim(1+t)^{-3}e^{-ct/r^2},
 \qquad t\ge0,\quad x,y\in B_{2r}.
\]
Every closed path that visits \(V\) has a visit in its first half
in at least one of its two time orientations. Time reversal and
the strong Markov property at the first such visit imply
\[
 T_V(t)\le2\sum_{x\in B_{2r}}\mathbb E_x\!\left[
 \mathbf 1_{\{H_V\le t/2,\,H_V<H_{B_{2r}^c}\}}
 p_{t-H_V}^{B_{2r}}(X_{H_V},x)\right].
\]
The remaining return time is at least \(t/2\). Applying the
heat-kernel bound and removing the condition
\(H_V<H_{B_{2r}^c}\) proves
\eqref{eq:six-loop-hitting-heat-bound}.
\end{proof}

Next, we bound the right-hand side of \eqref{eq:six-loop-hitting-heat-bound} in terms of \(\operatorname{cap}(V)\).

\begin{lemma}\label{s6:finite-time-hitting}
If \(V\) is
contained in a discrete box of side length \(s\ge1\), then
\begin{equation}\label{eq:six-finite-time-hitting}
 \sum_z\mathbb P_z(H_V\le t)
 \lesssim t\operatorname{cap}(V),\qquad t\ge s^2.
\end{equation}
\end{lemma}

\begin{proof}
We may assume that \(V\ne\varnothing\). Fix \(t\ge s^2\) and
write \(p_u(a,b)=\mathbb P_a(X_u=b)\).
Since \(|a-b|\le s\) for \(a,b\in V\), the continuous-time
local central limit theorem \cite[Section~2.5.2]{LL10} and
\eqref{eq:green-estimates} imply
\begin{equation}\label{eq:six-finite-time-green}
 \int_0^t p_u(a,b)\,du\gtrsim G(a,b),\qquad a,b\in V.
\end{equation}
Using \eqref{eq:six-finite-time-green} and the equilibrium identity
\eqref{eq:capacity-hitting-estimates}, we obtain, for every \(a\in V\),
\begin{equation}\label{eq:six-equilibrium-occupation}
 F(a):=\mathbb E_a\!\left[\int_0^t e_V(X_u)\,du\right]
 =\sum_{b\in V}e_V(b)\int_0^t p_u(a,b)\,du
 \gtrsim\sum_{b\in V}G(a,b)e_V(b)=1.
\end{equation}
On \(\{H_V\le t\}\), we have \(X_{H_V}\in V\) and
\([H_V,H_V+t]\subset[0,2t]\). By the strong Markov
property at \(H_V\) and \eqref{eq:six-equilibrium-occupation}, we get
\begin{equation}\label{eq:six-hitting-occupation}
 \begin{aligned}
 \mathbb E_z\!\left[\int_0^{2t}e_V(X_u)\,du\right]
 &\ge\mathbb E_z\!\left[
       \int_{H_V}^{H_V+t}e_V(X_u)\,du;\ H_V\le t\right]\\
 &=\mathbb E_z[F(X_{H_V});\ H_V\le t]
 \gtrsim\mathbb P_z(H_V\le t).
 \end{aligned}
\end{equation}
By symmetry, we have
\(\sum_zp_u(z,b)=\sum_zp_u(b,z)=1\).
Summing \eqref{eq:six-hitting-occupation} over \(z\), we obtain
\begin{equation}\label{eq:six-hitting-sum-bound}
 \sum_z\mathbb P_z(H_V\le t)
 \lesssim\int_0^{2t}\sum_{b\in V}e_V(b)
                    \sum_zp_u(z,b)\,du
 =2t\sum_{b\in V}e_V(b)
 =2t\operatorname{cap}(V).
\end{equation}
This concludes the proof.
\end{proof}

For an integer \(s\ge1\), let \(\mathcal Q_s\) be the partition
of \(\mathbb Z^6\) into the vertex boxes
\(sz+\{0,\ldots,s-1\}^6\), \(z\in\mathbb Z^6\). Applying
Lemma~\ref{s6:finite-time-hitting} to each piece of a finite set \(V\)
and taking a union bound, we obtain
\begin{equation}\label{eq:six-box-capacity}
 \sum_z\mathbb P_z(H_V\le t)
 \le\sum_{Q\in\mathcal Q_s}\sum_z
                 \mathbb P_z(H_{V\cap Q}\le t)
 \lesssim t\sum_{Q\in\mathcal Q_s}
                 \operatorname{cap}(V\cap Q),\qquad t\ge s^2.
\end{equation}
We choose the spatial partition according to the time scale \(t\),
matching the box size to the diffusive scale of the random walk.
This yields the following multiscale capacity bound.

\begin{lemma}\label{s6:duration}
For every integer \(r\ge2\) and \(V\subset B_{2r}\),
\begin{equation}\label{eq:six-duration-capacity-sum}
 \mathcal I_{2r}(V)
 \lesssim\sum_{j=0}^{\lceil\log_2r\rceil}
              \sum_{Q\in\mathcal Q_{2^j}}\operatorname{cap}(V\cap Q).
\end{equation}
\end{lemma}

\begin{proof}
By the first bound in Lemma~\ref{s6:duration-trace}, it suffices to estimate
\(\int_0^\infty tT_V(t)\,dt\).
For \(t\ge s^2\ge1\), use
\(\mathbb P_x(H_V\le t/2)\le\mathbb P_x(H_V\le t)\) in
\eqref{eq:six-loop-hitting-heat-bound}, then apply~\eqref{eq:six-box-capacity}
at time \(t\). We obtain
\begin{equation}\label{eq:six-duration-integrand-bound}
 tT_V(t)\lesssim t^{-1}e^{-ct/r^2}
                     \sum_{Q\in\mathcal Q_s}\operatorname{cap}(V\cap Q).
\end{equation}

Set \(N=\lceil\log_2r\rceil\). For \(0\le j<N\) and $t\in [4^j,4^{j+1}]$, 
set \(s=2^j\).
Then \(t\ge s^2\), and
\eqref{eq:six-duration-integrand-bound} applies with the same partition
\(\mathcal Q_{2^j}\) on this interval. Using \(e^{-ct/r^2}\le1\),
we obtain
\[
 \int_{4^j}^{4^{j+1}}tT_V(t)\,dt
 \lesssim\sum_{Q\in\mathcal Q_{2^j}}\operatorname{cap}(V\cap Q)
               \int_{4^j}^{4^{j+1}}\frac{dt}{t}
 =(\log4)\sum_{Q\in\mathcal Q_{2^j}}\operatorname{cap}(V\cap Q).
\]
Summing over
\(0\le j<N\) bounds the integral over \([1,4^N]\).
For \(t\ge4^N\), use the fixed scale \(s=2^N\). Then,
\[
 \int_{4^N}^{\infty}tT_V(t)\,dt
 \lesssim\sum_{Q\in\mathcal Q_{2^N}}\operatorname{cap}(V\cap Q)
               \int_{4^N}^{\infty}e^{-ct/r^2}\frac{dt}{t}
 \lesssim\sum_{Q\in\mathcal Q_{2^N}}\operatorname{cap}(V\cap Q).
\]

It remains to treat \(0<t<1\). By the strong Markov property,
on \(\{H_V\le t\}\) the conditional probability of no further
jump before time \(t\) is \(e^{-(t-H_V)}\). Hence
\[
 \mathbb P_x(X_t\in V)
 \ge\mathbb E_x[e^{-(t-H_V)};\ H_V\le t]
 \ge e^{-t}\mathbb P_x(H_V\le t).
\]
Summing over the starting vertex and using symmetry
of the heat kernel, we have
\[
 \sum_x\mathbb P_x(H_V\le t)
 \le e^t\sum_x\mathbb P_x(X_t\in V)
 =e^t\sum_{y\in V}\sum_xp_t(x,y)
 =e^t|V|\le e|V|.
\]
It follows from \eqref{eq:six-loop-hitting-heat-bound} that
\(\int_0^1 tT_V(t)\,dt\lesssim|V|\). The partition \(\mathcal Q_1\)
consists of singletons, and
\[
 \sum_{Q\in\mathcal Q_1}\operatorname{cap}(V\cap Q)
 =|V|/G(0,0).
\]
This is the \(j=0\) term of the desired bound and completes the proof.
\end{proof}

We next prove a local version of the weighted capacity estimate in
Proposition~\ref{s6:mean-capacity-ball}.
We write
$
 \widehat{\mathcal C}_r(0)
 =\mathcal C(0)\cap(B_{2r}\setminus B_{r/16})
$ for the outer part of the cluster $\mathcal C(0)$.

\begin{proposition}\label{s6:mean-capacity}
For all sufficiently large \(r\) and integers \(1\le s\le2r\),
\begin{equation}\label{eq:six-local-weighted-capacity-log}
 \mathbb E\!\left[|\phi_0|
       \sum_{Q\in\mathcal Q_s}\operatorname{cap}(\widehat{\mathcal C}_r(0)\cap Q)
       \right]
 \lesssim r^2\{(\log(s+2))^{-\frac12}+s^{-\frac12}\}.
\end{equation}
\end{proposition}

\begin{proof}
Set \(\rho=\lfloor\sqrt s\rfloor\). For each box
\(Q\in\mathcal Q_s\), write \(W=Q\cap(B_{2r}\setminus B_{r/16})\), and
expand the capacity into equilibrium masses as in the proof of
Proposition~\ref{s6:mean-capacity-ball}.
A vertex \(x\in W\) is an interior vertex if \(B(x,\rho)\subset W\), and a
boundary vertex otherwise. At an interior vertex, the
condition
\(|x|\ge8\rho\) holds for large \(r\). Therefore, as in \eqref{eq:six-mean-capacity-interior},
for every interior vertex  $x$,
\begin{equation}\label{eq:six-interior-equilibrium-mass}
\mathbb E\!\left[|\phi_0|
       e_{\widehat{\mathcal C}_r(0)\cap Q}(x);
       0\xleftrightarrow{\mathcal L}x\right]
 \lesssim G(0,x)b(\rho).
\end{equation}
Summing over the interior vertices gives a contribution bounded by a constant multiple of \(r^2b(\rho)\).

There are at most $C(r^6\rho/s+r^5\rho)$ boundary vertices. 
At each boundary vertex, the bound \(e\le1\) and the
weighted two-point identity \eqref{eq:weighted-two-point} control
the weighted equilibrium mass. Since \(G(0,x)\asymp r^{-4}\) on
the annulus $B_{2r}\setminus B_{r/16}$, the total boundary contribution is at most
\[
 Cr^{-4}(r^6\rho/s+r^5\rho)
 =Cr^2(\rho/s+\rho/r).
\]
Combining the interior and boundary estimates, we obtain
\[
 \mathbb E\!\left[|\phi_0|
       \sum_{Q\in\mathcal Q_s}\operatorname{cap}(\widehat{\mathcal C}_r(0)\cap Q)
       \right]
 \lesssim r^2\{b(\rho)+\rho/s+\rho/r\}.
\]
Now \(\rho/s\le s^{-\frac12}\) and
\(\rho/r\le2s^{-\frac12}\). For sufficiently large \(s\),
Proposition~\ref{s6:escape-bound}, with
\(\log\rho\asymp\log(s+2)\), proves
\eqref{eq:six-local-weighted-capacity-log}. For bounded \(s\),
the same bound follows directly from
\(\sum_Q\operatorname{cap}(\widehat{\mathcal C}_r(0)\cap Q)
\le|\widehat{\mathcal C}_r(0)|\) and
\eqref{eq:weighted-two-point}, after adjusting the constant.
\end{proof}

\subsection{Entropy estimate and upper bound}\label{s6:upper}

We first use Proposition~\ref{s6:mean-capacity} to derive an entropy bound.

\begin{lemma}\label{s6:entropy}
There exist constants \(C,r_0<\infty\) and
\(\delta_0\in(0,\frac1{32})\) such that the following holds.
Let \(\kappa\ge1\), \(0<\delta\le\delta_0\), and let
\(r\ge r_0\) be an integer with \(\delta r\ge L_0\).
For the protected killing \eqref{eq:protected-killing-rate} and
Algorithm~\ref{exp:root-cluster} in
\(\widetilde B_{2r}\),
\begin{equation}\label{eq:six-entropy-log-bound}
 \Delta_\varepsilon(r)
 \le C\kappa^4r^{-2}(\log r)^{\frac12},
 \qquad0<\varepsilon\le\varepsilon_0.
\end{equation}
Here \(L_0\) and \(\varepsilon_0\) are as in
Proposition~\ref{lem:query}, and \(C\) is independent of
\(\kappa,\delta,r,\varepsilon\).
\end{lemma}

\begin{proof}
As in the proof of Lemma~\ref{lem:below-six-renormalisation}, every
contributing candidate \(\gamma\) visits a vertex \(z\) with
\(|z|\ge r/8\). Since \(\operatorname{diam}\gamma\le\delta r\)
and \(\delta\le1/32\), all its vertices and their nearest
neighbors lie outside \(B_{r/16}\) for sufficiently large \(r\).
Let $\widehat{\mathcal C}_r^{+}(0)$ consist of $\widehat{\mathcal C}_r(0)$ and 
all its nearest neighbors.
Since \(\sqrt{2\widehat{\mathcal L}^{\le\delta r,\varepsilon}_{0}}\le|\phi_0|\),
for every \(\gamma\in\Gamma_\varepsilon^{\le\delta r}\) with
\(\boldsymbol J(\gamma)>0\), 
\begin{equation}\label{eq:six-query-probability-capacity-charge}
 R_\varepsilon(\gamma)
 \le\mathbb E[|\phi_0|
          \mathbf 1_{\{\gamma\cap\widehat{\mathcal C}_r^{+}(0)\ne\varnothing\}}].
\end{equation}
Substituting~\eqref{eq:six-query-probability-capacity-charge} into
\eqref{eq:exploration-entropy} and applying Tonelli's theorem, we obtain
\begin{equation}\label{eq:six-entropy-capacity-charge}
 \Delta_\varepsilon(r)\lesssim\kappa^4r^{-4}
     \mathbb E[|\phi_0|\mathcal I_{2r}(\widehat{\mathcal C}_r^{+}(0))].
\end{equation}
Applying Lemma~\ref{s6:duration} to
\(\widehat{\mathcal C}_r^+(0)\cap B_{2r}\), then using \eqref{eq:neighbor-capacity-comparison}
and Proposition~\ref{s6:mean-capacity}, we obtain
\[
 \begin{aligned}
 \Delta_\varepsilon(r)
 &\lesssim\kappa^4r^{-4}
       \sum_{j=0}^{\lceil\log_2r\rceil}
       \mathbb E\!\left[|\phi_0|
       \sum_{Q\in\mathcal Q_{2^j}}
       \operatorname{cap}(\widehat{\mathcal C}_r(0)\cap Q)\right]\\
 &\lesssim\kappa^4r^{-2}
       \sum_{j=0}^{\lceil\log_2r\rceil}
       \{(1+j)^{-\frac12}+2^{-j/2}\}\lesssim\kappa^4r^{-2}(\log r)^{\frac12}.
 \end{aligned}
\]
This finishes the proof.
\end{proof}

\begin{proof}[Proof of the upper bound in \eqref{eq:main-six}]
Let \(\kappa\ge1\) and \(0<\delta\le\delta_0\).
Enlarge \(r_0\) from Lemma~\ref{s6:entropy} so that
\(r_0\ge\max\{3,L_0/\delta,C\kappa\}\), where \(C\) is as in
Lemma~\ref{lem:killed-contraction}. For \(r\ge r_0\),
Lemmas~\ref{lem:query-transfer} and~\ref{s6:entropy} imply
\begin{equation}\label{eq:six-short-soup-comparison}
 w_6^{\le\delta r}(r)
 \le2w_6^{\boldsymbol k}(r)
       +C\kappa^4r^{-2}(\log r)^{\frac12}.
\end{equation}
Combining \eqref{eq:six-short-soup-comparison} with Lemmas~\ref{lem:large-loop-removal}
and~\ref{lem:killed-contraction}, we obtain, for \(r\ge r_0\),
\begin{equation}\label{eq:six-weighted-arm-recursion}
 w_6(r)\le Ce^{-c\kappa}w_6(r/4)+C_\delta r^{-2}
       +C\kappa^4r^{-2}(\log r)^{\frac12}.
\end{equation}

Now fix \(\kappa\) large enough that
\(16Ce^{-c\kappa}\le\frac12\), and fix any
\(0<\delta\le\delta_0\). Set \(h(r)=\sqrt{\log r}\) for
\(r\ge r_0\) and extend \(h\) constantly below \(r_0\).
After multiplication by \(r^2\), the recursion becomes
\[
 r^2w_6(r)\le\tfrac12(r/4)^2w_6(r/4)+Ch(r),\qquad r\ge r_0.
\]
Since \(w_6\) is nonincreasing and bounded by
\(\mathbb E[|\phi_0|]\), we apply Lemma~\ref{lem:contraction-iteration}
to obtain
\begin{equation}\label{eq:six-weighted-arm-upper}
 w_6(r)\lesssim r^{-2}(\log r)^{\frac12}.
\end{equation}
The FKG inequality, in the form
\eqref{eq:weighted-to-one-arm}, proves the upper bound in
\eqref{eq:main-six}.
\end{proof}

\begin{proof}[Proof of Theorem~\ref{thm:main}]
Combining this upper bound with the lower bound of
Subsection~\ref{sec:escape} proves \eqref{eq:main-six}, and hence the
up-to-constants estimate \eqref{eq:six-one-arm-asymptotic}.
Together with Section~\ref{sec:exploration}, this completes the proof
of Theorem~\ref{thm:main}.
\end{proof}

\section{Weighted escape probability}\label{sec:escape-proof}
 
In this section, we prove Proposition~\ref{s6:escape-bound}, namely 
\begin{equation}\label{eq:bR}
 b(R):=\mathbb E\!\left[
       e_{\mathcal K_R}(0)(2\Gamma_R(0))^{-\frac12}\right]\lesssim (\log R)^{-\frac12}.
\end{equation}
Recall that \(\Gamma_R\) is the occupation field of
\(\mathcal L_{\widetilde B_R}\) enhanced by an independent Brownian motion
\(Z\) from \(0\), stopped on leaving \(\widetilde B_R\), and that
\(\mathcal K_R\) is the discrete set of the resulting cluster containing \(0\).

Note that $\mathbb E[e_{\mathcal K_R}(0)]$ can be viewed as the (generalized) non-intersection probability of 
two independent random walks inside an independent critical loop soup. We refer to \cite{GNQ24,GLL26} for
a similar setup in two and three dimensions.
The idea for bounding \(b(R)\) is to test a geometrically separated family of scales from
$\sqrt{R}$ to $R$ to see if there is a common cluster hit by both random walks,
and the non-intersection event requires all of these tests to fail.
In Subsection~\ref{s6:annular-bridge}, we will employ a second-moment method to control the success probability 
at any intermediate scale \(s\) by
estimating the moments of the number of clusters hit by both walks.
Expressing all the estimates in terms of $b$, we end up with a renormalization inequality for 
$b$, which allows us to get the desired upper bound \eqref{eq:bR} by strong induction.

\subsection{Common-cluster counts}\label{s6:two-endpoint-comparison}\label{s6:annular-bridge}

Consider the open test box
\(\widetilde D_s=\{z\in\widetilde{\mathbb Z}^6:|z-se_1|<s/8\}\)
and set \(D_s=\widetilde D_s\cap\mathbb Z^6\).
Put
\[
 \mathcal T_s(u,v)=\sum_{z,w\in D_s}G(u,z)G(z,w)G(w,v).
\]
For each \(u\in D_s\), define its vertex cluster in \(\widetilde D_s\) by
\begin{equation}\label{eq:six-vertex-cluster}
 \mathcal C_{\widetilde D_s}(u)
 =\{u\}\cup\{v\in D_s:
       u\xleftrightarrow{\mathcal L_{\widetilde D_s}}v\}.
\end{equation}
Let \(\mathfrak C_s=\{\mathcal C_{\widetilde D_s}(u):u\in D_s\}\)
be the family of distinct such sets. These sets partition \(D_s\),
with isolated vertices included as singletons. Define the expected
sum of squared cluster capacities by
\begin{equation}\label{eq:six-common-cluster-count}
 T(s)=\mathbb E_{\widetilde D_s}\sum_{\mathcal C\in\mathfrak C_s}\operatorname{cap}(\mathcal C)^2.
\end{equation}
We also define the corresponding quantity weighted by \(\mathcal T_s(u,v)\),
\begin{equation}\label{eq:six-equilibrium-triangle-definition}
 Q(s)=\mathbb E_{\widetilde D_s}\sum_{\mathcal C\in\mathfrak C_s}
              \sum_{u,v\in\mathcal C}e_{\mathcal C}(u)e_{\mathcal C}(v)\mathcal T_s(u,v).
\end{equation}
Let \(X,Y\) be independent simple random walks from
\(a,b\in\partial B_{s/4}\), each stopped
on leaving \(B_{4s}\), and independent of \(\mathcal L_{\widetilde D_s}\).
Let \(N(s)\) count all distinct clusters in \(\mathfrak C_s\)
visited by both walks. The following lemma bounds the first moment and the second moment of \(N(s)\) in terms of \(T(s)\) and \(Q(s)\).

\begin{lemma}\label{s6:second-moment}
For all sufficiently large $s$ and \(a,b\in\partial B_{s/4}\),
\begin{equation}\label{eq:six-count-moments}
 \mathbb E N(s)\asymp s^{-8}T(s),\qquad
 \mathbb E[N(s)(N(s)-1)]\lesssim s^{-8}Q(s).
\end{equation}
\end{lemma}

\begin{proof}
For any nonempty deterministic \(A\subset D_s\), by last exit and
Harnack comparison~\cite[Propositions~4.6.4, 6.3.5, 6.5.1 and~6.5.4]{LL10}, we have
\begin{equation}\label{eq:six-cluster-entrance-estimates}
 \mathbb P_a(H_A<H_{B_{4s}^c})\asymp s^{-4}\operatorname{cap}(A),\qquad
 \mathbb P_a(X_{H_A}=u,H_A<H_{B_{4s}^c})\le C s^{-4}e_A(u).
\end{equation}
Applying the first estimate of \eqref{eq:six-cluster-entrance-estimates}
to each \(\mathcal C\in\mathfrak C_s\) and using the independence of $X$ and $Y$, we have
\begin{equation}\label{eq:six-common-cluster-first-moment}
 \begin{aligned}
 \mathbb E N(s)
 &=\mathbb E_{\widetilde D_s}\sum_{\mathcal C\in\mathfrak C_s}
       \mathbb P_a(H_{\mathcal C}<H_{B_{4s}^c})\mathbb P_b(H_{\mathcal C}<H_{B_{4s}^c})\\
 &\asymp s^{-8}\mathbb E_{\widetilde D_s}
       \sum_{\mathcal C\in\mathfrak C_s}\operatorname{cap}(\mathcal C)^2
   =s^{-8}T(s).
 \end{aligned}
\end{equation}

We next turn to the second moment.
Let \(L_X(z)\) be the number of visits by \(X\) to \(z\) before
it exits \(B_{4s}\).
For \(z\in D_s\), by the Markov property
and \eqref{eq:six-cluster-entrance-estimates}, we have
\begin{equation}\label{eq:six-marked-occupation}
 \begin{aligned}
 V_A^a(z)
 &:=\mathbb E_a\!\left[L_X(z)\mathbf{1}_{\{H_A<H_{B_{4s}^c}\}}\right]\\
 &=G_{B_{4s}\setminus A}(a,z)\mathbb P_z(H_A<H_{B_{4s}^c})
   +\sum_{u\in A}\mathbb P_a(X_{H_A}=u,H_A<H_{B_{4s}^c})G_{B_{4s}}(u,z)\\
 &\lesssim s^{-4}\sum_{u\in A}e_A(u)G(u,z).
 \end{aligned}
\end{equation}
Here we also used \(G_{B_{4s}\setminus A}(a,z)\lesssim s^{-4}\)
and the last-exit formula. The same bound holds with \(b\) in place of \(a\).

Write \(I_{\mathcal C}\) for the indicator that both walks $X$ and $Y$ visit a cluster
\(\mathcal C\in\mathfrak C_s\). Then
\(I_{\mathcal C}\le\sum_{z,w\in\mathcal C}L_X(z)L_Y(w)\).
Write $\mathcal C_v:=\mathcal C_{\widetilde D_s}(v)$ (see \eqref{eq:six-vertex-cluster}).
On the event $I_{\mathcal C_v}=1$, we have
\begin{equation}\label{eq:six-other-cluster-count}
 N(s)-1
 =\sum_{\substack{\mathcal C\in\mathfrak C_s\\\mathcal C\ne \mathcal C_v}}I_{\mathcal C}
 \le\sum_{\substack{\mathcal C\in\mathfrak C_s\\\mathcal C\ne \mathcal C_v}}
       \sum_{z,w\in\mathcal C}L_X(z)L_Y(w)
 =\sum_{z,w\in D_s\setminus \mathcal C_v}
       L_X(z)L_Y(w)
       \mathbf{1}_{\smash[t]{\{z\xleftrightarrow{\mathcal L_{\widetilde D_s}}w\}}}.
\end{equation}
Combining \eqref{eq:six-other-cluster-count}
and
\eqref{eq:six-marked-occupation}, we obtain
\begin{equation}\label{eq:six-conditional-cluster-moment}
 \begin{aligned}
 \mathbb E[I_{\mathcal C_v}(N(s)-1)\mid \mathcal C_v]
 &\le\sum_{z,w\in D_s\setminus \mathcal C_v}
       V_{\mathcal C_v}^a(z)V_{\mathcal C_v}^b(w)
       \mathbb P\!\left(z\xleftrightarrow{\mathcal L_{\widetilde D_s}}w
                        \mid \mathcal C_v\right)\\
 &\lesssim\sum_{z,w\in D_s}
       V_{\mathcal C_v}^a(z)V_{\mathcal C_v}^b(w)G(z,w)\\
 &\lesssim s^{-8}\sum_{u,u'\in \mathcal C_v}e_{\mathcal C_v}(u)e_{\mathcal C_v}(u')
       \sum_{z,w\in D_s}G(u,z)G(z,w)G(w,u')\\
 &= s^{-8}\sum_{u,u'\in \mathcal C_v}
       e_{\mathcal C_v}(u)e_{\mathcal C_v}(u')\mathcal T_s(u,u').
 \end{aligned}
\end{equation}
Multiplying \eqref{eq:six-conditional-cluster-moment} by \(|\mathcal C_v|^{-1}\),
taking expectations, and summing over \(v\in D_s\), we obtain
\begin{equation}\label{eq:six-common-cluster-factorial-moment}
 \begin{aligned}
 \mathbb E[N(s)(N(s)-1)]
 &=\mathbb E\!\left[\sum_{\mathcal C\in\mathfrak C_s}I_{\mathcal C}(N(s)-1)\right]
 =\sum_{v\in D_s}\mathbb E_{\widetilde D_s}\!\left[
       \frac{\mathbb E[I_{\mathcal C_v}(N(s)-1)\mid \mathcal C_v]}{|\mathcal C_v|}\right]\\
 &\lesssim s^{-8}\mathbb E_{\widetilde D_s}\!\left[
       \sum_{v\in D_s}\frac{1}{|\mathcal C_v|}
       \sum_{u,u'\in \mathcal C_v}e_{\mathcal C_v}(u)e_{\mathcal C_v}(u')\mathcal T_s(u,u')\right]\\
 &= s^{-8}\mathbb E_{\widetilde D_s}\!\left[
       \sum_{\mathcal C\in\mathfrak C_s}\sum_{u,u'\in\mathcal C}
       e_{\mathcal C}(u)e_{\mathcal C}(u')\mathcal T_s(u,u')\right]
   =s^{-8}Q(s).
 \end{aligned}
\end{equation}
This completes the proof.
\end{proof}

We next bound \(Q(s)\) in terms of the
weighted escape probability \(b\) defined in \eqref{eq:six-b-definition}.
Note that we can rewrite $Q(s)$ as
\begin{equation}\label{eq:six-equilibrium-triangle-pair-sum}
      Q(s)=\sum_{u,v\in D_s}\mathcal T_s(u,v)
       \mathbb E_{\widetilde D_s}\!\left[
       e_{\mathcal C_{\widetilde D_s}(u)}(u)
       e_{\mathcal C_{\widetilde D_s}(u)}(v);
       u\xleftrightarrow{\mathcal L_{\widetilde D_s}}v\right].
\end{equation}
Therefore, it suffices to prove a two-endpoint version of
Lemma~\ref{s6:local-bridge}. We state it for a general domain and later
apply it with \(\widetilde D=\widetilde D_s\).

\begin{lemma}\label{s6:two-endpoint-equilibrium}
Let \(\widetilde D\subset\widetilde{\mathbb Z}^6\) and 
\(u,v\in\widetilde D\cap\mathbb Z^6\). Let \(\rho\ge1\) be an integer.
Suppose that
\(\widetilde B(u,\rho),\widetilde B(v,\rho)\subset\widetilde D\)
and \(|u-v|\ge8\rho\).
Then
\[
 \mathbb E_{\widetilde D}\!\left[e_{\mathcal C_{\widetilde D}(u)}(u)e_{\mathcal C_{\widetilde D}(u)}(v);
                        u\xleftrightarrow{\mathcal L_{\widetilde D}}v\right]
 \lesssim G(u,v)b(\rho)^2.
\]
Here \(\mathcal C_{\widetilde D}(u)\) is defined as in
\eqref{eq:six-vertex-cluster}, with \(\widetilde D\) in place of \(\widetilde D_s\).
\end{lemma}

\begin{proof}
The proof is similar to that of Lemma~\ref{s6:local-bridge}, and we adopt the notation introduced there.
Set
\(\Gamma_{\widetilde D}^{\gamma}
 =\widehat{\mathcal L}_{\widetilde D}+\ell(\gamma)\) and
\(K_{\widetilde D}^{\gamma}
 =\mathcal C_{\Gamma_{\widetilde D}^{\gamma}}(u)\).
As in \eqref{eq:six-local-bridge-insertion}, we have
\begin{equation}\label{eq:six-two-endpoint-bridge-insertion}
 \mathbb E_{\widetilde D}[e_{\mathcal C_{\widetilde D}(u)}(u)e_{\mathcal C_{\widetilde D}(u)}(v);
                   u\xleftrightarrow{\mathcal L_{\widetilde D}}v]
 =\int\mu_{\widetilde D}^{u,v}(d\gamma)\,
 \mathbb E_{\widetilde D}\!\left[
 \frac{e_{K_{\widetilde D}^{\gamma}}(u)e_{K_{\widetilde D}^{\gamma}}(v)}
 {2\sqrt{\Gamma_{\widetilde D}^{\gamma}(u)
         \Gamma_{\widetilde D}^{\gamma}(v)}}\right].
\end{equation}
Here \(e_{K_{\widetilde D}^{\gamma}}(v)=0\) if
\(v\notin K_{\widetilde D}^{\gamma}\).

Retain the initial bridge segment up to its first exit from \(\widetilde B(u,\rho)\),
and the reversed terminal segment up to its first exit from \(\widetilde B(v,\rho)\),
with exit vertices \(a\) and \(c\). As in the proof of
Lemma~\ref{s6:local-bridge}, the Markov property and time reversal
show that their joint measure has density \(G_{\widetilde D}(a,c)\) relative
to the product law of two independent stopped Brownian motions.
Moreover, $ G_{\widetilde D}(a,c)\lesssim G(u,v)$.
Now restrict the loop soup to \(\widetilde B(u,\rho)\) and \(\widetilde B(v,\rho)\) 
and retain only these stopped
path segments. By monotonicity, the
independence under the product law, and translation invariance, 
the right-hand side of~\eqref{eq:six-two-endpoint-bridge-insertion}
is bounded by a constant multiple of \(G(u,v)b(\rho)^2\). The result follows.
\end{proof}

We can now bound \(Q(s)\) from above.
\begin{lemma}\label{s6:annealed-triangle}
For all sufficiently large $s$, we have 
\begin{equation}\label{eq:six-annealed-triangle}
 Q(s)\lesssim s^8b(\sqrt s)^2+s^{15/2}.
\end{equation}
\end{lemma}

\begin{proof}
First, by the Green convolution estimate~\eqref{eq:high-green-convolution}, we have
\begin{equation}\label{eq:six-triangle-log-bound}
 \mathcal T_s(u,v)\lesssim \log\left(2+\frac{s}{1+|u-v|}\right).
\end{equation}
Set \(\rho=\lfloor\sqrt s\rfloor\). Let \(A_s=\{u\in D_s:\operatorname{dist}(u,D_s^c)<2\rho\}\)
be the boundary layer. Partition the ordered pairs in \(D_s^2\)
into the disjoint sets
\[
 \begin{aligned}
 \mathcal P_s^{\mathrm{far}}
   &=\{(u,v)\in(D_s\setminus A_s)^2:|u-v|\ge8\rho\},\\
 \mathcal P_s^{\mathrm{bd}}
   &=\{(u,v)\in D_s^2:u\in A_s\text{ or }v\in A_s\},\\
 \mathcal P_s^{\mathrm{near}}
   &=\{(u,v)\in(D_s\setminus A_s)^2:|u-v|<8\rho\}.
 \end{aligned}
\]
Write \(q_s(u,v)\) for the expectation in
\eqref{eq:six-equilibrium-triangle-pair-sum}.
By Lemma~\ref{s6:two-endpoint-equilibrium}, we have
\begin{equation}\label{eq:six-triangle-far-pair-bound}
 q_s(u,v)\lesssim G(u,v)b(\rho)^2,
       \qquad (u,v)\in\mathcal P_s^{\mathrm{far}}.
\end{equation}
Using \(e_{\mathcal C}\le1\) and the two-point estimate
\eqref{eq:two-point}, including on the diagonal, we also obtain
\begin{equation}\label{eq:six-triangle-general-pair-bound}
 q_s(u,v)\le\mathbb P_{\widetilde D_s}
       (u\xleftrightarrow{\mathcal L_{\widetilde D_s}}v)
       \lesssim G(u,v),
       \qquad u,v\in D_s.
\end{equation}
By \eqref{eq:green-estimates}, uniformly in \(u\in D_s\),
\begin{equation}\label{eq:six-triangle-row-sum}
 \sum_{v\in D_s}G(u,v)
       \log\left(2+\frac{s}{1+|u-v|}\right)
 \lesssim s^2\sum_{k\ge0}4^{-k}(k+1)
 \lesssim s^2.
\end{equation}
Combining~\eqref{eq:six-triangle-log-bound},
\eqref{eq:six-triangle-far-pair-bound}, and~\eqref{eq:six-triangle-row-sum},
we obtain the following bound on the contribution of \(\mathcal P_s^{\mathrm{far}}\),
\[
 \sum_{(u,v)\in\mathcal P_s^{\mathrm{far}}}
       q_s(u,v)\mathcal T_s(u,v)
 \lesssim |D_s|s^2b(\rho)^2
 \lesssim s^8b(\rho)^2.
\]
Using \(|A_s|\lesssim s^5\rho\),
\eqref{eq:six-triangle-general-pair-bound}, and
\eqref{eq:six-triangle-row-sum}, we obtain
\[
 \sum_{(u,v)\in\mathcal P_s^{\mathrm{bd}}}
       q_s(u,v)\mathcal T_s(u,v)
   \lesssim |A_s|s^2\lesssim s^7\rho.
\]
For the close interior pairs, by~\eqref{eq:six-triangle-log-bound} and
\eqref{eq:six-triangle-general-pair-bound}, we have
\[
 \begin{aligned}
 \sum_{(u,v)\in\mathcal P_s^{\mathrm{near}}}
       q_s(u,v)\mathcal T_s(u,v)
 &\lesssim s^6\sum_{z\in B_{8\rho}}G(0,z)
       \log\left(2+\frac{s}{1+|z|}\right)\\
 &\lesssim s^6\rho^2\left(1+\log\frac{s}{\rho}\right)
 \lesssim s^7\log s.
 \end{aligned}
\]
Finally, summing the three contributions and using
\(\rho\le\sqrt s\) and \(\log s\lesssim\sqrt s\), we obtain
\[
 Q(s)\lesssim s^8b(\rho)^2+s^7\rho+s^7\log s
       \lesssim s^8b(\sqrt s)^2+s^{15/2}.
\]
This proves \eqref{eq:six-annealed-triangle}.
\end{proof}

Recall that $T(s)=\mathbb E_{\widetilde D_s}\sum_{\mathcal C\in\mathfrak C_s}\operatorname{cap}(\mathcal C)^2$ in \eqref{eq:six-common-cluster-count}.
We seek a lower bound on \(T(s)\) in terms of \(b(s/4)\).
To this end, we weight each cluster capacity by the sum of the absolute field values over that cluster.
For \(\mathcal C\in\mathfrak C_s\), write $ m(\mathcal C):=\sum_{y\in\mathcal C}|\phi_{\widetilde D_s,y}|$,
and define
\begin{equation}\label{eq:six-amplitude-sum}
\mathcal M(s)=\mathbb E_{\widetilde D_s}
       \sum_{\mathcal C\in\mathfrak C_s}\operatorname{cap}(\mathcal C)m(\mathcal C).
\end{equation}

\begin{lemma}\label{s6:capacity-amplitude-lower}
For all sufficiently large \(s\), we have
$
 T(s)\gtrsim s^{-8}\mathcal M(s)^2
$.
\end{lemma}

\begin{proof}
Conditional on \(\mathcal L_{\widetilde D_s}\), the signs of clusters are independent
by Lemma~\ref{thm:lupu-isomorphism}. Hence,
\[
 \mathbb E_{\widetilde D_s}[\phi_{\widetilde D_s,y}\phi_{\widetilde D_s,z}\mid\mathcal L_{\widetilde D_s}]
 =|\phi_{\widetilde D_s,y}\phi_{\widetilde D_s,z}|
       \mathbf{1}_{\smash[t]{\{y\xleftrightarrow{\mathcal L_{\widetilde D_s}}z\}}}.
\]
Expanding the square cluster by cluster and taking expectations gives
\begin{equation}\label{eq:six-amplitude-mass-moment}
 \begin{aligned}
 \mathbb E_{\widetilde D_s}\sum_{\mathcal C\in\mathfrak C_s}
       m(\mathcal C)^2
 &=\sum_{y,z\in D_s}\mathbb E_{\widetilde D_s}\!\left[
       |\phi_{\widetilde D_s,y}\phi_{\widetilde D_s,z}|
       \mathbf{1}_{\smash[t]{\{y\xleftrightarrow{\mathcal L_{\widetilde D_s}}z\}}}\right]\\
 &=\sum_{y,z\in D_s}\mathbb E_{\widetilde D_s}[\phi_{\widetilde D_s,y}\phi_{\widetilde D_s,z}]
  =\sum_{y,z\in D_s}G_{\widetilde D_s}(y,z).
 \end{aligned}
\end{equation}
The bound~\eqref{eq:green-estimates} then gives
\[
 \sum_{y,z\in D_s}G_{\widetilde D_s}(y,z)
 \le\sum_{y\in D_s}\sum_{z\in B(y,s/4)}G(y,z)
 \lesssim s^8.
\]
By the Cauchy--Schwarz inequality and~\eqref{eq:six-amplitude-mass-moment}, we have
\[
 \mathcal M(s)
 \le\left(\mathbb E_{\widetilde D_s}\sum_{\mathcal C\in\mathfrak C_s}
                    \operatorname{cap}(\mathcal C)^2\right)^{1/2}
       \left(\mathbb E_{\widetilde D_s}\sum_{\mathcal C\in\mathfrak C_s}
                    m(\mathcal C)^2\right)^{1/2}
 \lesssim s^4 T(s)^{1/2}.
\]
Rearranging gives the result.
\end{proof}

We next combine the one-sided bridge identity with a last-exit
decomposition to bound \(\mathcal M(s)\) and \(T(s)\) from below
in terms of \(b(s/4)\).

\begin{lemma}\label{s6:annular-capacity-lower}
There exists a constant \(c>0\) such that, for every sufficiently large
\(s\),
\begin{equation}\label{eq:six-annular-lower-bounds}
 \mathcal M(s)\ge c s^8b(s/4),\qquad
 T(s)\ge c s^8b(s/4)^2.
\end{equation}
\end{lemma}

\begin{proof}
We adopt the notation in the proof of Lemma~\ref{s6:local-bridge}.
For a path \(\gamma\) in \(\widetilde D_s\) starting at \(x\in D_s\), 
let \(\Gamma^\gamma\) and \(\mathcal C_{\Gamma^\gamma}(x)\) be the occupation field
and vertex cluster containing \(x\) obtained by adding \(\gamma\) to the loop soup. 
Set $F_x(\gamma)=e_{\mathcal C_{\Gamma^\gamma}(x)}(x)
[2\Gamma^\gamma(x)]^{-1/2}$.
As in \eqref{eq:six-local-bridge-insertion}, we have
\begin{equation}\label{eq:six-test-box-bridge-insertion}
 I_s(x,y)
 :=\mathbb E_{\widetilde D_s}\!\left[
       |\phi_{\widetilde D_s,y}|e_{\mathcal C_{\widetilde D_s}(x)}(x)
       \mathbf{1}_{\smash[t]{\{x\xleftrightarrow{\mathcal L_{\widetilde D_s}}y\}}}\right]
 =\int\mu_{\widetilde D_s}^{x,y}(d\gamma)\,
                   \mathbb E_{\widetilde D_s}[F_x(\gamma)].
\end{equation}
Expanding the capacity in \eqref{eq:six-amplitude-sum}, we obtain
\begin{equation}\label{eq:six-capacity-amplitude-expansion}
 \mathcal M(s)=\sum_{x,y\in D_s}I_s(x,y).
\end{equation}

Fix \(x\in B(se_1,s/32)\). Independently of the loop soup
\(\mathcal L_{\widetilde D_s}\), let \(Z^x\) be the Brownian motion from \(x\), stopped at
\(\zeta=H_{\widetilde D_s^c}(Z^x)\).
Since
\(\widetilde D_s\subset\widetilde B(x,s/4)\), we obtain by monotonicity and
translation invariance
\begin{equation}\label{eq:six-terminal-weight-comparison}
 \mathbb E F_x(Z^x_{[0,\zeta]})\ge b(s/4).
\end{equation}

For each integer \(t\in[s/64,s/32]\), define the ball
\(\mathrm{B}_t=\{y\in\mathbb Z^6:\|y-x\|_2\le t\}\),
where \(\|\cdot\|_2\) denotes the Euclidean norm. Its equilibrium
measure for the walk killed on leaving \(D_s\) is
\[
 e_{\mathrm{B}_t}^{D_s}(y)
 =\mathbf{1}_{\{y\in \mathrm{B}_t\}}
       \mathbb P_y(H_{D_s^c}<\widehat H_{\mathrm{B}_t}),
\]
where
\(\widehat H_{\mathrm{B}_t}\) is the first positive return time to
\(\mathrm{B}_t\). Let
\(\lambda_t=\sup\{u<\zeta:Z_u^x\in \mathrm{B}_t\}\) be the last visit
of $Z^x$ to \(\mathrm{B}_t\). 
By monotonicity, the last-exit decomposition~\cite[Appendix~A]{Pre23}, and \eqref{eq:six-test-box-bridge-insertion},
we obtain
\begin{equation}\label{eq:six-last-exit-weighted-bound}
 \mathbb E F_x(Z^x_{[0,\zeta]})
 \le\mathbb E F_x(Z^x_{[0,\lambda_t]})
 =\sum_{y\in D_s}e_{\mathrm{B}_t}^{D_s}(y)I_s(x,y).
\end{equation}
By~\cite[Proposition~6.4.2]{LL10}, we have
\begin{equation}\label{eq:six-annular-equilibrium-bound}
 e_{\mathrm{B}_t}^{D_s}(y)
 \lesssim s^{-1}
       \mathbf{1}_{\{t-1\le\|y-x\|_2\le t\}}.
\end{equation}
Summing \eqref{eq:six-last-exit-weighted-bound} over the integers
\(t\in[s/64,s/32]\) and using
\eqref{eq:six-annular-equilibrium-bound} and
\eqref{eq:six-terminal-weight-comparison}, we obtain
\begin{equation}\label{eq:six-annular-bridge-lower}
 \sum_{s/64-1\le\|y-x\|_2\le s/32}I_s(x,y)
 \gtrsim s^2\mathbb E F_x(Z^x_{[0,\zeta]})
 \ge s^2b(s/4).
\end{equation}
Summing \eqref{eq:six-annular-bridge-lower} over \(x\in B(se_1,s/32)\)
and using \eqref{eq:six-capacity-amplitude-expansion} proves
\(\mathcal M(s)\gtrsim s^8b(s/4)\).
By Lemma~\ref{s6:capacity-amplitude-lower}, we then have
\(T(s)\gtrsim s^{-8}\mathcal M(s)^2\gtrsim s^8b(s/4)^2\),
which completes the proof.
\end{proof}

\subsection{Multiscale tests and renormalization}\label{s6:relative-escape}

For a sufficiently large integer \(h\) and large \(R\),
choose the test scales
\begin{equation}\label{eq:six-test-scales}
 s_0=8h\lceil\sqrt R\rceil,
 \qquad s_j=32^j s_0\le R/100,\qquad 0\le j\le J,
\end{equation}
where \(J=\lfloor\log_{32}(R/(100s_0))\rfloor\).
Recall the setup in Lemma~\ref{s6:second-moment}.
At scale \(s_j\), write \(\mathbb P_{a,b}\) for the joint law of the
two independent walks started at \(a,b\), each stopped on leaving
\(B_{4s_j}\), and the independent loop soup \(\mathcal L_{\widetilde D_{s_j}}\).
For \(0\le j\le J\), define the common-intersection probability at scale $s_j$ by
\begin{equation}\label{eq:six-test-probability-definition}
 p_j:=\inf_{a,b\in\partial B_{s_j/4}}
                     \mathbb P_{a,b}(N(s_j)>0).
\end{equation}
Sequential conditioning over these scales yields the following estimate.
\begin{lemma}\label{s6:separated-scale-contraction}
There exists an integer \(h\ge1\) such that, for all sufficiently large integers \(R\),
\begin{equation}\label{eq:six-fixed-scale-closure}
 b(R)\le2b(\sqrt R)
                   \prod_{j=0}^J(1-p_j).
\end{equation}
\end{lemma}

\begin{proof}
Let \(X\) be the jump chain of the Brownian motion \(Z\)
used to define \(\Gamma_R\) in \eqref{eq:six-enhanced-occupation},
stopped on leaving \(\widetilde B_R\).
Let \(Y\) be another simple random walk from zero, independent of \(Z\)
and the loop soup $\mathcal L_{\widetilde B_R}$.
For each of \(X,Y,Z\), let \(\alpha_j\) be its first hitting time of
\(\partial B_{s_j/4}\), and let \(\beta_j\) be its first exit time from
\(B_{4s_j}\) after \(\alpha_j\).
Define the common-intersection event at scale $s_j$ by
\begin{equation}\label{eq:six-common-cluster-test}
 \mathcal A_j:=\left\{\exists\,\mathcal C\in\mathfrak C_{s_j}:
 X([\alpha_j,\beta_j])\cap\mathcal C\ne\varnothing,\quad
 Y([\alpha_j,\beta_j])\cap\mathcal C\ne\varnothing\right\}.
\end{equation}
Now choose \(h\) large enough that for every \(y\notin B_{h\sqrt R}\), we have
\(\mathbb P_y(H_{B_{\sqrt R}}<\infty)\le\tfrac12\).
Define the non-intersection event
$
 F_{\sqrt R}=\{Y_k\notin\mathcal K_{\sqrt R}\text{ for }1\le k\le H_{B_{h\sqrt R}^c}\}
$. Then, by monotonicity and the Markov property, we have 
\begin{equation}\label{eq:six-inner-nonintersection-weight}
 e_{\mathcal K_{\sqrt R}}(0)\ge\tfrac12\mathbb P(F_{\sqrt R}\mid\mathcal K_{\sqrt R}),
 \qquad \mathbb E[(2\Gamma_{\sqrt R}(0))^{-1/2}\mathbf1_{F_{\sqrt R}}]\le2b(\sqrt R).
\end{equation}
Moreover, \((2\Gamma_R(0))^{-1/2}\le (2\Gamma_{\sqrt R}(0))^{-1/2}\), and by definition,
\begin{equation}\label{eq:six-multiscale-nonintersection-bound}
 b(R)\le\mathbb E\!\left[
              (2\Gamma_{\sqrt R}(0))^{-1/2}\mathbf1_{F_{\sqrt R}}\prod_{j=0}^J\mathbf1_{\mathcal A_j^c}
                         \right].
\end{equation}
Let \(\mathcal H_j\) record the inner loop soup, the loop soups used in earlier
tests, and the two paths up to their entrance times,
\begin{equation}\label{eq:six-test-history}
 \mathcal H_j:=\sigma\!\left(
 \mathcal L_{\widetilde B_{\sqrt R}},
 (\mathcal L_{\widetilde D_{s_i}})_{0\le i<j},
 Z([0,\alpha_j]),
 Y([0,\alpha_j])
 \right).
\end{equation}
The test boxes and the inner box are pairwise disjoint. Moreover,
\(4s_i+1<s_j/4\) for \(i<j\), so all earlier tests and the inner weighted
event are \(\mathcal H_j\)-measurable. The loop soup in
\(\widetilde D_{s_j}\) is independent of \(\mathcal H_j\).
By the strong Markov property, we therefore have
\begin{equation}\label{eq:six-conditional-test-probability}
 \mathbb P(\mathcal A_j\mid\mathcal H_j)
 =\mathbb P_{X_{\alpha_j},Y_{\alpha_j}}(N(s_j)>0)
 \ge p_j.
\end{equation}
Conditioning successively in \eqref{eq:six-multiscale-nonintersection-bound}
using \eqref{eq:six-conditional-test-probability}, and then applying
\eqref{eq:six-inner-nonintersection-weight}, we obtain
\[
 b(R)\le\left(\prod_{j=0}^J(1-p_j)\right)
       \mathbb E[(2\Gamma_{\sqrt R}(0))^{-1/2}\mathbf1_{F_{\sqrt R}}]
 \le2b(\sqrt R)\prod_{j=0}^J(1-p_j).
\]
This completes the proof.
\end{proof}

We now combine all the previous estimates to prove the escape estimate by strong induction.

\begin{proof}[Proof of Proposition~\ref{s6:escape-bound}]

Fix a sufficiently large integer \(R_0\ge16\) and let
\(M\ge\log R_0\) be a constant to be chosen below.
We prove \(b(R)^2\le M/\log R\) by strong induction on integers \(R\ge2\).
For \(2\le R<R_0\), this follows from \(b\le1\) in
\eqref{eq:six-b-uniform-bound}.
Let \(R\ge R_0\), assume the claim for all smaller integers,
and suppose that \(b(R)^2>M/\log R\).

Use the test scales \(s_j\) in \eqref{eq:six-test-scales}
and the probabilities \(p_j\) in \eqref{eq:six-test-probability-definition}.
For all sufficiently large \(R\) and
\(0\le j\le J\), we have \(s_j\ge\sqrt R\), so both
\(\sqrt{s_j}\) and \(\sqrt R\) are at least \(R^{1/4}\).
Also \(s_j/4\le R\). By monotonicity, we therefore have
\(b(\sqrt{s_j}),\ b(\sqrt R)\le b(R^{1/4})\), \(b(s_j/4)\ge b(R)\), and \(b(R)\le b(R^{1/4})\).
By
Lemma~\ref{s6:annular-capacity-lower}, we have
\[
 \mathcal M(s_j)\ge c s_j^8b(s_j/4)\ge c s_j^8b(R),
 \qquad T(s_j)\ge c s_j^8b(s_j/4)^2\ge c s_j^8b(R)^2.
\]
Lemma~\ref{s6:annealed-triangle} gives
\[
 s_j^{-8}Q(s_j)\le C b(\sqrt{s_j})^2+C s_j^{-1/2}
             \le C(b(R^{1/4})^2+R^{-1/4}).
\]
By \eqref{eq:six-count-moments} in Lemma~\ref{s6:second-moment},
we thus have, uniformly in \(a,b\in\partial B_{s_j/4}\),
\[
 \begin{aligned}
 \mathbb E_{a,b}N(s_j)&\asymp s_j^{-8}T(s_j)\gtrsim b(R)^2,\\
 \mathbb E_{a,b}[N(s_j)(N(s_j)-1)]&\lesssim s_j^{-8}Q(s_j)
                                   \lesssim b(R^{1/4})^2+R^{-1/4}.
 \end{aligned}
\]
Applying the Paley--Zygmund inequality with
\eqref{eq:six-test-probability-definition} and the preceding moment bounds,
we obtain, for suitable constants \(c_1,C_1>0\),
\begin{equation}\label{eq:six-test-success-lower}
 \begin{aligned}
 p_j&\ge\inf_{a,b\in\partial B_{s_j/4}}
       \frac{(\mathbb E_{a,b}N(s_j))^2}
            {\mathbb E_{a,b}N(s_j)+\mathbb E_{a,b}[N(s_j)(N(s_j)-1)]}\\
    &\ge\frac{c_1^2b(R)^4}
              {c_1b(R)^2+C_1\bigl(b(R^{1/4})^2+R^{-1/4}\bigr)}
     \ge c\frac{b(R)^4}{b(R^{1/4})^2+R^{-1/4}},
 \end{aligned}
\end{equation}
where the last inequality uses \(b(R)^2\le b(R^{1/4})^2+R^{-1/4}\).
Using the induction hypothesis at \(R^{1/4}\) and
\(R^{-1/4}\le(\log R)^{-1}\), we obtain
\(b(R^{1/4})^2+R^{-1/4}\lesssim M/\log R\).
Combining \(b(R)^2>M/\log R\) with
\eqref{eq:six-test-success-lower}, we obtain
\(p_j\ge cM/\log R\) for every \(j\).
By \eqref{eq:six-test-scales}, \(J+1\gtrsim\log R\), and hence
\begin{equation}\label{eq:six-test-product-decay}
 \prod_{j=0}^J(1-p_j)
 \le\exp\left\{-\sum_{j=0}^Jp_j\right\}
 \le e^{-cM}.
\end{equation}
Combining \eqref{eq:six-fixed-scale-closure} with
\eqref{eq:six-test-product-decay} and the induction hypothesis
at \(\sqrt R\), we obtain
\[
 b(R)\le2b(\sqrt R)\prod_{j=0}^J(1-p_j)
   \le C\sqrt{\frac{M}{\log R}}\,e^{-cM}
   <\sqrt{\frac{M}{\log R}},
\]
where the last inequality holds for \(M\) sufficiently large,
independently of \(R\). This contradicts \(b(R)^2>M/\log R\),
completes the induction, and proves
\(b(R)\lesssim(\log R)^{-1/2}\).
\end{proof}

The one-arm upper bound now gives matching lower bounds for the weighted
mean capacity \(M(r)\) and the weighted escape probability \(b(r)\).

\begin{proposition}\label{s6:capacity-escape-two-sided}
For all sufficiently large \(r\),
\begin{equation}\label{eq:six-capacity-escape-two-sided}
 M(r)\asymp r^2(\log r)^{-\frac12},\qquad
 b(r)\asymp(\log r)^{-\frac12}.
\end{equation}
\end{proposition}

\begin{proof}
The upper bounds follow from~\eqref{eq:Mr} and
Proposition~\ref{s6:escape-bound}.
For the lower bound on \(M(r)\), Lemma~\ref{lem:capacity-arm-transfer}
and the one-arm upper bound in Theorem~\ref{thm:main} give
\(M(r)\gtrsim\theta_6(r)^{-1}
\gtrsim r^2(\log r)^{-\frac12}\).
Applying~\eqref{eq:six-mean-capacity-bound} at radius \(r^2\), we obtain
\((\log r)^{-\frac12}
\lesssim r^{-4}M(r^2)\lesssim b(r)+r^{-1}\).
This proves the
lower bound on \(b(r)\).
\end{proof}

\section{Macroscopic cluster counts in dimension six}
\label{sec:large-clusters}

In this section, we prove Theorem~\ref{thm:six-large-clusters}.
Throughout this section, 
\(\mathbb P_r\) and \(\mathbb E_r\) denote probability and
expectation for \(\mathcal L_{\widetilde B_r}\).
We use the vertex clusters \(\mathcal C^r(x)\), their family
\(\mathfrak C^r\), and the thresholds \(v_r^-\), \(v_r^+\)
introduced in Section~\ref{sec:main-results}.

\subsection{Volume moments}

We first establish volume-moment estimates from the two-point function
and the mean-capacity bound.

\begin{lemma}\label{lem:six-cluster-volume-moments}
For all sufficiently large \(r\), we have
\begin{equation}\label{eq:six-cluster-volume-moments}
 \sup_{z\in B_r}\mathbb E_r|\mathcal C^r(z)|
 \lesssim r^2.
\end{equation}
Moreover,
\begin{equation}\label{eq:six-root-volume-second-moment}
 \sup_{z\in B_r}\mathbb E_r|\mathcal C^r(z)|^2
 \lesssim r^6(\log r)^{-\frac12}.
\end{equation}
\end{lemma}

\begin{proof}
For \eqref{eq:six-cluster-volume-moments}, the two-point function
\eqref{eq:two-point} and the Green function estimate
\eqref{eq:green-estimates} give, for every \(z\in B_r\),
\[
 \mathbb E_r|\mathcal C^r(z)|
 =\sum_{y\in B_r}
       \mathbb P_r(z\xleftrightarrow{\mathcal L_{\widetilde B_r}}y)
 \lesssim\sum_{y\in B_r}G(z,y)\lesssim r^2.
\]
We next prove \eqref{eq:six-root-volume-second-moment}.
By \eqref{eq:capacity-fkg-comparison} and \eqref{eq:Mr}, we have
\begin{equation}\label{eq:six-unweighted-mean-capacity}
 \mathbb E\operatorname{cap}(\mathcal C(0)\cap B_s)
 \le \frac{M(s)}{\mathbb E|\phi_0|}
 \lesssim s^2(\log s)^{-\frac12}.
\end{equation}
Since \(\mathcal C^r(z)\subseteq\mathcal C(z)\cap B(z,3r)\)
for every \(z\in B_r\), the above also yields
\begin{equation}\label{eq:six-box-mean-capacity}
 \mathbb E_r\operatorname{cap}(\mathcal C^r(z) )
\lesssim r^2(\log r)^{-\frac12}.
\end{equation}

For \(x,y,z\in B_r\), let \(E_{x,y,z}^r\) be the event that
the three vertices belong to the same cluster of
\(\mathcal L_{\widetilde B_r}\).
Applying the FKG
inequality twice gives
\[
 \mathbb E_r\!\left[
      |\phi_{\widetilde B_r,x}\phi_{\widetilde B_r,y}|;
      E_{x,y,z}^r\right]
 \ge
 \mathbb E_r|\phi_{\widetilde B_r,x}|\,
 \mathbb E_r|\phi_{\widetilde B_r,y}|\,
 \mathbb P_r(E_{x,y,z}^r)
 \gtrsim\mathbb P_r(E_{x,y,z}^r).
\]
Applying \eqref{eq:bridge} in \(\widetilde B_r\), with \(F\)
the indicator that \(x\) and \(z\) lie in the same connected component of the set where the occupation field is positive, we obtain
\[
 \mathbb E_r\!\left[
      |\phi_{\widetilde B_r,x}\phi_{\widetilde B_r,y}|;
      E_{x,y,z}^r\right]
 =\mathbb E_r\!\left[
      \mu_{\widetilde B_r}^{x,y}
       \bigl(\gamma\cap\widetilde{\mathcal C}^r(z)\ne\varnothing\bigr)
      \right].
\]
Combining the preceding two estimates gives
\begin{equation}\label{eq:six-three-point-bridge}
 \mathbb P_r(E_{x,y,z}^r)
 \lesssim\mathbb E_r\!\left[
      \mu_{\widetilde B_r}^{x,y}
       (\gamma\text{ visits }\mathcal C^r(z))\right].
\end{equation}

For a nonempty finite discrete set \(A\), by the strong Markov property at the first visit to \(A\),
and the last-exit formula~\eqref{eq:capacity-hitting-estimates}, we obtain
\begin{equation}\label{eq:six-bridge-hit-capacity}
 \begin{aligned}
 \sum_{x,y\in B_r}\mu_{\widetilde B_r}^{x,y}(\gamma\text{ visits }A)
 &\le\sum_{x,y\in B_r}\mu^{x,y}(\gamma\text{ visits }A)\\
 &=\sum_{x\in B_r}\mathbb E_x\!\left[
       \sum_{y\in B_r}G(X_{H_A},y);H_A<\infty\right]\\
 &\le\biggl(\sup_{u\in A}\sum_{y\in B_r}G(u,y)\biggr)^2
       \operatorname{cap}(A)
 \lesssim r^4\operatorname{cap}(A).
 \end{aligned}
\end{equation}
Here \(X\) is the continuous-time random walk.
The last bound follows from \eqref{eq:green-estimates}.

Summing \eqref{eq:six-three-point-bridge} over \(x,y\in B_r\)
and using \eqref{eq:six-bridge-hit-capacity} and
\eqref{eq:six-box-mean-capacity} yields
\[
 \mathbb E_r|\mathcal C^r(z)|^2
 =\sum_{x,y\in B_r}\mathbb P_r(E_{x,y,z}^r)
 \lesssim r^4\mathbb E_r\operatorname{cap}(\mathcal C^r(z))
 \lesssim r^6(\log r)^{-\frac12}.
\]
This proves \eqref{eq:six-root-volume-second-moment} and completes the proof.
\end{proof}

Recall that \(v_r^-=ar^4/\sqrt{\log r}\).
For \(x\in B_{r/2}\), define
\[
 V_r(x)=\bigl|\mathcal C^r(x)\cap
              (B(x,r/4)\setminus B(x,r/8))\bigr|,
 \qquad
 \mathcal G_r=\{x\in B_{r/2}:V_r(x)\ge v_r^-\}.
\]
We call the vertices in \(\mathcal G_r\) good. Then, we combine the preceding
moment estimates with the one-arm upper bound to estimate the number of good vertices.

\begin{lemma}\label{lem:six-good-vertices}
For \(a\) sufficiently small and all sufficiently large \(r\), we have
\begin{equation}\label{eq:six-good-vertex-probability}
 \mathbb P_r(x\in\mathcal G_r)
 \gtrsim r^{-2}(\log r)^{\frac12},\qquad x\in B_{r/2}.
\end{equation}
Moreover, there is a constant \(c_1>0\) such that
\begin{equation}\label{eq:six-good-vertex-concentration}
 \mathbb P_r\!\left(|\mathcal G_r|
          \ge c_1r^4(\log r)^{\frac12}\right)
 \ge1-O\bigl((\log r)^{-1}\bigr).
\end{equation}
\end{lemma}

\begin{proof}
We first prove \eqref{eq:six-good-vertex-probability}.
Using~\eqref{eq:two-point} and~\eqref{eq:green-estimates}, 
we obtain
\(\mathbb E_r V_r(x)\asymp r^2\).
Since the annulus lies in \(B_r\),
\eqref{eq:six-root-volume-second-moment} also gives
\(\mathbb E_r V_r(x)^2\lesssim r^6(\log r)^{-\frac12}\).
The event \(\{V_r(x)>0\}\) implies a one-arm event from
\(x\) to distance \(\lfloor r/8\rfloor\). 
By the one-arm upper bound
\eqref{eq:main-six}, we also have
\[
 \mathbb E_r[V_r(x);V_r(x)<v_r^-]
 \le v_r^-\mathbb P_r(V_r(x)>0)\lesssim ar^2.
\]
Taking \(a\) sufficiently small, we obtain
\(\mathbb E_r[V_r(x);V_r(x)\ge v_r^-]\gtrsim r^2\).
The Cauchy--Schwarz inequality then gives
\[
 r^4\lesssim
 \mathbb E_r[V_r(x);V_r(x)\ge v_r^-]^2
 \le\mathbb E_r V_r(x)^2\,\mathbb P_r(x\in\mathcal G_r),
\]
which proves \eqref{eq:six-good-vertex-probability}. In particular,
\(\mathbb E_r|\mathcal G_r|\gtrsim r^4(\log r)^{\frac12}\).

We next prove \eqref{eq:six-good-vertex-concentration} by bounding
the variance of \(|\mathcal G_r|\). Splitting according to whether \(x\) and \(y\) lie in the same cluster,
and applying the BK inequality~\cite{vdB96} when their clusters are
distinct, we obtain
\begin{equation}\label{eq:six-good-vertex-joint-probability}
 \begin{aligned}
 \mathbb P_r(x,y\in\mathcal G_r)
 &\le \mathbb P_r(x\xleftrightarrow{\mathcal L_{\widetilde B_r}}y)
     +\mathbb P_r(x,y\in\mathcal G_r,\ y\notin\mathcal C^r(x))\\
 &\le \mathbb P_r(x\xleftrightarrow{\mathcal L_{\widetilde B_r}}y)
     +\mathbb P_r(x\in\mathcal G_r)\mathbb P_r(y\in\mathcal G_r).
 \end{aligned}
\end{equation}
Summing this covariance bound
and using
\eqref{eq:six-cluster-volume-moments}, we obtain
\[
 \operatorname{Var}_r(|\mathcal G_r|)
 \le\sum_{x,y\in B_{r/2}}
          \mathbb P_r(x\xleftrightarrow{\mathcal L_{\widetilde B_r}}y)
 \le\sum_{x\in B_{r/2}}\mathbb E_r|\mathcal C^r(x)|
 \lesssim r^8.
\]
Finally, Chebyshev's inequality gives
\[
 \mathbb P_r\!\left(|\mathcal G_r|<\tfrac12\mathbb E_r|\mathcal G_r|\right)
 \le\frac{4\operatorname{Var}_r(|\mathcal G_r|)}
              {(\mathbb E_r|\mathcal G_r|)^2}
 \lesssim\frac1{\log r},
\]
which proves \eqref{eq:six-good-vertex-concentration}.
\end{proof}

\subsection{Volume truncation and the cluster count}

\begin{proof}[Proof of Theorem~\ref{thm:six-large-clusters}]
Fix \(a\) as in Lemma~\ref{lem:six-good-vertices}. We will choose
\(A\ge1\) sufficiently large, independently of \(r\). Define
\[
 \begin{aligned}
 S_r^+&=\sum_{\mathcal C\in\mathfrak C^r}
       |\mathcal C|\mathbf1_{\{|\mathcal C|>v_r^+\}}
     =\sum_{x\in B_r}\mathbf1_{\{|\mathcal C^r(x)|>v_r^+\}},\\
 S_r^-&=\sum_{\mathcal C\in\mathfrak C^r}
       |\mathcal C|\mathbf1_{\{|\mathcal C|\ge v_r^-\}}
     =\sum_{x\in B_r}\mathbf1_{\{|\mathcal C^r(x)|\ge v_r^-\}}.
 \end{aligned}
\]
Thus \(S_r^+\) counts vertices in clusters above the upper volume
threshold, and \(S_r^-\) counts vertices in clusters at or above the
lower threshold.

We first bound both variances as in 
Lemma~\ref{lem:six-good-vertices}. For either choice
\(E_x=\{|\mathcal C^r(x)|>v_r^+\}\) or
\(E_x=\{|\mathcal C^r(x)|\ge v_r^-\}\), the events are increasing.
Applying the BK inequality again, as in \eqref{eq:six-good-vertex-joint-probability}, we obtain
\[
 \mathbb P_r(E_x\cap E_y)
 \le \mathbb P_r(x\xleftrightarrow{\mathcal L_{\widetilde B_r}}y)
       +\mathbb P_r(E_x)\mathbb P_r(E_y).
\]
Summing and using \eqref{eq:six-cluster-volume-moments} again gives
\begin{equation}\label{eq:six-volume-mass-variance}
 \operatorname{Var}_r(S_r^+)+\operatorname{Var}_r(S_r^-)
 \le 2\sum_{x,y\in B_r}
       \mathbb P_r(x\xleftrightarrow{\mathcal L_{\widetilde B_r}}y)
 \lesssim r^8.
\end{equation}

For the lower bound, applying Markov's inequality at each vertex
and summing, we obtain from \eqref{eq:six-cluster-volume-moments} that
\[
 \mathbb E_r S_r^+
 \le (v_r^+)^{-1}\sum_{x\in B_r}\mathbb E_r|\mathcal C^r(x)|
 \lesssim A^{-1}r^4(\log r)^{\frac12}.
\]
Fix \(A\) large enough that
\(\mathbb E_r S_r^+\le(c_1/4)r^4(\log r)^{1/2}\), where \(c_1\)
is the constant in \eqref{eq:six-good-vertex-concentration}.
By Chebyshev's inequality and \eqref{eq:six-volume-mass-variance},
\[
 \mathbb P_r\!\left(S_r^+>\frac{c_1}{2}r^4(\log r)^{\frac12}\right)
 \le\frac{16\operatorname{Var}_r(S_r^+)}{c_1^2r^8\log r}
 \lesssim\frac1{\log r}.
\]
Combining this with \eqref{eq:six-good-vertex-concentration} by a
union bound, we obtain that, with probability at least
\(1-O((\log r)^{-1})\), both
\[
 |\mathcal G_r|\ge c_1r^4(\log r)^{\frac12},
 \qquad
 S_r^+\le\frac{c_1}{2}r^4(\log r)^{\frac12}
\]
hold. On this event, clusters with volume greater than \(v_r^+\)
contain at most \(S_r^+\) good vertices, since \(S_r^+\) counts
all vertices in these clusters. Therefore,
\begin{equation}\label{eq:six-remaining-good-vertices}
 \#\{x\in\mathcal G_r:|\mathcal C^r(x)|\le v_r^+\}
 \ge |\mathcal G_r|-S_r^+
 \ge\frac{c_1}{2}r^4(\log r)^{\frac12}.
\end{equation}
Note that a good vertex \(x\) is connected to at least
\(v_r^-\) vertices in \(B(x,r/4)\setminus B(x,r/8)\subset B_r\),
and hence \(\operatorname{diam}(\mathcal C^r(x))\ge r/8\).
Each cluster containing a good vertex counted on the left-hand
side of \eqref{eq:six-remaining-good-vertices} is counted by \(N_r\).
Since each such cluster contains at most \(v_r^+\) good vertices,
by~\eqref{eq:six-remaining-good-vertices}, we have
\begin{equation}\label{eq:six-cluster-count-lower}
 N_r\ge\frac{(c_1/2)r^4(\log r)^{\frac12}}{v_r^+}
       =\frac{c_1}{2A}\log r
\end{equation}
on the same event.

For the upper bound, Markov's inequality similarly gives
\[
 \mathbb E_r S_r^-
 \le (v_r^-)^{-1}\sum_{x\in B_r}\mathbb E_r|\mathcal C^r(x)|
 \le C_2r^4(\log r)^{\frac12}
\]
for a fixed constant \(C_2>0\). Applying Chebyshev's
inequality and \eqref{eq:six-volume-mass-variance}, we obtain
\[
 \mathbb P_r\!\left(S_r^->2C_2r^4(\log r)^{\frac12}\right)
 \le\frac{\operatorname{Var}_r(S_r^-)}{C_2^2r^8\log r}
 \lesssim\frac1{\log r}.
\]
Every cluster counted by \(N_r\) contributes at least \(v_r^-\)
vertices to \(S_r^-\). Hence \(N_rv_r^-\le S_r^-\).
On the event \(\{S_r^-\le2C_2r^4(\log r)^{1/2}\}\), we therefore have
\begin{equation}\label{eq:six-cluster-count-upper}
 N_r\le\frac{2C_2r^4(\log r)^{\frac12}}{v_r^-}
      =\frac{2C_2}{a}\log r.
\end{equation}
By a union bound, \eqref{eq:six-cluster-count-lower} and
\eqref{eq:six-cluster-count-upper} hold simultaneously with probability
at least \(1-O((\log r)^{-1})\).
This proves \eqref{eq:six-large-cluster-count}.
\end{proof}

\bibliographystyle{alpha}
\begingroup
\small
\bibliography{references}
\endgroup
\end{document}